%% file: laplacian_eigenmaps_regression.tex
\documentclass{article}

\usepackage[utf8]{inputenc} 
\usepackage[T1]{fontenc}    
\usepackage{booktabs}       
\usepackage{nicefrac}       
\usepackage{microtype}      
\usepackage{times}          

\usepackage[round]{natbib}
\usepackage{amssymb,amsmath,amsthm,bbm}
\usepackage[margin=1in]{geometry}
\usepackage{verbatim,float,url,dsfont}
\usepackage{graphicx,subfigure,psfrag}
\usepackage{algorithm,algorithmic}
\usepackage{mathtools,enumitem}
\usepackage[colorlinks=true,citecolor=blue,urlcolor=blue,linkcolor=blue]{hyperref}
\usepackage{multirow}

\newtheorem{theorem}{Theorem}
\newtheorem{lemma}{Lemma}

\newtheorem{proposition}{Proposition}
\theoremstyle{definition}

\newtheorem{definition}{Definition}

\newtheorem*{assumption*}{\assumptionnumber}
\providecommand{\assumptionnumber}{}


\makeatletter
\newcommand*\rel@kern[1]{\kern#1\dimexpr\macc@kerna}
\newcommand*\widebar[1]{%
	\begingroup
	\def\mathaccent##1##2{%
		\rel@kern{0.8}%
		\overline{\rel@kern{-0.8}\macc@nucleus\rel@kern{0.2}}%
		\rel@kern{-0.2}%
	}%
	\macc@depth\@ne
	\let\math@bgroup\@empty \let\math@egroup\macc@set@skewchar
	\mathsurround\z@ \frozen@everymath{\mathgroup\macc@group\relax}%
	\macc@set@skewchar\relax
	\let\mathaccentV\macc@nested@a
	\macc@nested@a\relax111{#1}%
	\endgroup
}
\makeatother

\newcommand{\argmin}{\mathop{\mathrm{argmin}}}

\newcommand{\minimize}{\mathop{\mathrm{minimize}}}

\newcommand{\Ebb}{\mathbb{E}}
\newcommand{\Pbb}{\mathbb{P}}
\newcommand{\dotp}[2]{\langle #1, #2 \rangle}
\newcommand{\wt}[1]{\widetilde{#1}}
\newcommand{\wh}[1]{\widehat{#1}}
\newcommand{\mc}[1]{\mathcal{#1}}
\newcommand{\Reals}{\mathbb{R}} 
\newcommand{\Rd}{\Reals^d}
\newcommand{\wb}[1]{\widebar{#1}}
\newcommand{\floor}[1]{\left\lfloor #1 \right\rfloor}
\newcommand{\Var}{\mathrm{Var}}
\newcommand{\Cov}{\mathrm{Cov}}
\newcommand{\1}{\mathbf{1}}
\newcommand{\bj}{{\bf j}}
\newcommand{\restr}[2]{\ensuremath{\left.#1\right|_{#2}}}

\DeclareFontFamily{U}{mathx}{\hyphenchar\font45}
\DeclareFontShape{U}{mathx}{m}{n}{<-> mathx10}{}
\DeclareSymbolFont{mathx}{U}{mathx}{m}{n}
\DeclareMathAccent{\wc}{0}{mathx}{"71}

\begin{document}
	
\begin{center} {\Large{\bf{Minimax Optimal Regression over Sobolev Spaces \\
				\vspace{.2cm}
				via Laplacian Eigenmaps on Neighborhood Graphs}}}
	
	\vspace*{.3cm}
	
	{\large{
			\begin{center}
				Alden Green~~~~~ Sivaraman Balakrishnan~~~~~ Ryan J. Tibshirani\\
				\vspace{.2cm}
			\end{center}

			\begin{tabular}{c}
				Department of Statistics and Data Science \\
				Carnegie Mellon University
			\end{tabular}
			
			\vspace*{.2in}
			
			\begin{tabular}{c}
				\texttt{\{ajgreen,siva,ryantibs\}@stat.cmu.edu}
			\end{tabular}
	}}
	
	\vspace*{.2in}
	
	\today
	\vspace*{.2in}
\end{center}

\begin{abstract}
	In this paper we study the statistical properties of Principal Components Regression with Laplacian Eigenmaps (PCR-LE), a method for nonparametric regression based on Laplacian Eigenmaps (LE). PCR-LE works by projecting a vector of observed responses ${\bf Y} = (Y_1,\ldots,Y_n)$ onto a subspace spanned by certain eigenvectors of a neighborhood graph Laplacian. We show that PCR-LE achieves minimax rates of convergence for random design regression over Sobolev spaces. Under sufficient smoothness conditions on the design density $p$, PCR-LE achieves the optimal rates for both estimation (where the optimal rate in squared $L^2$ norm is known to be $n^{-2s/(2s + d)}$) and goodness-of-fit testing ($n^{-4s/(4s + d)}$). We also show that PCR-LE is \emph{manifold adaptive}: that is, we consider the situation where the design is supported on a manifold of small intrinsic dimension $m$, and give upper bounds establishing that PCR-LE achieves the faster minimax estimation ($n^{-2s/(2s + m)}$) and testing ($n^{-4s/(4s + m)}$) rates of convergence. Interestingly, these rates are almost always much faster than the known rates of convergence of graph Laplacian eigenvectors to their population-level limits; in other words, for this problem regression with estimated features appears to be much easier, statistically speaking, than estimating the features itself. We support these theoretical results with empirical evidence.
\end{abstract}

\input{main_text.tex}


\input{laplacian_eigenmaps_regression.bbl}\newpage

\appendix
\input{appendix.tex}

\end{document}

%% file: main_text.tex
\section{Introduction}
\label{sec:introduction}

Laplacian Eigenmaps (LE)~\citep{belkin03a} is a method for nonlinear dimensionality reduction and data representation. Given data points $\{X_1,\ldots,X_n\} \subset \Reals^d$, LE maps each $X_i$ to a vector $(v_{1,i},\ldots,v_{K,i})$ according to the following steps.
\begin{enumerate}
	\item First, LE forms a \emph{neighborhood graph} $G = (V,W)$ over the points $\{X_1,\ldots,X_n\}$. The graph $G$ is an undirected, weighted graph, with vertices $V = \{X_1,\ldots,X_n\}$, and weighted edges $W_{ij}$ which correspond to the proximity between points $X_i$ and $X_j$.
	\item Next, LE forms an (unweighted) \emph{graph Laplacian} matrix $L \in \Reals^{n \times n}$, a symmetric and diagonally dominant matrix with diagonal elements $L_{ii} = \sum_{j = 1}^{n} W_{ij}$, and off-diagonal elements $L_{ij} = -W_{ij}$. 
	\item Finally, LE takes the eigendecomposition $L = \sum_{k = 1}^{n} \lambda_k v_k v_k^{\top}$, and outputs the vectors $(v_{1,i},\ldots,v_{K,i}) \in \Reals^K$ for each $i = 1,\ldots,n$.
\end{enumerate} 
A natural way to use LE is by taking the collection of vectors $\{(v_{1,i},\ldots,v_{K,i})\}_{i = 1}^{n}$ to be features in a downstream regression algorithm. In this paper, we study a simple method along these lines: Principal Components Regression with Laplacian-Eigenmaps (PCR-LE), a method for nonparametric regression which operates by running ordinary least squares (OLS) using the features output by LE. Given pairs of design points and responses $(X_1,Y_1),\ldots, (X_n,Y_n)$, PCR-LE computes an estimate $\wh{f} \in \Reals^n$,
\begin{equation}
\label{eqn:pcr-le}
\wh{f} := \argmin_{f \in \mathrm{span}\{v_1,\ldots,v_K\}} \|{\bf Y} - f\|_2^2
\end{equation}
where ${\bf Y} = (Y_1,\ldots,Y_n) \in \Reals^n$ is the vector of responses and $\|\cdot\|_2$ denotes the usual Euclidean norm in $\Reals^n$. (For a formal definition of LE and PCR-LE, see Section~\ref{subsec:laplacian_eigenmaps}.)

LE has been practically very successful, and by now has been used for various statistical tasks such as spectral clustering, manifold learning, level-set estimation, semi-supervised learning, etc. At this point there exists a rich literature \citep{koltchinskii2000,belkin07,vonluxburg2008,burago2014,shi2015,singer2017,garciatrillos18,trillos2019, calder2019, cheng2021,dunson2021} explaining this practical success from a theoretical perspective. Loosely speaking, these works model the design points as being independent samples from a distribution $P$ with density $p$, and show that in this case the eigenvectors of the graph Laplacian $L$ are good empirical approximations of population-level objects. These population-level objects are eigenfunctions $\psi_k$---meaning solutions, along with eigenvalues $\rho_k$, to the equation $\Delta_P \psi_k = \rho_k \psi_k$--- of a density-weighted Laplacian operator defined via:
\begin{equation}
\label{eqn:density-weighted-laplace}
\Delta_Pf := -\frac{1}{p}~ \mathrm{div}(p^2 \nabla f).
\end{equation}  
(Here $\mathrm{div}$ stands for the divergence operator, and $\nabla$ for the gradient. See~\eqref{eqn:laplace_beltrami_eigenproblem} for the formal definition of eigenpairs $(\rho_k,\psi_k)$.) These eigenfunctions in turn characterize various interesting structural aspects of $p$, such as the location and number of high- and low-density regions, the shape and intrinsic dimension of its support, and so forth.

These aforementioned works justify LE as method for data representation, by establishing that each feature vector $(v_{1,i},\ldots,v_{K,i})$ serves an empirical approximation to an idealized representation $(\psi_1(X_i),\ldots,\psi_K(X_i))$. They also provide quantitative guarantees for the accuracy with which LE approximates this ideal representation. However, this theory does not focus on the statistical properties of PCR-LE for classical regression problems such as estimation and testing. That is the major question we address in this paper. We adopt the usual model of nonparametric regression with random design, where one observes independent pairs $(X_1,Y_1),\ldots,(X_n,Y_n)$ of design points and responses. We assume the design points $\{X_1,\ldots,X_n\}$ are sampled from an unknown distribution $P$ supported on $\mc{X} \subseteq \Rd$, and the responses follow a signal plus Gaussian noise model,
\begin{equation}
\label{eqn:model}
Y_i = f_0(X_i) + w_i, \quad w_i \sim N(0,1),
\end{equation}
with noise variables $w_i$  independent of design points $X_i$. The task is to learn the regression function $f_0$, which is unknown but assumed to belong to a Sobolev space $H^s(\mc{X})$. We consider two settings: one where $\mc{X}$ is a full-dimensional domain, and the other where $\mc{X}$ is a low-dimensional submanifold of $\Rd$. In each setting, we derive upper bounds which imply that the PCR-LE estimate $\wh{f}$, and a test using the statistic $T = \|\wh{f}\|_2^2$, are statistically optimal methods for two classical problems in nonparametric regression: estimation and goodness-of-fit testing.  


\paragraph{Sobolev spaces and spectral series regression.}
To analyze PCR-LE, we work in a classical situation where the regression function is assumed to belong to a (Hilbert-)Sobolev space. For an open domain $\mc{X} \subseteq \Rd$, the Sobolev space $H^s(\mc{X})$ consists of all functions $f \in L^2(\mc{X})$ which are $s$-times weakly differentiable, with all order-$s$ partial derivatives $D^{\alpha}f \in L^2(\mc{X})$. We study regression over Sobolev spaces in part because, generally speaking, the minimax rates are well-understood; as mentioned before, when the domain $\mc{X}$ is full-dimensional they are $n^{-2s/(2s + d)}$ for estimation, and $n^{-4s/(4s + d)}$ for testing. For this reason, regression over Sobolev spaces is a good setting in which to see whether PCR-LE measures up to more standard minimax optimal approaches, which have strong theoretical guarantees but are less often used in practice. We give a more specific comparison between PCR-LE and some of these more classical methods in Section~\ref{sec:discussion}.

We also view PCR-LE as being particularly well-suited for regression over Sobolev spaces due to their close connection with \emph{spectral series regression}. Spectral series regression computes generalized empirical Fourier coefficients $\wt{a}_k := \frac{1}{n}\sum_{i = 1}^{n} Y_i \psi_k(X_i)$, and truncates to the $K$-lowest frequency eigenfunctions of $\Delta_P$, producing the estimate
\begin{equation}
\label{eqn:population-level_spectral_series}
\wt{f}(x) = \sum_{k = 1}^{K} \wt{a}_k \psi_k(x).
\end{equation} 
Spectral series regression is intrinsically linked with Sobolev spaces. That is because under appropriate boundary conditions, a ball in the order-$s$ Sobolev space consists of functions $f = \sum_{k} a_k \psi_k \in L^2(\mc{X})$ for which the generalized Fourier coefficients $\{a_k\}_{k = 1}^{\infty}$ satisfy the decay condition $\sum_{k} a_k^2 \rho_k^s \leq C$ (See Section~\ref{subsec:spectral_projection} for more details). This decay condition justifies the truncated series estimator \eqref{eqn:population-level_spectral_series}, since it means the truncation will incur only a limited amount of bias for any $f_0 \in H^s(\mc{X})$. For this reason spectral series regression over Sobolev spaces has been well-studied---at least when $\mc{X} = [0,1]^d$---since at least \citet{rice1984},\footnote{And proposed much earlier in the context of density estimation by \citet{cencov1962}.} and its optimality properties are by this point generally well-understood.

PCR-LE serves as an empirical approximation to spectral series regression, since as already mentioned the eigenvectors $v_k$ are empirical approximations to the eigenfunctions $\psi_k$ of $\Delta_P$. Viewed in this light, a major advantage of PCR-LE is that it operates without needing knowledge of the design distribution $P$. This is an advantage because in our context $P$ is an unknown and potentially complex distribution: for example, it can be highly non-uniform, have a complicated support which may be a submanifold of $\Rd$, or both. In contrast, spectral series regression relies on diagonalizing the density-weighted Laplacian $\Delta_P$, and in our context must be viewed as an oracle method; to emphasize this we henceforth refer to the estimator defined in~\eqref{eqn:population-level_spectral_series} as \emph{population-level spectral series regression}. On the other hand, intuitively PCR-LE incurs some extra error by using an empirical approximation to the underlying basis $\{\psi_k\}_{k = 1}^{\infty}$: our work shows that in many cases, this extra error is not enough to change the overall rate of convergence.

\subsection{Main contributions}
Summarized succinctly, our main contribution is to theoretically analyze nonparametric regression with PCR-LE and establish upper bounds which imply that this method often achieves optimal rates of convergence over Sobolev spaces.

\paragraph{Rates of convergence: population-level spectral series regression.}
As we have already mentioned, the minimax optimal rates over Sobolev spaces are generally well-known, as are upper bounds for population-level spectral series methods which match these rates. However, we could not find precisely stated results applying to our setting, which is quite general in the following respects.
\begin{enumerate}
	\item We consider Sobolev spaces $H^s(\mc{X})$ for all combinations of $s$ and $d$. This includes the subcritical regime where the smoothness parameter $s$ satisfies $s < d/2$; in this regime $H^s(\mc{X})$ does not continuously embed into the space of continuous functions $C^0(\mc{X})$.
	\item We consider general design distributions $P$, which may satisfy certain regularity conditions but are not limited to being, say, the uniform distribution over $[0,1]^d$. 
\end{enumerate}
For completeness, we analyze population-level spectral series methods in this general setting, and establish upper bounds showing that such methods converge at the ``usual'' rates of $n^{-2s/(2s + d)}$ for estimation and $n^{-4s/(4s + d)}$ for testing. This analysis relies heavily on certain asymptotic properties of the continuum eigenfunctions $\psi_k$ and eigenvalues $\rho_k$, which hold for quite general second-order differential operators $\mc{L}$ including the density-weighted Laplacian $\mc{L} = \Delta_P$.

\paragraph{Rates of convergence: PCR-LE.}
The rest of our results consist of various upper bounds on the rates of convergence for the PCR-LE estimator $\wh{f}$, and a test using the statistic $\wh{T} = \|\wh{f}\|_2^2$. These upper bounds quantify two important properties of PCR-LE: first, that it can take advantage of smooth higher-order derivatives, and second that it can adapt to low intrinsic dimension of the design distribution, each in an optimal manner. We consider two models for the design distribution $P$, the flat Euclidean and manifold models (See Section~\ref{subsec:regression_laplacian_eigenmaps} for the formal definitions of these models). In the first model, the design distribution $P$ has support $\mc{X}$ which is a full-dimensional set in $\Rd$. In this case, our main contributions are as follows:
\begin{itemize}
	\item  Over a ball in the Sobolev space $H^{s}(\mc{X})$, we establish that the PCR-LE estimator $\wh{f}$ has in-sample mean-squared error on the order of $n^{-2s/(2s + d)}$, for any number of derivatives $s \in \mathbb{N}$ and dimension $d$ (Theorems~\ref{thm:laplacian_eigenmaps_estimation_fo} and~\ref{thm:laplacian_eigenmaps_estimation_ho}).
	\item We show that a test based on the statistic $\|\wh{f}\|_2^2$ has a squared critical radius on the order of $n^{-4s/(4s + d)}$, for any number of derivatives $s \in \mathbb{N}$ and dimension $d \in \{1,2,3,4\}$ (Theroems~\ref{thm:laplacian_eigenmaps_testing_fo} and~\ref{thm:laplacian_eigenmaps_testing_ho}).
\end{itemize}
We then consider the behavior of PCR-LE when the data satisfies a \emph{manifold hypothesis}, meaning the design distribution is supported on an (unknown) domain $\mc{X}$ which is a submanifold of $\Rd$ of intrinsic dimension $m \in \mathbb{N}, m < d$. In this case, our main contributions are as follows:
\begin{itemize}
	\item Over a ball in the Sobolev space $H^{s}(\mc{X})$, the PCR-LE estimator $\wh{f}$ has in-sample mean squared error of at most $n^{-2s/(2s + m)}$, when $s \in \{1,2,3\}$ and for any $m \in \mathbb{N}$ (Theorem~\ref{thm:laplacian_eigenmaps_estimation_manifold}). 
	\item A test based on the statistic $\|\wh{f}\|_2^2$ has a squared critical radius on the order of $n^{-4s/(4s + m)}$, when $s \in \{1,2,3\}$ and $m \in \{1,2,3,4\}$ (Theorem~\ref{thm:laplacian_eigenmaps_testing_manifold}).
\end{itemize}
To the best of our knowledge, the minimax rates for nonparametric regression with random design over unknown manifolds have only been worked out for H\"{o}lder classes, and even in this case the calculations are only for $s \leq 2$ bounded derivatives~\citep{bickel2007,yang2016}. Our upper bounds confirm that these rates are the same for Sobolev spaces---in estimation, when loss is measured in empirical norm---for the values of $s$ and $m$ mentioned above.

In all these cases, our bounds also depend optimally on the radius $M$ of the Sobolev ball under consideration. However, for some values of $s$ (number of derivatives) and $d$ (dimension), there do exist gaps between our upper bounds on the error of PCR-LE and the minimax rates. Although we do not give corresponding lower bounds verifying the tightness of our analysis, we believe these gaps reflect the true behavior of the method rather than some looseness in our analysis, and we comment more on this at relevant parts in the text. For completeness, we summarize all of our upper bounds---those which match the minimax rates, and those which do not---in Tables~\ref{tbl:estimation_rates} and~\ref{tbl:testing_rates}.

\begin{table}
	\begin{center}
		\begin{tabular}{p{.2\textwidth} | p{.14\textwidth} p{.12\textwidth} }
			Smoothness order & Flat Euclidean (Model~\ref{def:model_flat_euclidean}) & Manifold (Model~\ref{def:model_manifold}) \\
			\hline
			$s \leq 3$ & ${\bf n^{-2s/(2s + d)}}$ & ${\bf n^{-2s/(2s + m)}}$ \\
			$s > 3$  & ${\bf n^{-2s/(2s + d)}}$ & $n^{-6/(6 + m)}$
		\end{tabular}
	\end{center}
	\caption{Summary of PCR-LE estimation rates over Sobolev balls. Bold font marks minimax optimal rates. In each case, rates hold for all $d \in \mathbb{N}$ (under Model~\ref{def:model_flat_euclidean}), and for all $m \in \mathbb{N}, 1 < m < d$ (under Model~\ref{def:model_manifold}). Although we suppress it for simplicity, in all cases when the PCR-LE estimator is optimal, the dependence of the error rate on the radius $M$ of the Sobolev ball is also optimal.}
	\label{tbl:estimation_rates}
\end{table}
\begin{table}
	\begin{center}
		\begin{tabular}{p{.175\textwidth} p{.175\textwidth} | p{.14\textwidth} p{.12\textwidth} }
			Smoothness order & Dimension & Flat Euclidean (Model~\ref{def:model_flat_euclidean}) & Manifold (Model~\ref{def:model_manifold}) \\
			\hline
			\multirow{2}{*}{$s = 1$} & $\dim(\mc{X}) < 4$ & ${\bf n^{-4s/(4s + d)}}$ & ${\bf n^{-4s/(4s + m)}}$ \\
			& $\dim(\mc{X}) \geq 4$ & ${\bf n^{-1/2}}$ & ${\bf n^{-1/2}}$ \\
			\hline
			\multirow{3}{*}{$s = 2$ or $3$} & $\dim(\mc{X}) \leq 4$  & ${\bf n^{-4s/(4s + d)}}$ & ${\bf n^{-4s/(4s + m)}}$ \\
			& $4 <\dim(\mc{X}) < 4s$  & $n^{-2s/(2(s - 1) + d)}$ & $n^{-2s/(2(s - 1) + m)}$\\
			& $\dim(\mc{X}) \geq 4s$ & ${\bf n^{-1/2}}$ & ${\bf n^{-1/2}}$ \\
			\hline
			\multirow{3}{*}{$s > 3$} & $\dim(\mc{X}) \leq 4$ & ${\bf n^{-4s/(4s + d)}}$ & $n^{-12/(12 + d)}$ \\
			& $4 < \dim(\mc{X}) < 4s$ & $n^{-2s/(2(s - 1) + d)}$ & $n^{-6/(4 + m)}$ \\
			& $\dim(\mc{X}) \geq 4s$ & ${\bf n^{-1/2}}$ & ${\bf n^{-1/2}}$ \\
		\end{tabular}
	\end{center}
	\caption{Summary of PCR-LE testing rates over Sobolev balls. Bold font marks minimax optimal rates. Rates when $d > 4s$ assume that $f_0 \in L^4(\mc{X})$, and depend on $\|f_0\|_{L^4(\mc{X})}$. Although we suppress it for simplicity, in all cases when othe PCR-LE test is optimal, the dependence of the error rate on the radius $M$ of the Sobolev ball is also optimal.}
	\label{tbl:testing_rates}
\end{table}

\paragraph{Perspective: regression error versus feature reconstruction.}
We now pause for a moment, to emphasize that in a certain respect the aforementioned rates of convergence for PCR-LE are quite surprising. Remember that PCR-LE is a regression method using features (eigenvectors $v_k$ of the graph Laplacian $L$) which are themselves empirical estimates of population-level quantities (eigenfunctions $\psi_k$ of the density-weighted Laplacian $\Delta_P$). It seems reasonable to expect that the error of PCR-LE should be decomposed into two parts: first, the error with which these empirically-derived features estimate their continuum limits; second, the error with which, given ideal population-level features, the regression function is learned.

Crucially, our analysis \emph{does not} work in this way. This is important because all known upper bounds on the rates at which $v_k \to \psi_k$ as $n \to \infty$ are much slower than the minimax rates for regression over Sobolev classes. For instance, the best currently known upper bound on the empirical $L^2$ error $\frac{1}{n}\sum_{i = 1}^{n}(\sqrt{n} v_{k,i} - \psi_k(X_i))^2$ is only on the order of $n^{-2/(4 + d)}$~\citep{cheng2021}, which is slower than the minimax estimation rate over $H^s(\mc{X})$ for any $s \in \mathbb{N}, s \geq 1$.\footnote{To make matters worse, PCR-LE, when deployed optimally, does not use a single eigenvector $v_k$ for a fixed index $k \in \mathbb{N}$, but rather many eigenvectors $v_1,\ldots,v_K$ with $K$ growing in $n$. As $K$ grows larger, the rate at which $v_K \to \psi_K$ gets slower, since the population-level object being estimated is less regular; see~\citep{burago2014,trillos2019}.} Although this upper bound may not reflect the true rate of convergence of graph Laplacian eigenvectors---this is still an active area of research, and no lower bounds are known---it seems very unlikely that the true rate matches the minimax estimation rate $n^{-2s/(2s + d)}$, which after all approaches the dimension-free rate $1/n$ for large values of $s$. The bottom line is that the rate at which graph Laplacian eigenvectors are known to converge to density-weighted Laplacian eigenfunctions is too slow to explain the upper bounds we establish for PCR-LE.

Instead of relying on convergence of eigenvectors to eigenfunctions, our analysis proceeds via a bias-variance decomposition at the level of the graph. As usual for OLS estimates, the variance term depends only on the degrees of freedom $\mathrm{df}(\wh{f}) = \mathrm{tr}(V_KV_K^{\top}) = K$. More surprisingly, the bias can also be upper bounded without appealing to concentration of eigenvectors $v_1,\ldots,v_K$ around eigenfunctions $\psi_1,\ldots,\psi_K$; for instance, we show in Lemma~\ref{lem:fixed_graph_estimation} that for estimation the squared bias is at most on the order of $f_0^{\top} L^s f_0/(n\lambda_{K + 1}^s)$. 

Ultimately our upper bound on the error of PCR-LE is determined entirely by a pair of graph functionals: the quadratic form $f_0^{\top}L^s f_0$, and the graph Laplacian eigenvalue $\lambda_{K + 1}$. This brings a couple of advantages:
\begin{itemize}
	\item First, it eliminates the need to analyze convergence of eigenvectors to eigenfunctions, which is critical in order to get sufficiently fast rates of convergence for PCR-LE, as we have already explained. Instead, we only have to consider these two graph functionals, both of which are known to converge at faster rates than graph Laplacian eigenvectors. 
	\item Second, in order to obtain upper bounds on $\|\wh{f} - f_0\|_n^2$ we do not require that these graph functionals themselves converge to population-level limits, but only that they be stochastically bounded on the right order. The latter is a much weaker requirement. 
\end{itemize}
To derive our upper bounds on the error of PCR-LE, we directly analyze the quadratic form $f_0^{\top}L^s f_0$ and the eigenvalue $\lambda_{K + 1}$, using some existing results as well as deriving some new ones which may be of independent interest.

To summarize, our work demonstrates, broadly speaking, that regression using estimated features can be analyzed independently from the estimation error of the features themselves. Regression using learned features---that is, a feature representation derived from the data itself---is a general and widely applied paradigm, and we believe this observation may have consequences outside of its application to PCR-LE in this work.

\subsection{Related work}

\paragraph{Laplacian smoothing.}
In a previous paper~\citep{green2021}, we (the authors) considered an alternative method for nonparametric regression via neighborhood graphs: \emph{Laplacian smoothing}, defined as the solution to the following optimization problem,
\begin{equation}
\label{eqn:laplacian_smoothing}
\minimize_{f \in \Reals^n} \|{\bf Y} - f\|_2^2 + \lambda f^{\top} L f.
\end{equation}
Laplacian smoothing is penalized method for regression, where the penalty functional $f^{\top} L f$ serves as a discrete approximation to the continuum functional $J(f) := \int \|\nabla f(x)\|^2 p^2(x) \,dx$ \citep{bousquet03}. In the univariate setting ($d = 1$), this casts Laplacian smoothing as a discrete and density-weighted alternative to a first-order thin-plate spline estimator, which is defined as the solution to
\begin{equation}
\label{eqn:thin_plate_spline}
\minimize_{f \in H^1(\Reals)} \frac{1}{n}\sum_{i = 1}^{n} (Y_i - f(X_i))^2 + J(f).
\end{equation}
When $d = 1$ the first-order thin-plate spline estimator enjoys excellent theoretical properties, such as being minimax optimal over the first-order Sobolev space $H^1(\Reals)$. However, when $d \geq 2$ the story changes dramatically: the problem ~\eqref{eqn:thin_plate_spline} is in fact not even well-posed.\footnote{This can be explained by reference to the Sobolev Embedding Theorem, since it is an implication of this theorem that convergence of a sequence of functions $\{f_N\}_{N \in \mathbb{N}} \to f$ in first-order Sobolev norm implies pointwise convergence only when $d = 1$.} In contrast, in this previous paper, we showed that Laplacian smoothing was a well-posed and consistent estimator for any (fixed) dimension $d$, and achieved minimax optimal rates for estimation and testing so long as $d \in \{1,2,3,4\}$.

However, Laplacian smoothing neither takes advantage of smooth higher-order derivatives, nor is it provably optimal over $H^1(\mc{X})$ for dimensions $d \geq 5$. One of our motivations for considering PCR-LE was to find an estimator which addressed these deficiencies. In this work we indeed establish that PCR-LE has much stronger optimality properties than those we derived for Laplacian smoothing, or indeed those known for any other method of regression using neighborhood graphs. 

One way to interpret this difference between PCR-LE and Laplacian smoothing is to view the latter as a ridge regression problem. This follows from writing the Laplacian smoothing penalty as a (weighted) ridge penalty in the spectral domain,  $f^{\top} L f  = \sum_{k = 1}^{n} \lambda_k(v_k^{\top}f)^2$. \citet{dhillon2013} establish conditions under which principal components regression can have smaller risk than ridge regression using the same set of features. Viewed in this light, our work shows this phenomenon occurs when the features are eigenvectors of a neighborhood graph Laplacian and the estimand is a function in Sobolev space. It also establishes that principal components regression can obtain the minimax rate of convergence even when ridge fails to do so. Interestingly, this is not the case if the function class in question is an RKHS~\citep{dicker2017}, and further motivates the study of regression over Sobolev spaces in the subcritical regime, where surprising new phenomena emerge.

\paragraph{Other related work.}

Much of the work regarding regression using neighborhood graph Laplacians deals with \emph{semi-supervised learning}, where in addition to the labeled data $(X_1,Y_1),\ldots,(X_n,Y_n)$ one observes unlabeled points $(X_{n + 1},\ldots,X_{N})$, and the task is to produce an estimate at labeled and unlabeled points alike. To this end, the landmark paper of \cite{zhu2003semisupervised} proposed to interpolate the observed values by~\emph{harmonic extension}, i.e. compute the Laplacian matrix $L_N$  corresponding to a graph formed over all design points $X_1,\ldots,X_N$, and then solve the constrained problem
\begin{equation*}
\minimize_{f \in \Reals^N} f^{\top} L_N f \quad \mathrm{subject\,\,to}~~~ f_i = Y_i~~\textrm{for $i = 1,\ldots,n$.}
\end{equation*}
Conventional wisdom says that harmonic extension is sensible only when the responses are noiseless, $Y_i = f_0(X_i)$, and that in the noisy setting one should instead solve the penalized formulation
\begin{equation}
\label{eqn:graph_laplacian_regularization_ssl}
\minimize_{f \in \Reals^N} \sum_{i = 1}^{n}(Y_i - f_i)^2 + \lambda f^{\top} L_N f.
\end{equation}
Notwithstanding their intuitive appeal, both the constrained and penalized problems have issues when $d > 1$ and $n/N \to 0$: the estimates tend towards degeneracy, meaning they are ``spiky'' at labeled data points and close to constant everywhere else~\citep{nadler09,calder2019b, calder2020}. One solution to this problem is to instead use Laplacian Eigenmaps for semi-supervised learning (SSL-LE), i.e. compute the eigendecomposition $L_N = \sum_{k = 1}^{N} \lambda_k u_k u_k^{\top}$ and, letting $U \in \Reals^{n \times K}$ be the matrix with entries $U_{ik} = u_{k,i}$ and columns $U_1,\ldots,U_K$, solve the problem
\begin{equation}
\label{eqn:laplacian_eigenmaps_ssl}
\minimize_{f \in \mathrm{span}\{U_1,\ldots,U_K\}} \sum_{i = 1}^{n}(Y_i - f_i)^2.
\end{equation}
\cite{zhou2011,lee2016} analyze SSL-LE in a particular asymptotic regime where the number of labeled points $n$ is held fixed while the number of unlabeled points $N - n \to \infty$. They show that the SSL-LE estimator achieves minimax optimal rates---as a function of the number of labeled points $n$---over Sobolev spaces. However, in the particular asymptotic regime when $n$ is fixed and $N - n \to \infty$, the $n$ lowest-frequency eigenvectors of the graph Laplacian $L_N$ all converge to their continuum limits. Consequently, the SSL-LE estimator converges to the population-level spectral series estimator, and the analysis of SSL-LE reduces to that of the population-level method. As we have already explained, the supervised setting (where $N = n$) we consider in this work is very different, and analyzing PCR-LE necessitates an entirely different approach, 

In this supervised setting, there has been relatively little work regarding \emph{random design} regression with neighborhood graph Laplacians . Aside from our own work on Laplacian smoothing, summarized above, we highlight two other related papers: \citet{lee2016}, who analyze a variant of PCR-LE, but derive suboptimal rates of convergence, and \citet{trillos2020}, who study Laplacian smoothing and establish the uniform upper bound $\max_{i = 1,\ldots,n}|\wc{f}(X_i) - f_0(X_i)| \leq C n^{-2/(2 + d)}$ under the assumption $f_0 \in C^2(\mc{X})$, which is slower than the minimax rate $n^{-4/(4 + d)}$ for this function class. 

Most work on supervised learning using graphs adopts a \emph{fixed design} perspective, treating the design points $X_1 = x_1,\ldots,X_n = x_n$ as vertices of a fixed graph, and carrying out inference with respect to the conditional mean vector $(f_0(x_1),\ldots,f_0(x_n))$. In this setting, matching upper and lower bounds have been established that certify the optimality of graph-based methods for estimation \citep{wang2016,hutter2016,sadhanala16,sadhanala17,kirichenko2017,kirichenko2018}) and testing \citep{sharpnack2010identifying,sharpnack2013b,sharpnack2013,sharpnack2015} over different ``function'' classes (in quotes because these classes really model the $n$-dimensional vector of evaluations). This setting is quite general, because the graph need not be a geometric graph defined on a vertex set which belongs to Euclidean space. On the other hand, depending on the data collection process, it may be unnatural to model the design points as being a priori fixed, and the estimand as being a vector which exhibits a discrete notion of ``smoothness'' over this fixed design. Instead, we adopt the \emph{random design} perspective, and seek to estimate a function that we assume exhibits a more classical notion of smoothness. 

\paragraph{Roadmap.}
We now outline the structure of the rest of this paper. In Section~\ref{sec:setup_main_results}, we give our formal modeling assumptions, and precisely define the PCR-LE estimator and test we study. Propositions~\ref{prop:spectral_series_estimation} and ~\ref{prop:spectral_series_testing}, in Section~\ref{subsec:spectral_projection}, show that under rather general (nonparametric) conditions on the design distribution, population-level spectral series methods achieve minimax rates of convergence over Sobolev classes. Then in Sections~\ref{sec:minimax_optimal_laplacian_eigenmaps} and ~\ref{sec:manifold_adaptivity} we give our main upper bounds on the error of PCR-LE. These upper bounds (summarized above) hold under similarly general conditions, and imply that the PCR-LE estimator and test are also minimax rate-optimal. In Section~\ref{sec:experiments} we examine the empirical behavior of PCR-LE, and show that even at moderate sample sizes PCR-LE is competitive with population-level spectral series regression. We conclude with some discussion in Section~\ref{sec:discussion}. 

\paragraph{Notation.}
We now introduce some notation; for ease of reference, we include a table summarizing notation in Appendix~\ref{sec:notation_table}.

We frequently refer to various classical function classes, starting with the Lebesgue space $L^2(\mc{X})$, defined differently depending on whether $\mc{X} \subseteq \Rd$ is a full-dimensional open set or a compact Riemannian manifold. When $\mc{X} \subseteq \Rd$ is a full-dimensional open set, letting $\,d\nu$ denote the Lebesgue measure, the space $L^2(\mc{X})$ refers to the set of $\nu$-measurable functions $f$ for which $\|f\|_{L^2(\mc{X})}^2 := \int f^2 \,d\nu  < \infty$. When $\mc{X}$ is a compact Riemannian manifold, letting $\,d\mu$ denote the volume form induced by the embedding of $\mc{X}$ into $\Rd$, the space $L^2(\mc{X})$ refers to the set of $\mu$-measurable functions $f$ for which $\|f\|_{L^2(\mc{X})}^2 := \int f^2 \,d\mu  < \infty$. We also define an inner-product over these spaces: for a measure $P$ which admits a density $p$ with respect to $\nu$, we define $\dotp{f}{g}_P := \int f(x)g(x)p(x) \,d\nu(x)$; likewise, if $P$ admits a density $p$ with respect to $\mu$, $\dotp{f}{g}_P := \int f(x)g(x)p(x) \,d\mu(x)$. We refer to the norm $\|f\|_P^2 := \dotp{f}{f}_{P}$ as $L^2(P)$-norm.

We use $C^k(\mc{X})$ to refer to functions which are $k$ times continuously differentiable in $\mc{X}$, either for some integer $k \geq 1$ or for $k = \infty$. We let $C_c^{\infty}(\mc{X})$ represent those functions in $C^{\infty}(\mc{X})$ with support $V$ compactly contained in $\mc{X}$, meaning $\wb{V}$ is compact and $\wb{V} \subseteq \mc{X}$. We write $\partial f/\partial r_i$ for the partial derivative of $f$ in the $i$th standard coordinate of $\Rd$, and use the multi-index notation $D^{\alpha}f := \partial^{|\alpha|}f/\partial^{\alpha_1}x_1\ldots\partial^{\alpha_d}x_d$ for multi-indices $\alpha \in \Reals^d$. Recall that for a given multi-index $\alpha \in \mathbb{N}^d$, a function $f$ is \emph{$\alpha$-weakly differentiable} if there exists some $h \in L^1(\mc{X})$ such that
\begin{equation*}
\int_{\mc{X}} h g = (-1)^{|\alpha|} \int_{\mc{X}} f D^{\alpha}g, \quad \textrm{for every $g \in C_c^{\infty}(\mc{X})$.}
\end{equation*}
If such a function $h$ exists, it is the $\alpha$th weak partial derivative of $f$, and denoted by $D^{\alpha}f := h$. For functions $f$ which are $|\alpha|$-times classically differentiable, this coincides with the classical definition of derivative, and so we use the same notation for both. 

We write $\|\cdot\| = \|\cdot\|_2$ for Euclidean norm, $|\cdot| = \|\cdot\|_1$ for $\ell_1$ norm, and $d_{\mc{X}}(x',x)$ for the geodesic distance between points $x$ and $x'$ on a manifold $\mc{X}$. Then for a given $\delta > 0$, $B(x,\delta)$ is the radius-$\delta$ ball with respect to Euclidean distance, whereas $B_{\mc{X}}(x,\delta)$ is the radius-$\delta$ ball with respect to geodesic distance. Letting $T_x(\mc{X})$ be the tangent space at a point $x \in \mc{X}$, we write $B_m(v,\delta) \subset T_x(\mc{X})$ for the radius-$\delta$ ball centered at $v \in T_x(\mc{X})$.

For sequences $(a_n)$ and $(b_n)$, we use the asymptotic notation $a_n \lesssim b_n$ to mean that there exists a number $C$ such that $a_n \leq C b_n$ for all $n$. We write $a_n \asymp b_n$ when $a_n \lesssim b_n$ and $b_n \lesssim a_n$. On the other hand we write $a_n = o(b_n)$ when $\lim a_n/b_n = 0$, and likewise $a_n = \omega(b_n)$ when $\lim a_n/b_n = \infty$. Finally $a \vee b := \max\{a,b\}$ and $a \wedge b := \min\{a,b\}$.

\section{Preliminaries}
\label{sec:setup_main_results}

We begin in Sections~\ref{subsec:regression_laplacian_eigenmaps}-\ref{subsec:laplacian_eigenmaps} by precisely defining the models (random design points, Sobolev regression functions) and methods (Principal Components Regression with Laplacian Eigenmaps) under consideration. Then in Section~\ref{subsec:spectral_projection}, we analyze the behavior of population-level spectral series methods.

\subsection{Nonparametric regression over Sobolev spaces}
\label{subsec:regression_laplacian_eigenmaps}

As mentioned, we will always operate in the usual setting of nonparametric regression with random design. We observe independent random samples $(X_1,Y_1),\ldots,(X_n,Y_n)$, where the design points $X_1,\ldots,X_n$ are sampled from a distribution $P$ with support $\mc{X} \subseteq \Rd$, and the responses follow~\eqref{eqn:model}. We now formulate two models for the design distribution $P$ and regression function $f_0$: the \emph{flat Euclidean} and \emph{manifold} models.

\paragraph{Flat Euclidean model.}
In Definitions~\ref{def:model_flat_euclidean}-\ref{def:zero_trace_sobolev_space}, we collect the assumptions we make when working under the flat Euclidean model. We begin by giving some regularity conditions on the design.

\begin{definition}[Flat Euclidean model]
	\label{def:model_flat_euclidean}
	The support $\mc{X}$ of the design distribution $P$ is an open, connected, and bounded subset of $\Rd$, with Lipschitz boundary. The distribution $P$ admits a Lipschitz density $p$ with respect to the $d$-dimensional Lebesgue measure $\nu$, which is bounded away from $0$ and $\infty$,
	\begin{equation*}
	0 < p_{\min} \leq p(x) \leq p_{\max} < \infty, \quad \textrm{for all $x \in \mc{X}$.}
	\end{equation*}
\end{definition}
At various points we will also assume that the density $p \in C^k(\mc{X})$. On the other hand, we model the regression function as belonging to an order-$s$ Sobolev space, and being bounded in Sobolev norm.
\begin{definition}[Sobolev space on a flat Euclidean domain]
	\label{def:sobolev_space}
	For an integer $s \geq 1$, a function $f \in L^2(\mc{X})$ belongs to the Sobolev space $H^s(\mc{X})$ if for all $|\alpha| \leq s$, the weak derivatives $D^{\alpha}f$ exist and satisfy $D^{\alpha}f \in L^2(\mc{X})$. The $j$th order semi-norm for $f \in H^s(\mc{X})$ is $|f|_{H^j(\mc{X})} := \sum_{|\alpha| = j}\|D^{\alpha}f\|_{L^2(\mc{X})}$, and the corresponding norm
	\begin{equation*}
	\|f\|_{H^s(\mc{X})}^2 := \|f\|_{L^2(\mc{X})}^2 + \sum_{j = 1}^{s} |f|_{H^j(\mc{X})}^2,
	\end{equation*}
	induces the Sobolev ball
	\begin{equation*}
	H^s(\mc{X};M) := \bigl\{f \in H^s(\mc{X}): \|f\|_{H^s(\mc{X})} \leq M\bigr\}.
	\end{equation*} 
\end{definition}
When $s > 1$ we will also assume that $f_0$ satisfies a zero-trace boundary condition. Recall that $H^s(\mc{X})$ can alternatively be defined as the completion of $C^{\infty}(\mc{X})$ in the Sobolev norm $\|\cdot\|_{H^s(\mc{X})}$. The zero-trace Sobolev spaces are defined in a similar fashion, as the completion of $C_c^{\infty}(\mc{X})$ in the same norm.

\begin{definition}[Zero-trace Sobolev space]
	\label{def:zero_trace_sobolev_space}
	A function $f \in H^s(\mc{X})$ belongs to the zero-trace Sobolev space $H_0^s(\mc{X})$ if there exists a sequence $f_1,f_2,\ldots$ of functions in $C_c^{\infty}(\mc{X})$ such that
	\begin{equation*}
	\lim_{k \to \infty}\|f_k - f\|_{H^s(\mc{X})} = 0.
	\end{equation*}
	The normed ball $H_0^{s}(\mc{X};M) := H_0^{s}(\mc{X}) \cap H^{s}(\mc{X};M)$.
\end{definition}
Boundary conditions play an important role in the analysis of spectral methods, as we explain further in Section~\ref{subsec:spectral_projection}. For now, we limit ourselves to pointing out that for functions $f \in C^\infty(\mc{X})$, the zero-trace condition can be stated more concretely, as implying that $\partial^{k}f/\partial{\bf n}^k(x) = 0$ for each $k = 0,\ldots,s - 1$, and for all $x \in \partial\mc{X}$. (Here $\partial/(\partial {\bf n})$ is the partial derivative operator in the direction of the normal vector $\mathbf{n}$.)

\subsubsection{Manifold model}
As in the flat Euclidean case, we start with some regularity conditions on the design. One such condition will be on the reach $R$ of the manifold $\mc{X}$, which we recall is defined as follows:
\begin{equation*}
R := \Bigl\{\sup_{r > 0}: \forall z \in \Rd, \inf_{x \in \mc{X}} \|z - x\| \leq r, ~\exists! y \in \mc{X}~\mathrm{s.t.}~\|z - y\| = \inf_{x \in \mc{X}} \|z - x\|\Bigr\}.
\end{equation*}
In words, the reach is the largest radius of a ball which can be rolled around the manifold $\mc{X}$.
\begin{definition}[Manifold model]
	\label{def:model_manifold}
	The support $\mc{X}$ of the design distribution $P$ is a closed, connected, and smooth Riemannian manifold (without boundary) embedded in $\Rd$, of intrinsic dimension $1 \leq m < d$, and with a positive reach $R > 0$. The design distribution $P$ admits a Lipschitz density $p$ with respect to the volume form $d\mu$ induced by the Riemannian structure of $\mc{X}$, which is bounded away from $0$ and $\infty$,
	\begin{equation*}
	0 < p_{\min} \leq p(x) \leq p_{\max} < \infty, \quad \textrm{for all $x \in \mc{X}$.}
	\end{equation*}
\end{definition}

There are several equivalent ways to define Sobolev spaces on smooth Riemannian manifolds. We will stick with a definition that parallels our setup in the flat Euclidean setting as much as possible. To do so, we first recall the notion of partial derivatives on a manifold, which are defined with respect to a local coordinate system. Letting $r_1,\ldots,r_m$ be the standard basis of $\Reals^m$, for a given chart $(\phi,U)$ (meaning an open set $U \subseteq \mc{X}$, and a smooth mapping $\phi: U \to \Reals^m$) we write $\phi =: (x_1,\ldots,x_m)$ in local coordinates, meaning $x_i = r_i \circ \phi$. Then we define the partial derivative $\partial f/\partial x_i$ of a function $f: \mc{X} \to \Reals$ at $x \in U$ to be
\begin{equation*}
\frac{\partial f}{\partial x_i}(x) := \frac{\partial(f \circ \phi^{-1})}{\partial r_i}\bigl(\phi(x)\bigr).
\end{equation*}
The right hand side should be interpreted in the weak sense of derivative. As before, we use the multi-index notation $D^{\alpha}f := \partial^{|\alpha|}f/\partial^{\alpha_1}x_1\ldots\partial^{\alpha_m}x_m$. 

\begin{definition}[Sobolev space on a manifold]
	\label{def:sobolev_space_manifold}
	A function $f \in L^2(\mc{X})$ belongs to the Sobolev space $H^{s}(\mc{X})$ if for all $|\alpha| \leq s$, the weak derivatives $D^{\alpha}f$ exist and satisfy  $D^{\alpha}f \in L^2(\mc{X})$. The $j$th order semi-norm $|f|_{H^j(\mc{X})}$, the norm $\|f\|_{H^s(\mc{X})}$, and the ball $H^s(\mc{X};M)$ are all defined as in Definition~\ref{def:sobolev_space}.
\end{definition}
The partial derivatives $D^{\alpha}f$ clearly depend on the choice of local coordinates, and so will the resulting Sobolev norm~$\|f\|_{H^s(\mc{X})}$. However, for our purposes the important thing is that regardless of the choice of local coordinates the resulting norms will be equivalent\footnote{Recall that norms $\|\cdot\|_1$ and $\|\cdot\|_2$ on a space $\mc{F}$ are said to be equivalent if there exist constants $c$ and $C$ such that
	\begin{equation*}
	c \|f\|_1 \leq \|f\|_2 \leq C \|f\|_1 \quad \textrm{for all $f \in \mc{F}$.}
	\end{equation*}} 
and so the ultimate Sobolev space $H^s(\mc{X})$ is independent of local coordinates. For more information regarding manifolds and Sobolev spaces defined thereupon, see~\cite{lee2013} and~\cite{hebey1996}.

\subsection{Principal Components Regression with Laplacian Eigenmaps (PCR-LE)}
\label{subsec:laplacian_eigenmaps}
We now formally define the estimator and test statistic we study. Both are derived from eigenvectors of a graph Laplacian.  For a positive, symmetric kernel $\eta: [0,\infty) \to [0,\infty)$, and a radius parameter $\varepsilon > 0$, let $G = ([n],W)$ be the neighborhood graph formed over the design points $\{X_1,\ldots,X_n\}$, with a weighted edge $W_{ij} = \eta(\|X_i - X_j\|/\varepsilon)$ between vertices $i$ and $j$. Then the 
\emph{neighborhood graph Laplacian} $L_{n,\varepsilon}: \Reals^n \to \Reals$ is defined by its action on vectors $u \in \Reals^n$ as
\begin{equation}
\label{eqn:neighborhood_graph_laplacian}
\bigl(L_{n,\varepsilon}u\bigr)_i := \frac{1}{n\varepsilon^{2 + \mathrm{dim}(\mc{X})}} \sum_{j = 1}^{n} \bigl(u_i - u_j\bigr) \eta\biggl(\frac{\|X_i - X_j\|}{\varepsilon}\biggr).
\end{equation}
(Here $\mathrm{dim}(\mc{X})$ stands for the dimension of $\mc{X}$. It is equal to $d$ under the assumptions of Model~\ref{def:model_flat_euclidean}, and equal to $m$ under the assumptions of Model~\ref{def:model_manifold}. The pre-factor $(n\varepsilon^{2 + \mathrm{dim}(\mc{X})})^{-1}$ ensures non-degenerate stable limits as $n \to \infty, \varepsilon \to 0$). Note that $(n\varepsilon^{\dim(\mc{X}) + 2}) \cdot L_{n,\varepsilon} = D - W$, where $D \in \Reals^{n \times n}$ is the diagonal degree matrix, $D_{ii} = \sum_{i = 1}^{n} W_{ij}$.

The graph Laplacian is a positive semi-definite matrix, and admits the eigendecomposition $L_{n,\varepsilon} = \sum_{k = 1}^{n} \lambda_k v_k v_k^{\top}$, where for each $k \in \{1,\ldots,n\}$ the eigenvalue-eigenvector pair $(\lambda_k,v_k)$ satisfies
\begin{equation*}
L_{n,\varepsilon}v_k = \lambda_k v_k, \quad \|v_k\|_2^2 = 1.
\end{equation*}
We will assume without loss of generality that each eigenvalue $\lambda$ of $L_{n,\varepsilon}$ has algebraic multiplicity $1$, and so we can index the eigenpairs $(\lambda_1,v_1),\ldots,(\lambda_n,v_n)$ in ascending order of eigenvalue, $0 = \lambda_1 < \ldots < \lambda_n$. 

The PCR-LE estimator $\wh{f}$ defined in~\eqref{eqn:pcr-le} simply projects the response vector ${\bf Y}$ onto the first $K$ eigenvectors of $L_{n,\varepsilon}$. Since the eigenvectors of the graph Laplacian are orthonormal with respect to the Euclidean inner product on $\Reals^n$, we can more simply write this as
\begin{equation}
\label{eqn:laplacian_eigenmaps_estimator}
\wh{f} = V_K V_K^{\top} {\bf Y},
\end{equation} 
where $V_K \in \Reals^{n \times K}$ is the matrix with $k$th column $V_{K,k} = v_k$. The PCR-LE test statistic is defined with respect to the empirical norm $\|f\|_{n}^2 := \frac{1}{n}\sum_{i = 1}^{n} \bigl(f(X_i)\bigr)^2$ of $\wh{f}$:\footnote{Here and throughout, when there is no chance of confusion we will identify vectors $f \in \Reals^n$ with functions $f: \{X_1,\ldots,X_n\} \to \Reals, f(X_i) = f_i$.}
\begin{equation}
\label{eqn:laplacian_eigenmaps_test}
\wh{T} := \|\wh{f}\|_n^2 = \frac{1}{n} {\bf Y}^{\top} V_K V_K^{\top} {\bf Y},
\end{equation}
and can be used in the \emph{signal detection} problem to distinguish whether or not $f_0 = 0$.

\subsection{Spectral series regression over Sobolev spaces}
\label{subsec:spectral_projection}
We now establish some upper bounds on the error of population-level spectral series regression when $f_0 \in H^s(\mc{X})$, which imply that such methods achieve optimal rates of convergence for both estimation and testing. The upper bounds we establish are ``usual'' in the sense that they match the rates $n^{-2s/(2s + d)}$ (estimation) and $n^{-4s/(4s + d)}$ (testing) which are already known in many cases. However, they are unusual in that we treat both the case where $s < d/2$ and thus the Sobolev space $H^s(\mc{X})$ does not continuously embed into $C(\mc{X})$, and the case where $P$ is not the uniform distribution over the unit cube. The upper bounds given in this section serve two purposes: first, to clarify what the rates are in these less-typically studied settings; second, to show that even in this general setting, population-level spectral series regression can always obtain the optimal rates. This latter point is important since the method we focus on for the most part, PCR-LE, is an empirical approximation to population-level spectral series regression.

\subsubsection{Spectrally defined Sobolev spaces}
Let $\mc{X}$ be an open domain which satisfies the conditions of Model~\ref{def:model_flat_euclidean}. Recalling the density-weighted Laplacian $\Delta$, defined in~\eqref{eqn:density-weighted-laplace}, we consider the eigenvector equation with Neumann boundary conditions,
\begin{equation}
\label{eqn:laplace_beltrami_eigenproblem}
\Delta_P\psi = \rho \psi, \quad \frac{\partial}{\partial{\bf n}}\psi = 0~~\textrm{on $\partial \mc{X}$.}
\end{equation}
Under Model~\ref{def:model_flat_euclidean}, the eigenvector equation~\eqref{eqn:laplace_beltrami_eigenproblem} has enumerable solutions $(\rho_1,\psi_1),(\rho_2,\psi_2),\ldots$, sorted as usual in ascending order of eigenvalue~\citep{garciatrillos18}. These eigenvalues and eigenfunctions can be used to give a spectral definition of Sobolev spaces. Consider the ellipsoid
\begin{equation}
\label{eqn:sobolev_ellipsoid}
\mc{H}^{s}(\mc{X}) := \Bigl\{\sum_{k = 1}^{\infty} a_k \psi_k \in L^2(\mc{X}):  \sum_{k = 1}^{\infty} a_k^2 \rho_k^s \leq M^2 \Bigr\},
\end{equation}
equipped with the norm $\|\sum_{k = 1}^{\infty} a_k \psi_k\|_{\mc{H}^s(\mc{X})}^2 = \sum_{k = 1}^{\infty} a_k^2 \rho_k^s$. Under appropriate regularity conditions $\mc{H}^s(\mc{X})$ consists of functions $f \in H^s(\mc{X})$ which also satisfy some additional boundary conditions. For instance, assuming Model~\ref{def:model_flat_euclidean}, $p \in C^{\infty}(\mc{X})$ and $\partial \mc{X} \in C^{1,1}$, \citet{dunlop2020} show that for any $s \geq 1$, the ellipsoid $\mc{H}^{2s}(\mc{X})$ satisfies
\begin{equation}
\label{eqn:sobolev_ellipsoid_to_sobolev_ball}
\mc{H}^{2s}(\mc{X}) = 
\biggl\{f \in H^{2s}(\mc{X}): \frac{\partial \Delta_P^rf}{\partial {\bf n}} = 0~\textrm{on}~\partial\mc{X},~~\textrm{for all $0 \leq r \leq s - 1$} \biggr\},
\end{equation}
and likewise $\mc{H}^{2s + 1}(\mc{X}) = \mc{H}^{2s}(\mc{X}) \cap H^{2s + 1}(\mc{X})$ for any $s \geq 0$. Additionally, the norms $\|\cdot\|_{\mc{H}^s(\mc{X})}$ and $\|\cdot\|_{H^s(\mc{X})}$ are equivalent.


\subsubsection{Estimation with spectral series regression}
Recall the population-level spectral series estimator $\wt{f}$ defined in~\eqref{eqn:population-level_spectral_series}. We now give an upper bound on the risk of $\wt{f}$, when loss is measured in $L^2(P)$ norm.
\begin{proposition}
	\label{prop:spectral_series_estimation}
	Suppose data is observed according to Model~\ref{def:model_flat_euclidean}. Assume additionally that $\partial \mc{X} \in C^{1,1}$, $p \in C^{\infty}(\mc{X})$, $f_0 \in \mc{H}^{s}(\mc{X};M)$ and $\|f_0\|_P^2 \leq 1$. Then there exists a constant $C$ which does not depend on $f_0$, $M$ or $n$ such that the following statement holds: if the population-level spectral series estimator $\wt{f}$ is computed with parameter $K = \max\{\floor{M^2n}^{d/(2s + d)},1\}$, then
	\begin{equation}
	\label{eqn:spectral_series_estimation}
	\Ebb\bigl[\|\wt{f} - f_0\|_P^2\bigr] \leq C \min\Bigl\{M^2\bigl(M^2n\bigr)^{-2s/(2s + d)}, \frac{1}{n}\Bigr\}.
	\end{equation}
\end{proposition}
When the Sobolev ball radius $n^{-1/2} \lesssim M$, the upper bound in~\eqref{eqn:spectral_series_estimation} is on the order of $M^2(M^2n)^{-2s/(2s + d)}$. This is well known to be the minimax rate of estimation over the Sobolev classes $H^s([0,1]^d;M)$ when $s > d/2$; see e.g.~\cite{gyorfi2006,wasserman2006,tsybakov08} and references therein, and specifically Theorem~3.2 of~\cite{gyorfi2006} for a matching lower bound in the context of nonparametric regression with random design. On the other hand there seems to have been much less study of minimax rates over $H^s([0,1]^d;M)$ when $s < d/2$. In this \emph{subcritical} regime, the Sobolev space contains functions without continuous representatives, and certain questions become more subtle; see our remark after Theorem~\ref{thm:laplacian_eigenmaps_estimation_ho}. However, Proposition~\ref{prop:spectral_series_estimation} confirms that in this regime the minimax rates (with loss measured in squared-$L^2(P)$ norm) are still on the order of $M^2(M^2n)^{-2s/(2s + d)}$, since a matching lower bound follows from the known estimation rates over $C^s([0,1]^d) \subseteq H^s([0,1]^d)$~\citep{stone1980}.

\subsubsection{Testing with spectral series regression}
In the goodness-of-fit testing problem, one asks for a test function---formally, a Borel measurable function $\phi$ that takes values in $\{0,1\}$--- which can distinguish between the hypotheses
\begin{equation}
\mathbf{H}_0: f_0 = f_0^{\star}, ~~\textrm{versus}~~ \mathbf{H}_a: f_0 \in \mc{H}^{s}(\mc{X};M) \setminus \{f_0^{\star}\}.
\end{equation} 
To fix ideas, here and throughout we focus on the signal detection problem, meaning the special case where $f_0^{\star} = 0$.\footnote{This is without loss of generality since all the test statistics we consider are easily modified to handle the case when $f_0^{\ast}$ is not $0$, by simply subtracting $f_0^{\ast}(X_i)$ from each observation $Y_i$, with no change in the analysis.} For more background on nonparametric goodness-of-fit testing problems, see~\cite{ingster2012}.

For the signal detection problem, the population-level spectral series test $\wt{\varphi} := \1\{\wt{T} \geq K/N + \sqrt{2K/an^2}\}$ has bounded Type I error, $\Ebb_{0}[\wt{\varphi}] \leq a (1 + o(1))$. Proposition~\ref{prop:spectral_series_testing} gives an upper bound on the Type II error that holds over all $f_0 \in \mc{H}^s(\mc{X};M)$ for which $\|f_0\|_P^2$ is sufficiently large.
\begin{proposition}
	\label{prop:spectral_series_testing}
	Suppose data is observed according to Model~\ref{def:model_flat_euclidean}, and that the density $p$ is known.  Suppose additionally that $\partial \mc{X} \in C^{1,1}$, $p \in C^{\infty}(\mc{X})$, $f_0 \in \mc{H}^{s}(\mc{X};M)$ for some $s > d/4$, and $\|f_0\|_{L^4(\mc{X})}^4 \leq 1$. Then there exists a constant $C$ which does not depend on $f_0$, $M$ or $n$ such that the following statement holds: if the population-level spectral series test $\wt{\varphi}$ is computed with parameter $K = \max\{\floor{M^2n}^{2d/(4s + d)},1\}$, and if
	\begin{equation}
	\label{eqn:spectral_series_testing}
	\|f_0\|_P^2 \geq C\min\Bigl\{M^2(M^2n)^{-4s/(4s + d)}, \frac{1}{n}\Bigr\}
	\end{equation}
	then the Type II error is upper bounded, $\Ebb_{f_0}[1 - \wt{\varphi}] \leq b$.
\end{proposition}
Assuming again that $n^{-1/2} \lesssim M$, the right hand side of~\eqref{eqn:spectral_series_testing} is $M^2(M^2n)^{-4s/(4s + d)}$, matching the usual minimax critical radius over Sobolev spaces (see e.g.~\cite{guerre02,ingster2009,ingster2012}. Specifically, \cite{ingster2009} show that the minimax squared critical radius is on the order of $n^{-4s/(4s + d)}$ when $M = 1$, and simple alterations of their analysis imply the rate $M^2(M^2n)^{-4s/(4s + d)}$ for general $M$.) On the other hand, when $s \leq d/4$ the minimax regression testing rates over $H^s(\mc{X})$ are not known. If one explicitly assumes $f_0 \in L^4(\mc{X};1)$---note that $H^s(\mc{X})$ does not continuously embed into $L^4(\mc{X})$ when $s \leq d/4$---then the minimax critical radius for regression testing is on the order of $n^{-1/2}$ \citep{guerre02}, and is achieved by a test using the naive statistic $\|{\bf Y}\|_n^2$. In other words, the regression testing problem over Sobolev spaces fundamentally changes when $s \leq d/4$, and hereafter when we discuss testing we will limit our consideration to $s > d/4$. 

The main takeaway from Propositions~\ref{prop:spectral_series_estimation} and~\ref{prop:spectral_series_testing} is that population-level spectral series methods for regression achieve optimal rates of convergence, when the regression function $f_0$ is Sobolev smooth and the design distribution $P$ is known a priori and satisfies an appropriate notion of smoothness.\footnote{The assumption $p \in C^{\infty}(\mc{X})$ could likely be weakened, but since this would not substantially add to the main points of Propositions~\ref{prop:spectral_series_estimation} and~\ref{prop:spectral_series_testing}, we do not pursue the details further.} We reiterate that when the design distribution is unknown, these methods have to be treated as oracle methods, in contrast to PCR-LE. As we will see, PCR-LE achieves comparable rates of convergence when $p$ is sufficiently smooth but potentially unknown.

Of course, it is worth pointing out that other methods besides PCR-LE are statistically optimal for nonparametric regression even when $p$ is unknown. We comment more on some of these in Section~\ref{sec:discussion}, after we have derived our major results regarding PCR-LE.

\section{Minimax Optimality of PCR-LE}
\label{sec:minimax_optimal_laplacian_eigenmaps}

In this section we give upper bounds on the error of PCR-LE in the flat Euclidean setting, where we observe data $(X_1,Y_1),\ldots,(X_n,Y_n)$ according to Model~\ref{def:model_flat_euclidean}. We will divide our theorem statements based on whether the regression function $f_0$ belongs to the first order Sobolev class $H^1(\mc{X})$ or a higher-order Sobolev class ($H_0^{s}(\mc{X})$ for some $s > 1$), since the details of the two settings are somewhat different.

\subsection{First-order Sobolev classes}
\label{sec:first_order_sobolev_classes}
We begin by assuming $f_0 \in H^1(\mc{X}; M)$. We show that $\wh{f}$ and a test based on $\wh{T}$ are minimax optimal, for all values of $d$ for which the minimax rates are known, and under no additional assumptions (beyond those of Model~\ref{def:model_flat_euclidean}) on the design distribution $P$.

\subsubsection{Estimation with PCR-LE} 
PCR-LE depends on the kernel $\eta$ and two tuning parameters, the graph radius $\varepsilon$ and number of eigenvectors $K$. We will need to make some assumptions on each.
\begin{enumerate}[label=(K\arabic*)]
	\setcounter{enumi}{0}
	\item
	\label{asmp:kernel_flat_euclidean}
	The kernel function $\eta$ is a nonincreasing function supported on $[0,1]$. Its restriction to $[0,1]$ is Lipschitz, and $\eta(1) > 0$. Additionally, it is normalized so that
	\begin{equation*}
	\int_{\Rd} \eta(\|z\|) \,dz = 1,
	\end{equation*}
	and we assume \smash{$\sigma_{\eta} := \frac{1}{d}\int_{\Rd} \|x\|^2 \eta(\|x\|) \,dx < \infty$}.
\end{enumerate}
\begin{enumerate}[label=(P\arabic*)]
	\setcounter{enumi}{0}
	\item 
	\label{asmp:parameters_estimation_fo} 
	For constants $c_0$ and $C_0$, the graph radius $\varepsilon$ and the number of eigenvectors $K$ satisfy the following inequalities:
	\begin{equation}\\
	\label{eqn:radius_fo} 
	C_0\biggl(\frac{\log n}{n}\biggr)^{1/d} \leq \varepsilon \leq c_0\min\{1,K^{-1/d}\},
	\end{equation}
	and 
	\begin{equation}
	\label{eqn:eigenvector_estimation_fo} 
	K = \min\Bigl\{\floor{(M^2n)^{d/(2 + d)}} \vee 1, n\Bigr\}.
	\end{equation}
\end{enumerate}
We comment on these assumptions after stating our first main theorem, regarding the estimation error of PCR-LE. The proof of this theorem, along with the proofs of all subsequent results, can be found in the Appendix.
\begin{theorem}
	\label{thm:laplacian_eigenmaps_estimation_fo}
	Suppose Model~\ref{def:model_flat_euclidean}, and additionally $f_0 \in H^1(\mc{X},M)$. There are constants $c,C$ and $N$ (not depending on $f_0$, $M$ or $n$), such that the following statement holds for all $n \geq N$ and any $\delta \in (0,1)$: if the PCR-LE estimator $\wh{f}$ is computed with a kernel $\eta$ satisfying~\ref{asmp:kernel_flat_euclidean}, and parameters $\varepsilon$ and $K$ satisfying~\ref{asmp:parameters_estimation_fo}, then
	\begin{equation}
	\label{eqn:laplacian_eigenmaps_estimation_fo}
	\|\wh{f} - f_0\|_n^2 \leq C\Bigl(\frac{1}{\delta}M^2(M^2n)^{-2/(2 + d)} \wedge 1\Bigr) \vee \frac{1}{n},
	\end{equation}
	with probability at least $1 - \delta - Cn\exp(-cn\varepsilon^d) - \exp(-K)$.
\end{theorem}
From~\eqref{eqn:laplacian_eigenmaps_estimation_fo} it follows immediately that when $n^{-1/2} \lesssim M \lesssim n^{1/d}$, then with constant probability $\|\wh{f} - f_0\|_n^2 \lesssim M^2(M^2n)^{-2/(2 + d)}$, matching the minimax estimation rate over Sobolev classes.

Some other remarks:
\begin{itemize}
	\item \emph{Radius of the Sobolev ball.} 
	When $M = o(n^{-1/2})$ then computing PCR-LE with $K = 1$ achieves the parametric rate $\|\wh{f} - f_0\|_n^2 \lesssim n^{-1}$, and the zero-estimator $\wh{f} = 0$ achieves the better rate $\|\wh{f} - f_0\|_n^2 \lesssim M^2$. However, we do not know what the minimax rate is in this regime. On the other hand, when $M = \omega(n^{1/d})$, then computing PCR-LE with $K = n$ achieves the rate $\|\wh{f} - f_0\|_n^2 \lesssim 1$, which is better than the rate in~\eqref{eqn:spectral_series_estimation}. This is because we are evaluating error in-sample rather than out-of-sample. However, in truth these are edge cases, which do not fall neatly into the framework of nonparametric regression. 
	
	\item \emph{In-sample error.} Since the PCR-LE estimator is defined only at the design points $\{X_1,\ldots,X_n\}$, we use the empirical norm $\|\cdot\|_n^2$ as our estimation loss. Depending on the problem at hand, it may be more interesting to consider loss in $L^2(P)$ norm, for instance because this loss is intrinsically tied to the prediction error $\Ebb_{X \sim P}(\wh{f}(X) - f_0(X))^2$. In some preliminary analysis, we have considered applying a generic kernel smoother $T$ to the PCR-LE estimate, so as to produce a bona-fide function $T\wh{f}: \Rd \to \Reals$; the error $\|T\wh{f} - f_0\|_P^2$ is provably on the order of the minimax estimation rate $n^{-2/(2 + d)}$. Indeed, our analysis is completely independent of the structure of $\wh{f}$, and could apply to any estimator $\wh{g} \in \Reals^n$ sufficiently close to $f_0$ in empirical norm. We intend to pursue the matter further in future work.
	
	\item \emph{Meaning of pointwise evaluation.} There is one subtlety introduced by the use of in-sample mean squared error. Since elements $f \in H^s(\mc{X})$ are equivalence classes, defined only up to a set of measure zero, one cannot really speak of the pointwise evaluation $f_0(X_i)$, as we do by defining our target of estimation to be $(f_0(X_1),\ldots,f_0(X_n))$, until one selects a representative of each equivalence class $f$. Implicitly, we will always pick the \emph{precise representative} $f_0^{\ast} \in f_0$ (as defined in~\cite{evans15}), and the notation ``$f_0(X_i)$'' should always be interpreted as $f_0^{\ast}(X_i)$. To be clear, however, it does not really matter which representative we choose, since all versions agree except on a set of measure zero, and so any two $g_0,h_0 \in f_0$ satisfy $g_0(X_i) = h_0(X_i)$ for all $i = 1,\ldots,n$ almost surely. For this reason we can write $f_0(X_i)$ without fear of ambiguity or confusion. 
	
	\item \emph{Tuning parameters}. The assumptions placed on the kernel function $\eta$ are needed for technical reasons. They can likely be weakened, although we note that they are already fairly general. The lower bound on $\varepsilon$ imposed by~\eqref{eqn:radius_fo} is on the order of the connectivity threshold, the smallest radius for which the resulting graph will still be connected with high probability. On the other hand, as we will see in Section~\ref{subsec:analysis}, the upper bound on $\varepsilon$ is needed to ensure that the graph eigenvalue $\lambda_K$ is of at least the same order as the continuum eigenvalue $\rho_K$; this is essential in order to obtain a tight upper bound on the bias of $\wh{f}$.  Finally, we set $K = \floor{(M^2n)^{d/(2 + d)}}$ (when possible) to optimally trade-off bias and variance, as is typical.
	
	In practice, one typically tunes hyper-parameters by cross-validation. However, because the estimator $\wh{f}$ is defined only in-sample, cross-validation cannot be used to tune parameters for PCR-LE. As previously mentioned, we are currently considering methods for extrapolation of $\wh{f}$ via kernel smoothing, which should allow for cross-validation and other sample-splitting techniques.
	
	
	\item \emph{High-probability guarantees}. The upper bound given in~\eqref{eqn:laplacian_eigenmaps_estimation_fo} holds with ``constant probability'', meaning with probability $1 - \delta - o(1)$. Under the stronger assumption that $f_0$ is $M$-Lipschitz, we can establish the same guarantee~\eqref{eqn:laplacian_eigenmaps_estimation_fo} with probability $1 - \delta^2/n - Cn\exp(-cn\varepsilon^d) - \exp(-K)$; in other words, we can give a high probability guarantee (for details see~\cite{green2021}). In this case a routine calculation shows that $\Ebb[\|\wh{f} - f_0\|_n^2]$ will also be on the some order as~\eqref{eqn:laplacian_eigenmaps_estimation_fo}. We also suspect that high-probability guarantees will hold so long as $\|\nabla f\|_{L^q(\mc{X})}$ is bounded for some sufficiently large $q < \infty$, but it remains an open question whether such guarantees can be obtained in the Sobolev case ($q = 2$) which is the focus of this work. 
\end{itemize}

\subsubsection{Testing with PCR-LE} 
Consider the test $\varphi := \1\{\wh{T} \geq t_{a}\}$, where $t_{a}$ is the threshold
\begin{equation*}
t_{a} := \frac{K}{n} + \frac{1}{n}\sqrt{\frac{2K}{a}}.
\end{equation*}
This choice of threshold $t_{a}$ guarantees that $\varphi$ is a level-$a$ test. As we show in Theorem~\ref{thm:laplacian_eigenmaps_testing_fo}, when $d < 4$, $\varepsilon$ and $K$ are chosen appropriately, and the alternative $f_0$ has is sufficiently well-separated from $0$, the test $\varphi$ has Type II error of at most $b$.

\begin{enumerate}[label=(P\arabic*)]
	\setcounter{enumi}{1}
	\item 
	\label{asmp:parameters_testing_fo}
	The graph radius $\varepsilon$ satisfies~\eqref{eqn:radius_fo}, and the number of eigenvectors 
	\begin{equation}
	\label{eqn:eigenvector_testing_fo}
	K = \min\Bigl\{\floor{(M^2n)^{2d/(4 + d)}} \vee 1, n\Bigr\}.
	\end{equation}
\end{enumerate}
\begin{theorem}
	\label{thm:laplacian_eigenmaps_testing_fo}
	Fix $a,b \in (0,1)$. Suppose Model~\ref{def:model_flat_euclidean}. Then $\mathbb{E}_0[\varphi] \leq a$, i.e $\varphi$ is a level-$a$ test. Suppose additionally $f_0 \in H^1(\mc{X};M)$, and that $d < 4$. Then there exist constants $C$ and $N$ that do not depend on $f_0$, such that the following statement holds for all $n \geq N$: if the PCR-LE test $\varphi$ is computed with a kernel $\eta$ satisfying~\ref{asmp:kernel_flat_euclidean}, and parameters $\varepsilon$ and $K$ satisfying~\ref{asmp:parameters_testing_fo}, and if $f_0$ satisfies
	\begin{equation}
	\label{eqn:laplacian_eigenmaps_testing_criticalradius_fo}
	\|f_0\|_P^2 \geq C\biggl(\Bigl(M^2(M^2n)^{-4/(4 + d) } \wedge n^{-1/2}\Bigr)\biggl[\sqrt{\frac{1}{a}} + \frac{1}{b}\biggr] \vee \frac{M^2}{b n^{2/d}} \biggr) \vee \frac{1}{n},
	\end{equation}
	then $\Ebb_{f_0}[1 - \varphi] \leq b$.
\end{theorem}
Although~\eqref{eqn:laplacian_eigenmaps_testing_criticalradius_fo} involves taking the maximum of several different terms, the important takeaway of Theorem~\ref{thm:laplacian_eigenmaps_testing_fo} is that if $n^{-1/2} \lesssim M \lesssim n^{(4 - d)/4d}$, then $\varphi$ has small worst-case risk as long as $f_0$ is separated from $0$ by at least $M^2(M^2n)^{-4/(4 + d)}$. This implies that $\varphi$ is a minimax rate-optimal test over $H^1(\mc{X};M)$ when $d \in \{1,2,3\}$. As mentioned previously, when $d \geq 4$ the first order Sobolev space $H^1(\mc{X})$ does not continuously embed into $L^4(\mc{X})$, and in this case the optimal rates for regression testing over Sobolev spaces are unknown.

\subsection{Higher-order Sobolev classes}
\label{sec:higher_order_sobolev_classes}
We now consider the situation where the regression function displays some higher-order regularity, $f_0 \in H_0^s(\mc{X};M)$. We show that the PCR-LE estimator and test continue to be optimal for all orders of $s$, as long as the design density is itself also sufficiently regular, $p \in C^{s - 1}(\mc{X})$. In estimation, this is the case for any dimension $d$, whereas in testing it is the case only when $d \leq 4$. 

\subsubsection{Estimation with PCR-LE} 
In order to show that $\wh{f}$ is an optimal estimator over $H_0^s(\mc{X};M)$, we will require that $\varepsilon$ be meaningfully larger than the lower bound in~\ref{asmp:parameters_estimation_fo}.
\begin{enumerate}[label=(P\arabic*)]
	\setcounter{enumi}{2}
	\item 
	\label{asmp:parameters_estimation_ho}
	For constants $c_0$ and $C_0$, the graph radius $\varepsilon$ and number of eigenvectors $K$ satisfy
	\begin{equation}
	\label{eqn:radius_ho}
	C_0\max\biggl\{\biggl(\frac{\log}{n}\biggr)^{1/d}, (M^2n)^{-1/(2(s - 1) + d)}\biggr\} \leq \varepsilon \leq c_0\min\{1, K^{-1/d}\}
	\end{equation}
	and
	\begin{equation*}
	K = \min\Bigl\{\floor{(M^2n)^{d/(2s + d)}} \vee 1,n\Bigr\}
	\end{equation*}
\end{enumerate}
Crucially, when $n$ is sufficiently large the two conditions in~\ref{asmp:parameters_estimation_ho} are not mutually exclusive.
\begin{theorem}
	\label{thm:laplacian_eigenmaps_estimation_ho}
	Suppose Model~\ref{def:model_flat_euclidean}, and additionally $f_0 \in H_0^s(\mc{X},M)$ and $p \in C^{s - 1}(\mc{X})$. There exist constants $c,C$ and $N$ that do not depend on $f_0$, such that the following statement holds all for all $n$ larger than $N$ and for any $\delta \in (0,1)$: if the PCR-LE estimator $\wh{f}$ is computed with a kernel $\eta$ satisfying~\ref{asmp:kernel_flat_euclidean}, and parameters $\varepsilon$ and $K$ satisfying~\ref{asmp:parameters_estimation_ho}, then
	\begin{equation}
	\label{eqn:laplacian_eigenmaps_estimation_ho}
	\|\wh{f} - f_0\|_n^2 \leq C\Bigl(\frac{1}{\delta}M^2(M^2n)^{-2s/(2s + d)} \wedge 1\Bigr) \vee \frac{1}{n},
	\end{equation}
	with probability at least $1 - \delta - Cn\exp(-cn\varepsilon^d) - \exp(-K)$.
\end{theorem}
Theorem~\ref{thm:laplacian_eigenmaps_estimation_ho}, in combination with Theorem~\ref{thm:laplacian_eigenmaps_estimation_fo}, implies that in the flat Euclidean setting PCR-LE is a minimax rate-optimal estimator over Sobolev classes, for all values of $s$ and $d$. Some other remarks:
\begin{itemize}
	\item \emph{Sub-critical Sobolev spaces}. Theorems~\ref{thm:laplacian_eigenmaps_estimation_fo} and~\ref{thm:laplacian_eigenmaps_estimation_ho} do not require that the smoothness index $s$ of the Sobolev space satisfy $s > d/2$, a condition often seen in the literature. In the sub-critical regime $s \leq d/2$, the Sobolev space $H^s(\mc{X})$ is quite irregular. It is not a Reproducing Kernel Hilbert Space (RKHS), nor does it continuously embed into $C^0(\mc{X})$, much less into any H\"{o}lder space. As a result, for certain versions of the nonparametric regression problem---e.g. when loss is measured in $L^{\infty}$ norm, or when the design points $\{X_1,\ldots,X_n\}$ are assumed to be fixed---in a minimax sense even consistent estimation is not possible. Likewise, certain estimators are ``off the table'', most notably RKHS-based methods such as thin-plate splines of degree $k \leq d/2$. Nevertheless, for random design regression with error measured in squared $L^2(P)$-norm, the population-level spectral series estimator $\wt{f}$ obtains the standard minimax rates $n^{-2s/(2s + d)}$ for all values of $s$ and $d$. Theorems~\ref{thm:laplacian_eigenmaps_estimation_fo} and \ref{thm:laplacian_eigenmaps_estimation_ho} show that the same is true with respect to PCR-LE, when error is measured in empirical norm.
	\item \emph{Smoothness of design density}. As promised,  Theorem~\ref{thm:laplacian_eigenmaps_estimation_ho} shows that PCR-LE achieves optimal rates of convergence so long as the unknown design density $p$ is sufficiently smooth, $p \in C^{s - 1}(\mc{X})$. The requirement $p \in C^{s - 1}(\mc{X})$ is essential to showing that $\wh{f}$ enjoys the faster minimax rates of convergence when $s > 1$,  as we discuss in Section~\ref{subsec:analysis}. 
	\item \emph{Computational considerations}. The lower bound on $\varepsilon$ in~\ref{asmp:parameters_estimation_ho} will result in a dense neighborhood graph $G$, meaning the average degree of $G$ will grow polynomially in the sample size $n$ as $n \to \infty$. As compared to a sparse $G$, this results in more non-zero entries in the graph Laplacian, and increases the computational burden involved in computing $\wh{f}$. To address this issue, in Appendix~\ref{sec:computational_considerations} we review some approaches to \emph{spectral sparsification}, in which one efficiently computes a sparse graph $\wc{G}$ that approximates $G$ in a spectral sense. The hope is that the PCR-LE estimator $\wc{f}$ computed with respect to the sparsified graph $\wc{G}$ has similar statistical properties as $\wh{f}$, while being much faster to compute. To that end, we provide upper bounds on $\|\wc{f} - f_0\|_n^2$, which show that under mild conditions on $\wc{G}$---provably achieved by many spectral sparsification algorithms---the estimator $\wc{f}$ achieves the same rates of convergence as $\wh{f}$. 
\end{itemize}

\subsubsection{Testing with PCR-LE}  The test $\varphi$ can adapt to the higher-order smoothness of $f_0$, when $\varepsilon$ and $K$ are chosen correctly.
\begin{enumerate}[label=(P\arabic*)]
	\setcounter{enumi}{3}
	\item 
	\label{asmp:parameters_testing_ho}
	The graph radius $\varepsilon$ satisfies~\eqref{eqn:radius_ho}, and the number of eigenvectors
	\begin{equation}
	\label{eqn:eigenvector_testing_ho}
	K = \min\Bigl\{\floor{(M^2n)^{2d/(4s + d)}} \vee 1, n\Bigr\}.
	\end{equation}
\end{enumerate}
When $d \leq 4$ and $n$ is sufficiently large, it is possible to choose $\varepsilon$ and $K$ such that both~\eqref{eqn:radius_ho} and~\eqref{eqn:eigenvector_testing_ho} are satisfied, and our next theorem establishes that in this situation $\varphi$ is an optimal test.
\begin{theorem}
	\label{thm:laplacian_eigenmaps_testing_ho}
	Fix $a,b \in (0,1)$. Suppose Model~\ref{def:model_flat_euclidean}. Then $\mathbb{E}_0[\varphi] \leq a$, i.e $\varphi$ is a level-$a$ test. Suppose additionally $f_0 \in H_0^s(\mc{X},M)$, that $p \in C^{s-1}(\mc{X})$, and that $d \leq 4$. Then there exist constants $c,C$ and $N$ that do not depend on $f_0$, such that the following statement holds for all $n \geq N$: if the PCR-LE test $\varphi$ is computed with a kernel $\eta$ satisfying~\ref{asmp:kernel_flat_euclidean}, and parameters $\varepsilon$ and $K$ satisfying~\ref{asmp:parameters_testing_ho}, and if $f_0$ satisfies
	\begin{equation}
	\label{eqn:laplacian_eigenmaps_testing_criticalradius_ho}
	\|f_0\|_P^2 \geq \frac{C}{b}\biggl(\Bigl(M^2(M^2n)^{-4s/(4s + d) } \wedge n^{-1/2}\Bigr)\biggl[\sqrt{\frac{1}{a}} + \frac{1}{b}\biggr] \vee \frac{M^2}{b n^{2s/d}} \biggr) \vee \frac{1}{n},
	\end{equation}
	then $\Ebb_{f_0}[1 - \varphi] \leq b$.
\end{theorem}
Similarly to the first-order case, the main takeaway from Theorem~\ref{thm:laplacian_eigenmaps_testing_ho} is that when $n^{-1/2} \lesssim M \lesssim n^{(4s - d)/4d}$, then $\varphi$ is a minimax rate-optimal test over $H_0^s(\mc{X})$. However, unlike the first-order case, when $4 < d < 4s$ the minimax testing rate over $H_0^s(\mc{X})$ is still on the order of $M^2(M^2n)^{-4s/(4s + d)}$, but we can no longer claim that $\varphi$ is an optimal test in this regime.
\begin{theorem}
	\label{thm:laplacian_eigenmaps_testing_ho_suboptimal}
	Under the same setup as Theorem~\ref{thm:laplacian_eigenmaps_estimation_ho}, but with $4 < d < 4s$. If the PCR-LE test $\varphi$ is computed with a kernel $\eta$ satisfying~\ref{asmp:kernel_flat_euclidean}, number of eigenvectors $K$ satisfying~\eqref{eqn:eigenvector_testing_ho}, and $\varepsilon = (M^2n)^{-1/(2(s - 1) + d)}$, and if 
	\begin{equation}
	\label{eqn:laplacian_eigenmaps_testing_criticalradius_ho_suboptimal}
	\|f_0\|_P^2 \geq \frac{C}{b}\biggl(\Bigl(M^2(M^2n)^{-2s/(2(s - 1) + d) } \wedge n^{-1/2}\Bigr)\biggl[\sqrt{\frac{1}{a}} + \frac{1}{b}\biggr] \vee \frac{M^2}{b n^{2s/d}} \biggr) \vee \frac{1}{n},
	\end{equation}
	then $\Ebb_{f_0}[1 - \varphi] \leq b$.
\end{theorem}
Focusing on the special case where $M \asymp 1$, Theorem~\ref{thm:laplacian_eigenmaps_testing_ho_suboptimal} says that $\varphi$ has small Type II error whenever $\|f_0\|_P^2 \gtrsim n^{-2s/(2(s - 1) + d)}$ and $4 < d < 4s$. This is smaller than the estimation rate $n^{-2s/(2s + d)}$, but larger than the minimax squared critical radius $n^{-4s/(4s + d)}$. 

At a high level, it is intuitively reasonable that PCR-LE should have more difficulty achieving the minimax rates of convergence for testing, as opposed to estimation. To obtain the faster rates of convergence for testing, PCR-LE must use many more eigenvectors than are necessary for estimation, including some eigenvectors which correspond to very large eigenvalues. It is known that the approximation properties of eigenvectors corresponding to large eigenvalues are very poor~\citep{burago2014,trillos2019}, and when $d > 4$ this prevents us from establishing that PCR-LE is an optimal test. That being said, although we suspect $\varphi$ is truly suboptimal when $d > 4$, our analysis relies on an upper bound on testing bias. Since we do not prove a matching lower bound, we cannot rule out that the test $\varphi$ is optimal for all $s < d/4$. We leave the matter to future work.

\subsection{Analysis of PCR-LE}
\label{subsec:analysis}

We now outline the high-level strategy we follow when proving each of Theorems~\ref{thm:laplacian_eigenmaps_estimation_fo}-\ref{thm:laplacian_eigenmaps_testing_ho_suboptimal}. We analyze the estimation error of $\wh{f}$, and the testing error of $\varphi$, by first conditioning on the design points $\{X_1,\ldots,X_n\}$ and deriving \emph{design-dependent} bias and variance terms. For estimation, we show that with probability at least $1 - \exp(-K)$,
\begin{equation}
\label{eqn:estimation_biasvariance}
\|\wh{f} - f_0\|_n^2 \leq \underbrace{\frac{\dotp{L_{n,\varepsilon}^s f_0}{f_0}_n}{\lambda_{K}^s}}_{\textrm{bias}} + \underbrace{\frac{5K}{n} \vphantom{\frac{\dotp{L^s f_0}{f_0}_n}{\lambda_{K}^s}}}_{\textrm{variance}}.
\end{equation}
For testing, we show that $\varphi$ (which is a level-$a$ test by construction) also has small Type II Error, $\Ebb_{f_0}[1 - \varphi] \leq b/2$, if 
\begin{equation}
\label{eqn:testing_biasvariance}
\|f_0\|_n^2 \geq  \underbrace{\frac{\dotp{L_{n,\varepsilon}^s f_0}{f_0}_n}{\lambda_{K}^s}}_{\textrm{bias}} + \underbrace{32\frac{\sqrt{2K}}{n}\biggl[\sqrt{\frac{1}{a}} + \frac{1}{b}\biggr]}_{\textrm{variance}}.
\end{equation}
These design-dependent bias-variance decompositions are reminiscent of the more classical bias-variance decompositions typical in the analysis of population-level spectral series methods (for instance~\eqref{pf:spectral_series_estimation_3} and~\eqref{pf:spectral_series_test}), but different in certain key respects. Comparing~\eqref{eqn:estimation_biasvariance} and~\eqref{eqn:testing_biasvariance} to~\eqref{pf:spectral_series_estimation_3} and~\eqref{pf:spectral_series_test}, we see that two continuum objects in the latter pair of bounds, the Sobolev norm $\|f_0\|_{\mc{H}^s(\mc{X})}^2$ and the eigenvalue $\rho_k^s$, have been replaced by graph-based analogues: the graph Sobolev seminorm $\dotp{L_{n,\varepsilon}^sf_0}{f_0}_n$ and the graph Laplacian eigenvalue $\lambda_k^s$. These latter quantities, along with the empirical squared norm $\|f_0\|_n^2$, are random variables that depend on the random design points $\{X_1,\ldots,X_n\}$. We proceed to establish suitable upper and lower bounds on these quantities that hold in probability. 

\paragraph{Estimates of graph Sobolev seminorms.}
In Proposition~\ref{prop:graph_seminorm_fo} we restate an upper bound on the first-order graph Sobolev semi-norm $\dotp{L_{n,\varepsilon}f}{f}_n$ from~\cite{green2021}. 
\begin{proposition}[Lemma 1 of~\cite{green2021}]
	\label{prop:graph_seminorm_fo}
	Suppose Model~\ref{def:model_flat_euclidean}, and additionally $f \in H^1(\mc{X})$. There exist constants $c,C$ that do not depend on $f$ or $n$ such that the following statement holds for any $\delta \in (0,1)$: if $\eta$ satisfies~\ref{asmp:kernel_flat_euclidean} and $\varepsilon < c$, then
	\begin{equation}
	\label{eqn:graph_seminorm_fo}
	\dotp{L_{n,\varepsilon}f}{f}_n \leq \frac{C}{\delta} \|f\|_{H^1(\mc{X})}^2,
	\end{equation}
	with probability at least $1 - \delta$.
\end{proposition}
Proposition~\ref{prop:graph_seminorm_fo} follows by upper bounding the expectation $\Ebb\dotp{L_{n,\varepsilon}f}{f}_n = \dotp{L_{P,\varepsilon}f}{f}_P$---where $L_{P,\varepsilon}$ is the non-local Laplacian operator defined in~\eqref{eqn:nonlocal_laplacian}---by (a constant times) the squared Sobolev norm $\|f\|_{H^1(\mc{X})}^2$, and then applying Markov's inequality.

In this work, we establish that an analogous bound holds for the graph Sobolev seminorm $\dotp{L_{n,\varepsilon}^sf}{f}_n$, when $s > 1$.
\begin{proposition}
	\label{prop:graph_seminorm_ho} 
	Suppose Model~\ref{def:model_flat_euclidean}, and additionally that $f \in H_0^s(\mc{X})$ and $p \in C^{s - 1}(\mc{X})$. Then there exist constants $c$ and $C$ that do not depend on $f_0$ or $n$ such that the following statement holds for any $\delta \in (0,1)$: if $\eta$ satisfies~\ref{asmp:kernel_flat_euclidean} and $Cn^{-1/(2(s - 1) + d)} < \varepsilon < c$, then
	\begin{equation}
	\label{eqn:graph_seminorm_ho}
	\dotp{L_{n,\varepsilon}^s f}{f}_n \leq \frac{C}{\delta} \|f\|_{H^s(\mc{X})}^2 ,
	\end{equation}
	with probability at least $1 - \delta$.
\end{proposition}
We now summarize the techniques used to prove Proposition~\ref{prop:graph_seminorm_ho}, which will help explain the role played by our conditions on $f_0$, $p$ and $\varepsilon$. To upper bound $\dotp{L_{n,\varepsilon}^sf}{f}_n$ in terms of $\|f\|_{H^s(\mc{X})}^2$, we introduce an intermediate quantity: the \emph{non-local Sobolev seminorm} $\dotp{L_{P,\varepsilon}^sf}{f}_{P}$. This seminorm is defined with respect to the iterated non-local Laplacian $L_{P,\varepsilon}^s = L_{P,\varepsilon} \circ \cdots \circ L_{P,\varepsilon}$, where $L_{P,\varepsilon}$ is a non-local approximation to $\Delta_P$, 
\begin{equation}
\label{eqn:nonlocal_laplacian}
L_{P,\varepsilon}f(x) := \frac{1}{\varepsilon^{d + 2}}\int_{\mc{X}}\bigl(f(z) - f(x)\bigr) \eta\biggl(\frac{\|z - x\|}{\varepsilon}\biggr) \,dP(x).
\end{equation}
Then the proof of Proposition~\ref{prop:graph_seminorm_ho} proceeds according to the following steps.
\begin{enumerate}
	\item \emph{Bound on pure bias terms.} First we note that $\dotp{L_{n,\varepsilon}^s f}{f}_n$ is itself a biased estimate of the non-local seminorm $\dotp{L_{P,\varepsilon}^sf}{f}_{P}$. This is because $\dotp{L_{n,\varepsilon}^s f}{f}_n$ is a $V$-statistic, meaning it is the sum of an unbiased estimator of $\dotp{L_{P,\varepsilon}^sf}{f}_{P}$ (in other words, a $U$-statistic) plus some higher-order, pure bias terms. We show that these pure bias terms are negligible when $\varepsilon = \omega(n^{-1/(2(s - 1) + d)})$. 
	\item \emph{Convergence in the interior.} For $x$ sufficiently far from the boundary of $\mc{X}$---precisely $x \in \mc{X}$ such that $B(x,j\varepsilon) \subseteq \mc{X}$---we show that $L_{P,\varepsilon}^jf(x) \to \sigma_{\eta}^j \Delta_P^jf(x)$ as $\varepsilon \to 0$. Here $j = (s - 1)/2$ when $s$ is odd and $j = (s - 2)/2$ when $s$ is even. This step bears some resemblance to the analysis of the bias term in kernel smoothing, and requires that $p \in C^{s-1}(\mc{X})$.
	\item \emph{Boundedness at the boundary.} On the other hand for $x$ sufficiently near the boundary of $\mc{X}$, $L_{P,\varepsilon}^jf(x)$ does not in general converge to $\sigma_{\eta}^j\Delta_P^jf(x)$. Instead, we use the zero-trace property of $f$ to show that $L_{P,\varepsilon}^jf(x)$ is small.
	\item \emph{Putting together the pieces.} Finally, we combine the results of the previous two steps to deduce an upper bound on $\dotp{L_{P,\varepsilon}^sf}{f}_{P}$ in terms of the squared Sobolev norm $\|f\|_{H^s(\mc{X})}^2$. The nature of this last step depends on whether $s$ is an even or an odd integer. Roughly speaking, when $s$ is odd, letting $j = (s - 1)/2$, we show that  
	\begin{equation*}\dotp{L_{P,\varepsilon}^sf}{f}_P = \dotp{L_{P,\varepsilon}L_{P,\varepsilon}^jf}{L_{P,\varepsilon}^jf}_P \approx \sigma_{\eta}^{2j}\dotp{L_{P,\varepsilon}\Delta_P^jf}{\Delta_P^j}_{P} \lesssim \sigma_{\eta}^{2j + 1} \dotp{\Delta_P^jf}{f}_{P}.
	\end{equation*}
	When $s$ is even, letting $j = (s - 2)/2$, we show that
	\begin{equation*}
	\dotp{L_{P,\varepsilon}^sf}{f}_P = \|L_{P,\varepsilon}L_{P,\varepsilon}^{j}f\|_{P}^2 \approx \sigma_{\eta}^{2j}\|L_{P,\varepsilon} \Delta_P^jf\|_P^2 \lesssim \|\Delta_P^{j + 1}f\|_P^2.
	\end{equation*}
	In either case, the desired upper bound~\eqref{eqn:graph_seminorm_ho} follows from the boundedness of the density $p$.
\end{enumerate}
It is worth pointing out that we do not need to establish the pointwise estimate $L_{P,\varepsilon}^sf \to \sigma_{\eta}^s\Delta_{P}^sf$ in $L^2(P)$ norm. If we had such an estimate, it would immediately follow that $\dotp{L_{P,\varepsilon}^sf}{f}_{P} \to \sigma_{\eta}^s\dotp{\Delta_P^sf}{f}_{P}$. Unfortunately, we assume only that $f$ has $s$ bounded derivatives, while $\Delta_P^s$ is an order-$2s$ differential operator; thus in general $L_{P,\varepsilon}^sf$ may not approach $\sigma_{\eta}^s\Delta_{P}^sf$ as $\varepsilon \to 0$. Instead we opt for the slightly more complicated approach outlined above, in which we only ever need show that $L_{P,\varepsilon}^jf(x) \to \sigma_{\eta}^j \Delta_P^jf(x)$ for some $j < s/2$. 

\paragraph{Neighborhood graph eigenvalues.}
On the other hand, several recent works \citep{burago2014,garciatrillos18,calder2019} have analyzed the convergence of $\lambda_{k}$ towards $\rho_k$. They provide explicit bounds on the relative error $|\lambda_{k} - \rho_k|/\rho_k$, which show that the relative error is small for sufficiently large $n$ and small $\varepsilon$. Crucially, these guarantees hold simultaneously for all $1 \leq k \leq K$ as long as $\\rho_K = O(\varepsilon^{-2})$. These results are actually stronger than are necessary to establish Theorems~\ref{thm:laplacian_eigenmaps_estimation_fo}-\ref{thm:laplacian_eigenmaps_testing_ho}---in order to get rate-optimality, we need only show that for the relevant values of $K$, $\lambda_{K}/\lambda_K(P) = \Omega_P(1)$---but unfortunately they all assume $P$ is supported on a manifold without boundary (i.e. they assume Model~\ref{def:model_manifold} rather than Model~\ref{def:model_flat_euclidean}). 

In the case where $\mc{X}$ is assumed to have a boundary, the graph Laplacian $L_{n,\varepsilon}$ is a reasonable approximation of the operator $\Delta_P$ only at points $x \in \mc{X}$ for which $B(x,\varepsilon) \subseteq \mc{X}$. In contrast, at points $x$ near the boundary of $\mc{X}$, the graph Laplacian is known to approximate a different operator altogether \citep{belkin2012}.\footnote{This is directly related to the boundary bias of kernel smoothing, since the graph Laplacian can be viewed as a kernel-based estimator of $\Delta_P$.} This renders analysis of $\lambda_k$ substantially more challenging, since its continuum limit is not $\rho_k$.  Rather than analyzing the convergence of $\lambda_k$, we will instead use Lemma~2 of \cite{green2021}, whose assumptions match our own, and who give a weaker bound on the ratio $\lambda_k/\rho_k$ that will nevertheless suffice for our purposes. 

\begin{proposition}[Lemma~2 of \cite{green2021}]
	\label{prop:graph_eigenvalue}
	Suppose Model~\ref{def:model_flat_euclidean}. Then there exist constants $c$ and $C$ such that the following statement holds: if $\eta$ satisfies~\ref{asmp:kernel_flat_euclidean} and $C(\log n/n)^{1/d} < \varepsilon < c$, then
	\begin{equation}
	\label{eqn:graph_eigenvalue}
	\lambda_k \geq c \cdot \min\Bigl\{\rho_k, \frac{1}{\varepsilon^{2}} \Bigr\} \quad \textrm{for all $1 \leq k \leq n$,}
	\end{equation}
	with probability at least $1 - Cn\exp\{-c n\varepsilon^d\}$. 
\end{proposition}
By our assumptions on $P$, $\lambda_0(\Delta_P) = \lambda_0 = 0$. Furthermore, Weyl's Law~\eqref{eqn:weyl} tells us that under Model~\ref{def:model_flat_euclidean}, $k^{2/d} \lesssim \rho_k \lesssim k^{2/d}$ for all $k \in \mathbb{N}, k > 1$. Combining these statements with~\eqref{eqn:graph_eigenvalue}, we conclude that $\lambda_{K} = \Omega_P(K^{2/d})$ so long as $K \lesssim \varepsilon^{-d}$. 

\paragraph{Empirical norm.}
Finally, in Proposition~\ref{prop:empirical_norm_sobolev} we establish that a one-sided bound of the form $\|f_0\|_n^2 \gtrsim \|f_0\|_P^2$ whenever $\|f_0\|_P^2$ is itself sufficiently large.
\begin{proposition}
	\label{prop:empirical_norm_sobolev}
	Suppose Model~\ref{def:model_flat_euclidean}, and additionally that $f \in H^s(\mc{X},M)$ for some $s > d/4$. There exist constants $c$ and $C$ that do not depend on $f_0$ or $n$ such that the following statement holds for any $\delta > 0$:  if
	\begin{equation}
	\label{eqn:empirical_norm_sobolev_1}
	\|f\|_{P} \geq C M \biggl(\frac{1}{\delta n}\biggr)^{s/d}
	\end{equation}
	then with probability at least $1 - \exp\{-(cn \wedge 1/\delta)\}$,
	\begin{equation}
	\label{eqn:empirical_norm_sobolev}
	\|f\|_n^2 \geq \frac{1}{2} \|f_0\|_P^2.
	\end{equation}
\end{proposition}
To prove Proposition~\ref{prop:empirical_norm_sobolev}, we use a Gagliardo-Nirenberg interpolation inequality (see e.g. Theorem~12.83 of~\citep{leoni2017}) to control the $4$th moment of $f \in H^s(\mc{X})$ in terms of $\|f\|_P$ and $|f|_{H^s(\mc{X})}$, then invoke a one-sided Bernstein's inequality as in \cite[Section 14.2]{wainwright2019}. Note carefully that the statement~\eqref{eqn:empirical_norm_sobolev} is \emph{not} a uniform guarantee over all $f \in H^s(\mc{X};M)$. Indeed, such a statement cannot hold in the sub-critical regime ($2s \leq d$).\footnote{This is because in the sub-critical regime, for any set of points $\{x_1,\ldots,x_n\}$ there exists a sequence of functions $\{f_k:k \in \mathbb{N}\} \subset H^s(\mc{X};1)$ satisfying $f_k(x_i) = 1$ for each $i = 1,\ldots,n$---and therefore $\|f_k\|_n^2 = 1$---but for which $\|f_k\|_P^2 \to 0$ as $k \to \infty$.} Fortunately, a pointwise bound---meaning a bound that holds with high probability for a single $f \in H^s(\mc{X})$---is sufficient for our purposes.

Finally, invoking the bounds of Propositions~\ref{prop:graph_seminorm_fo}-\ref{prop:empirical_norm_sobolev} inside the bias-variance tradeoffs~\eqref{eqn:estimation_biasvariance} and~\eqref{eqn:testing_biasvariance} and then choosing $K$ to balance bias and variance (when possible), leads to the conclusions of Theorems~\ref{thm:laplacian_eigenmaps_estimation_fo}-\ref{thm:laplacian_eigenmaps_testing_ho_suboptimal}.


\section{Manifold Adaptivity}
\label{sec:manifold_adaptivity}

In this section we consider the manifold setting, where $(X_1,Y_1),\ldots,(X_n,Y_n)$ are observed according to Model~\ref{def:model_manifold}. In this setting, it is known that the minimax rates depend only on the intrinsic dimension $m$; more specifically,~\cite{bickel2007,ariascastro2018} show that for functions with H\"{o}lder smoothness $s$, the minimax estimation rate is $n^{-2s/(2s + m)}$ and the testing rate is $n^{-4s/(4s + m)}$.\footnote{Although~\cite{ariascastro2018} considers density testing, usual arguments regarding equivalence of experiments~\citep{brown96} imply that the same rates apply to regression testing.} On the other hand, a theory has been developed \citep{niyogi2008finding,belkin03,belkin05,niyogi2013,balakrishnan2012minimax,balakrishnan2013cluster} establishing that the neighborhood graph $G$ can ``learn'' the manifold $\mc{X}$ in various senses, so long as $\mc{X}$ is locally linear.  We build on this work by showing that when $P$ is supported on a manifold $\mc{X}$ and $f_0 \in H^s(\mc{X})$, PCR-LE achieves the sharper minimax estimation and testing rates.

\subsection{Upper bounds}
Unlike in the flat-Euclidean case, since Model~\ref{def:model_manifold} assumes that $\mc{X}$ is without boundary it is easy to deal with the first-order $(s = 1)$ and higher-order $(s > 1)$ cases all at once. A more important distinction between the results of this section and those of Section~\ref{sec:minimax_optimal_laplacian_eigenmaps} is that we will establish PCR-LE is optimal only when the regression function $f_0 \in H^s(\mc{X};M)$ for $s \in \{1,2,3\}$. Otherwise, this section will proceed in a similar fashion to Section~\ref{sec:higher_order_sobolev_classes}.

\subsubsection{Estimation with PCR-LE} 
To ensure that $\wh{f}$ is an in-sample minimax rate-optimal estimator, we choose the kernel function $\eta$, graph radius $\varepsilon$ and number of eigenvectors $K$ as in~\ref{asmp:parameters_estimation_ho}, except with ambient dimension $d$ replaced by the intrinsic dimension $m$.

\begin{enumerate}[label=(P\arabic*)]
	\setcounter{enumi}{4}
	\item 
	\label{asmp:kernel_manifold}
	The kernel function $\eta$ is a nonincreasing function supported on a subset of $[0,1]$. Its restriction to $[0,1]$ is Lipschitz, and $\eta(1/2) > 0$. Additionally, it is normalized so that
	\begin{equation*}
	\int_{\Reals^m} \eta(\|z\|) \,dz = 1,
	\end{equation*}
	and we assume \smash{$\int_{\Reals^m} \|x\|^2 \eta(\|x\|) \,dx < \infty$}.
	\item 
	\label{asmp:parameters_estimation_manifold}
	For constants $c_0,C_0$, the graph radius $\varepsilon$ and number of eigenvectors $K$ satisfy
	\begin{equation}
	\label{eqn:radius_estimation_manifold}
	C_0\max\biggl\{\biggl(\frac{\log}{n}\biggr)^{1/m}, n^{-1/(2(s - 1) + m)}\biggr\} \leq \varepsilon \leq c_0\min\{1, K^{-1/m}\}.
	\end{equation}
	Additionally,
	\begin{equation*}
	K = \min\Bigl\{\floor{(M^2n)^{m/(2s + m)}} \wedge 1,n \Bigr\}.
	\end{equation*}
\end{enumerate}

\begin{theorem}
	\label{thm:laplacian_eigenmaps_estimation_manifold}
	Suppose Model~\ref{def:model_manifold}, and additionally $f_0 \in H^s(\mc{X},M)$ and $p \in C^{s - 1}(\mc{X})$ for $s \leq 3$. There exist constants $c,C$ and $N$ that do not depend on $f_0$, such that the following statement holds all for all $n$ larger than $N$ and for any $\delta \in (0,1)$: if the PCR-LE estimator $\wh{f}$ is computed with a kernel $\eta$ satisfying~\ref{asmp:kernel_manifold}, and parameters $\varepsilon$ and $K$ satisfying~\ref{asmp:parameters_estimation_manifold}, then
	\begin{equation}
	\label{eqn:laplacian_eigenmaps_estimation_manifold}
	\|\wh{f} - f_0\|_n^2 \leq C\Bigl(\frac{1}{\delta}M^2(M^2n)^{-2s/(2s + m)} \wedge 1\Bigr) \vee \frac{1}{n},
	\end{equation}
	with probability at least $1 - \delta - Cn\exp(-cn\varepsilon^m) - \exp(-K)$.
\end{theorem}

\subsubsection{Testing with PCR-LE}
Likewise, to construct a minimax optimal test using $\wh{T}$, we choose $\varepsilon$ and $K$ as in~\ref{asmp:parameters_testing_fo}, except with the ambient dimension $d$ replaced by the intrinsic dimension $m$.
\begin{enumerate}[label=(P\arabic*)]
	\setcounter{enumi}{5}
	\item 
	\label{asmp:parameters_testing_manifold}
	The graph radius $\varepsilon$ satisfies~\eqref{eqn:radius_estimation_manifold}, and the number of eigenvectors
	\begin{equation*}
	K = \min\Bigl\{\floor{(M^2n)^{2m/(4s + m)}} \wedge 1,n \Bigr\}.
	\end{equation*}
\end{enumerate}

\begin{theorem}
	\label{thm:laplacian_eigenmaps_testing_manifold}
	Fix $a,b \in (0,1)$. Suppose Model~\ref{def:model_manifold}. Then $\mathbb{E}_0[\varphi] \leq a$, i.e $\varphi$ is a level-$a$ test. Suppose additionally $f_0 \in H^s(\mc{X},M)$, that $p \in C^{s-1}(\mc{X})$, and that $s \leq 3$ and $m \leq 4$. Then there exist constants $c$, $C$ and $N$ that do not depend on $f_0$, such that the following statement holds for all $n$ larger than $N$: if the PCR-LE test $\varphi$ is computed with a kernel $\eta$ satisfying~\ref{asmp:kernel_manifold}, and parameters $\varepsilon$ and $K$ satisfying~\ref{asmp:parameters_testing_manifold}, and if $f_0$ satisfies
	\begin{equation}
	\label{eqn:laplacian_eigenmaps_testing_criticalradius_manifold}
	\|f_0\|_P^2 \geq \frac{C}{b}\biggl(\Bigl(M^2(M^2n)^{-4s/(4s + m) } \wedge n^{-1/2}\Bigr)\biggl[\sqrt{\frac{1}{a}} + \frac{1}{b}\biggr] \vee \frac{M^2}{b n^{2s/m}} \biggr) \vee \frac{1}{n},
	\end{equation}
	then $\Ebb_{f_0}[1 - \varphi] \leq b$.
\end{theorem}
Focusing on the case $M \asymp 1$,\footnote{To the best of our knowledge, the minimax rates for general $M$ in the manifold setting have not been worked out.} the upper bounds in Theorems~\ref{thm:laplacian_eigenmaps_estimation_manifold} and~\ref{thm:laplacian_eigenmaps_testing_manifold} imply that PCR-LE attain the optimal rates of convergence over Sobolev balls $H^s(\mc{X})$ for $s \in \{1,2,3\}$. 

Unlike in the full-dimensional case, in the manifold setting our upper bounds on the estimation and testing error of PCR-LE do not match the minimax rate when $s \geq 4$.  In this case, the containment $H^s(\mc{X};1) \subset H^{3}(\mc{X};1)$ implies that the PCR-LE estimator $\wh{f}$ has in-sample mean-squared error of at most on the order of $n^{-6/(6 + m)}$, and that the PCR-LE test has small Type II error whenever $\|f_0\|_P^2 \gtrsim n^{-12/(12 + m)}$; however, these are slower than the minimax rates. 

We now explain this difference between the flat Euclidean and manifold settings. At a high level, thinking of the graph $G$ as an estimate of the manifold $\mc{X}$, we incur some error by using Euclidean distance rather than geodesic distance to form the edges of $G$. This is in contrast with the full-dimensional setting, where the Euclidean metric exactly coincides with the geodesic distance for all points $x,z \in \mc{X}$ that are sufficiently close to each other and far from the boundary of $\mc{X}$. This extra error incurred in the manifold setting by using the ``wrong distance'' dominates when $s \geq 4$. 

As this explanation suggests, by building $G$ using the geodesic distance one could avoid this error, and might obtain superior rates of convergence. However this is not an option for us, as we assume $\mc{X}$---and in particular its geodesics---are unknown. Likewise, a population-level spectral series estimator using eigenfunctions of the manifold Laplace-Beltrami operator, will achieve the minimax rate for all values of $s$ and $m$; but this is undesirable for the same reason---we do not want to assume that $\mc{X}$ is known. It is not clear whether this gap between population-level spectral series regression and the PCR-LE estimator is real, or a product of loose upper bounds. 

Finally, as in the full-dimensional case, when the intrinsic dimension $m > 4$ we cannot choose the graph radius $\varepsilon$ and number of eigenvectors $K$ to optimally balance bias and variance.  Instead, reasoning as in the proof of Theorem~\ref{thm:laplacian_eigenmaps_testing_ho_suboptimal} shows that when $1 \leq s \leq 3$, the PCR-LE test has critical radius as given by~\eqref{eqn:laplacian_eigenmaps_testing_criticalradius_ho_suboptimal}, but with the ambient dimension $d$ replaced by $m$.

\subsection{Analysis}
The high-level strategy used to prove Theorems~\ref{thm:laplacian_eigenmaps_estimation_manifold} and~\ref{thm:laplacian_eigenmaps_testing_manifold} is the same as in the flat-Euclidean setting. More specifically, we will use precisely the same bias-variance decompositions~\eqref{eqn:estimation_biasvariance} (for estimation) and~\eqref{eqn:testing_biasvariance} (for testing). The difference will be that our bounds on the graph Sobolev seminorm $\dotp{L_{n,\varepsilon}^sf_0}{f_0}_n$, graph eigenvalue $\lambda_K$, and empirical norm $\|f_0\|_n^2$ will now always depend on the intrinsic dimension $m$, rather than the ambient dimension $d$. The precise results we use are contained in Propositions~\ref{prop:graph_seminorm_manifold}-\ref{prop:empirical_norm_sobolev_manifold}.
\begin{proposition}
	\label{prop:graph_seminorm_manifold} 
	Suppose Model~\ref{def:model_manifold}, and additionally that $f_0 \in H^s(\mc{X};M)$ and $p \in C^{s - 1}(\mc{X})$ for $s = 1,2$ or $3$. Then there exist constants $c_0,C_0$ and $C$ that do not depend on $f_0$, $n$ or $M$ such that the following statement holds for any $\delta \in (0,1)$: if $\eta$ satisfies~\ref{asmp:kernel_manifold} and $C_0n^{-1/(2(s - 1) + m)} < \varepsilon < c_0$, then
	\begin{equation}
	\label{eqn:graph_seminorm_manifold}
	\dotp{L_{n,\varepsilon}^s f}{f}_n \leq \frac{C}{\delta} \|f\|_{H^s(\mc{X})}^2,
	\end{equation}
	with probability at least $1 - 2\delta$.
\end{proposition}

As discussed previously, when $\mc{X}$ is a domain without boundary and $\Delta_P$ is the manifold weighted Laplace-Beltrami operator, appropriate bounds on the graph eigenvalues $\lambda_k$ have already been derived in \citep{burago2014,trillos2019,garciatrillos19}. The precise result we need is a direct consequence of Theorem 2.4 of~\citep{calder2019}.
\begin{proposition}[\textbf{c.f Theorem 2.4 of~\citep{calder2019}}]
	\label{prop:graph_eigenvalue_manifold}
	Suppose Model~\ref{def:model_manifold}. Then there exist constants $c$ and $C$ such that the following statement holds: if $\eta$ satisfies~\ref{asmp:kernel_manifold} and $C(\log n/n)^{1/m} < \varepsilon < c$, then
	\begin{equation}
	\label{eqn:graph_eigenvalue_manifold}
	\lambda_k \geq c \cdot \min\Bigl\{\rho_k, \frac{1}{\varepsilon^{2}} \Bigr\} \quad \textrm{for all $1 \leq k \leq n$,}
	\end{equation}
	with probability at least $1 - Cn\exp\{-c n\varepsilon^d\}$. 
\end{proposition}
(For the specific computation used to deduce Proposition~\ref{prop:graph_eigenvalue_manifold} from Theorem 2.4 of~\citep{calder2019}, see~\cite{green2021}.)

Finally, we have the following lower bound on the empirical norm $\|f\|_n$ under the hypotheses of Model~\ref{def:model_manifold}. 
\begin{proposition}
	\label{prop:empirical_norm_sobolev_manifold}
	Suppose Model~\ref{def:model_manifold}, and additionally that $f_0 \in H^s(\mc{X},M)$ for some $s > m/4$. There exists a constant $C$ that does not depend on $f_0$ such that the following statement holds for all $\delta > 0$:  if
	\begin{equation}
	\label{eqn:empirical_norm_sobolev_manifold_1}
	\|f_0\|_{P} \geq \frac{C M}{\delta^{s/m}}n^{-s/m},
	\end{equation}
	then with probability at least $1 - \exp\{-(cn \wedge 1/\delta)\}$,
	\begin{equation}
	\label{eqn:empirical_norm_sobolev_manifold}
	\|f_0\|_n^2 \geq \frac{1}{2} \|f_0\|_P^2.
	\end{equation}
\end{proposition}
We prove Proposition~\ref{prop:empirical_norm_sobolev_manifold} in a parallel manner to its flat Euclidean counterpart (Proposition~\ref{prop:empirical_norm_sobolev}), by first using a Gagliardo-Nirenberg inequality to upper bound the $L^4(\mc{X})$ norm of a Sobolev function defined on a compact Riemannian manifold, and then applying a one-sided Bernstein's inequality. Finally, combining Propositions~\ref{prop:graph_seminorm_manifold}-\ref{prop:empirical_norm_sobolev_manifold} with the conditional-on-design bias-variance decompositions~\eqref{eqn:estimation_biasvariance} and \eqref{eqn:testing_biasvariance} leads to the conclusions of Theorems~\ref{thm:laplacian_eigenmaps_estimation_manifold} and~\ref{thm:laplacian_eigenmaps_testing_manifold}. 

\section{Experiments}
\label{sec:experiments}

In this section we empirically demonstrate that the PCR-LE estimator and test are reasonably good alternatives to population-level spectral series methods, even at moderate sample sizes $n$. In order to compare the two approaches, in our experiments we stick to the simple case where the design distribution $P$ is the uniform distribution over $\mc{X} = [-1,1]^d$, and we have simple closed-form expressions for the eigenfunctions of $\Delta_P$. In general, it is not easy to analytically compute these eigenfunctions, which is part of the appeal of LE and PCR-LE.

\begin{figure*}[b]
	\includegraphics[width=.245\textwidth]{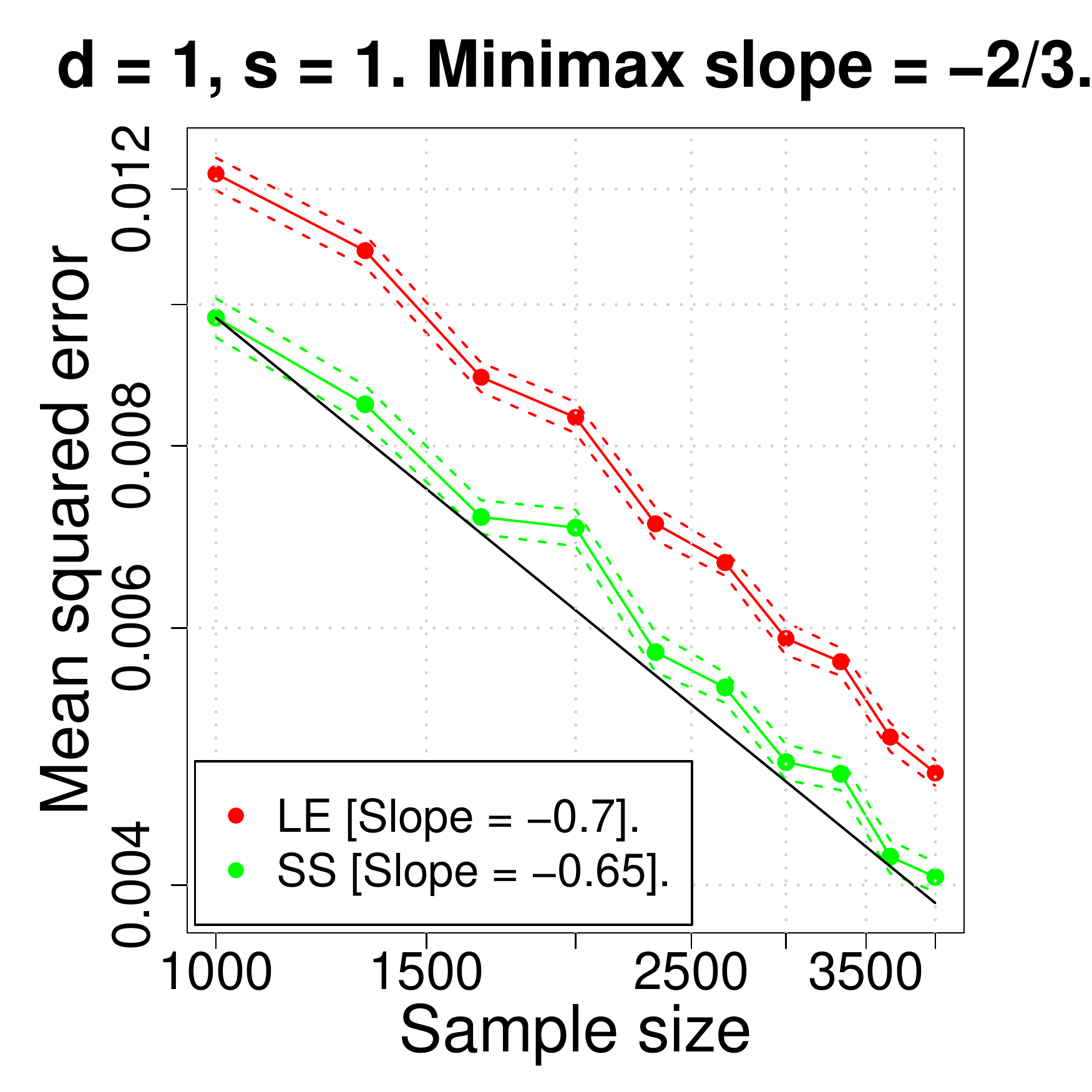}
	\includegraphics[width=.245\textwidth]{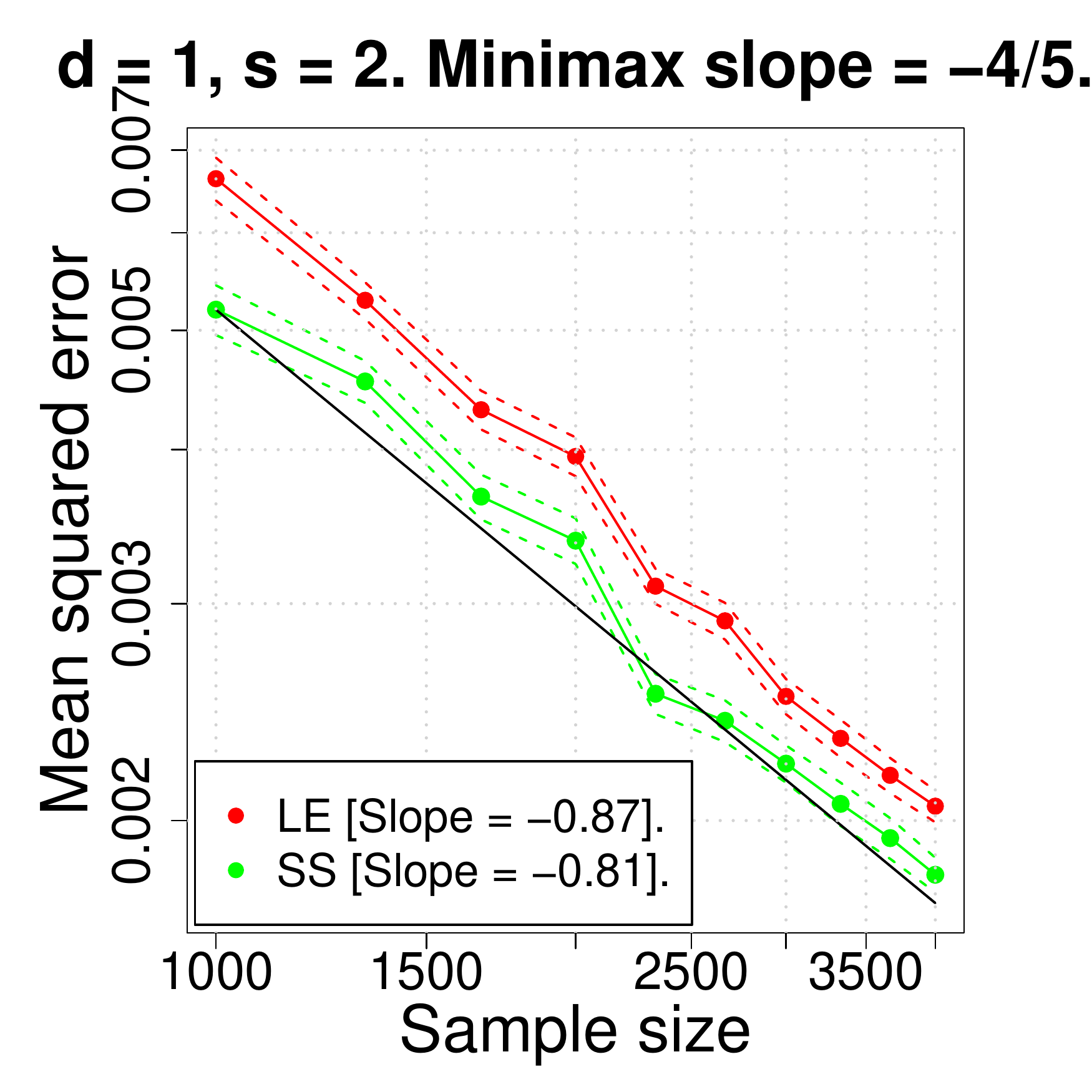}
	\includegraphics[width=.245\textwidth]{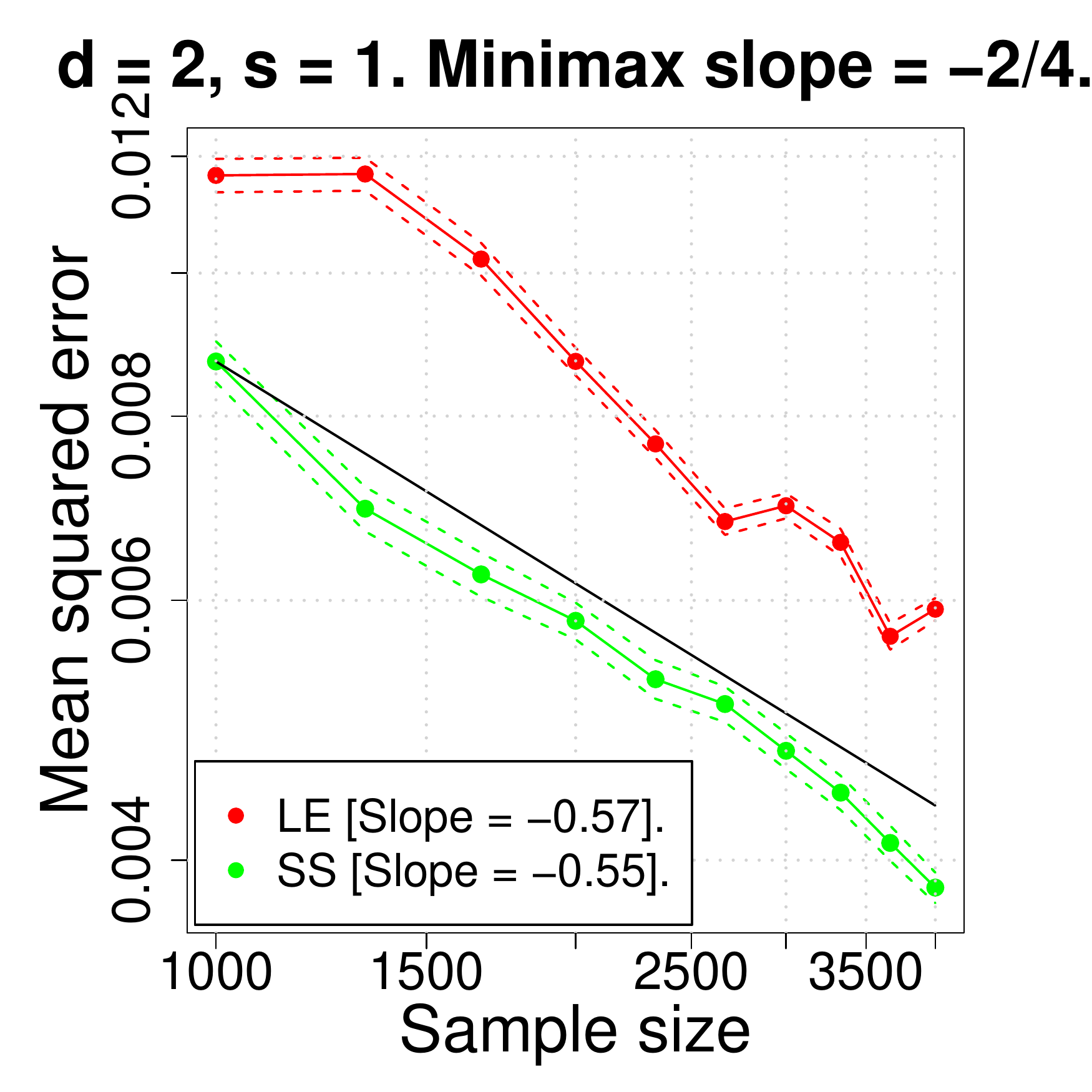}
	\includegraphics[width=.245\textwidth]{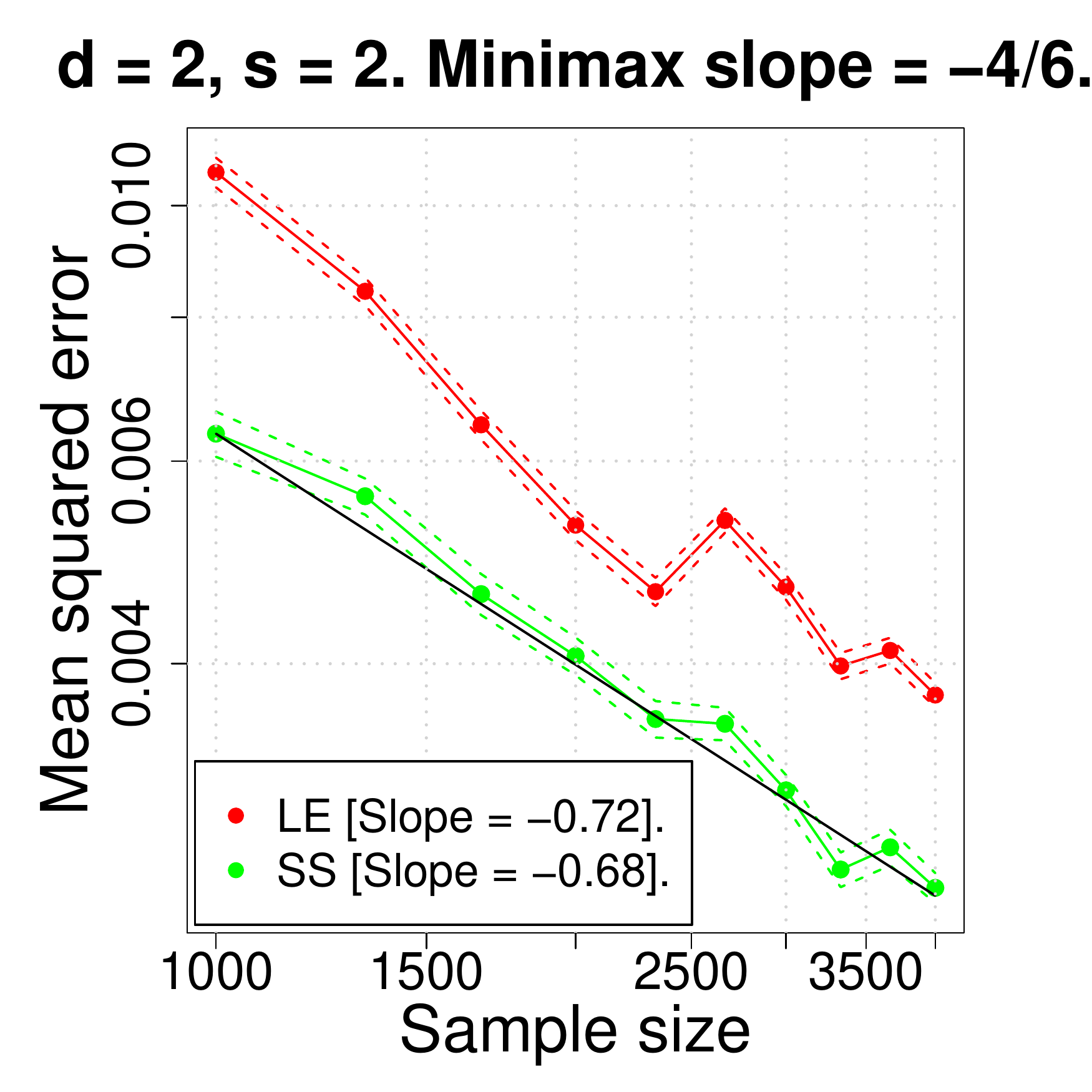}
	\caption{In-sample mean squared error (mse) of PCR-LE (\texttt{LE}) vs. population-level spectral series (\texttt{SS}) estimator, as a function of sample size $n$. Each plot is on the log-log scale, and the results are averaged over 400 repetitions. All estimators are tuned for optimal average mse, separately at each value of $n$. The black line shows the minimax rate (in slope only; the intercept is chosen to match the observed error).}
	\label{fig:fig1}
\end{figure*}

\paragraph{Estimation.}
In our first experiment, we compare the mean-squared error of the PCR-LE estimator $\wh{f}$  to that of its population-level counterpart $\wt{f}$. We vary the sample size from $n = 1000$ to $n = 4000$; sample $n$ design points $\{X_1,\ldots,X_n\}$ from the uniform distribution on the cube $[-1,1]^d$; and sample responses $Y_i$ according to~\eqref{eqn:model} with regression function $f_0 = M/\rho_K^{s/2} \cdot \psi_K$ for $K \asymp n^{d/(2s + d)}$ (the pre-factor $M/\rho_K^{s/2}$ is chosen so that $|f_0|_{H^s(\mc{X})}^2 = M^2$). In Figure~\ref{fig:fig1} we show the in-sample mean-squared error of the two estimators as a function of $n$, for different dimensions $d$ and order of smoothness $s$. We see that both estimators have mean-squared error converging to zero at roughly the minimax rate. While the unsurprisingly population-level spectral series estimator has the smaller error, generally speaking the error of PCR-LE approaches that of the population-level spectral series method as $n$ gets larger. 

\paragraph{Testing.} 
In our second experiment, we compare the PCR-LE test $\varphi$ against the  population-level spectral series test $\wt{\varphi}$. The setup is generally the same as that of our first experiment, but to get an empirical estimate of the critical radius the details are necessarily somewhat more complicated. First we take $\mc{F} = \{M/\rho_k^{s/2} \psi_k\}_{k = 1}^{n}$ to be a discrete subset of $H^1(\mc{X};M)$. Then, for each $f_0 \in \mc{F}$, we run a given test $\phi$ (either the PCR-LE test $\phi = \varphi$, or the population-level spectral series test $\phi = \wt{\varphi}$) and record whether it was a false negative or true positive. We repeat this process over $100$ replications, giving a Monte Carlo estimate of the type II error $E_{f_0}[1 - \phi]$ for each $f_0 \in \mc{F}$. Finally, we take the smallest value of $\|f_0\|_P^2$ such $E_{f_0}[1 - \phi] \leq b$ as our estimate of the critical radius of $\phi$. 

\begin{figure*}[b]
	\includegraphics[width=.245\textwidth]{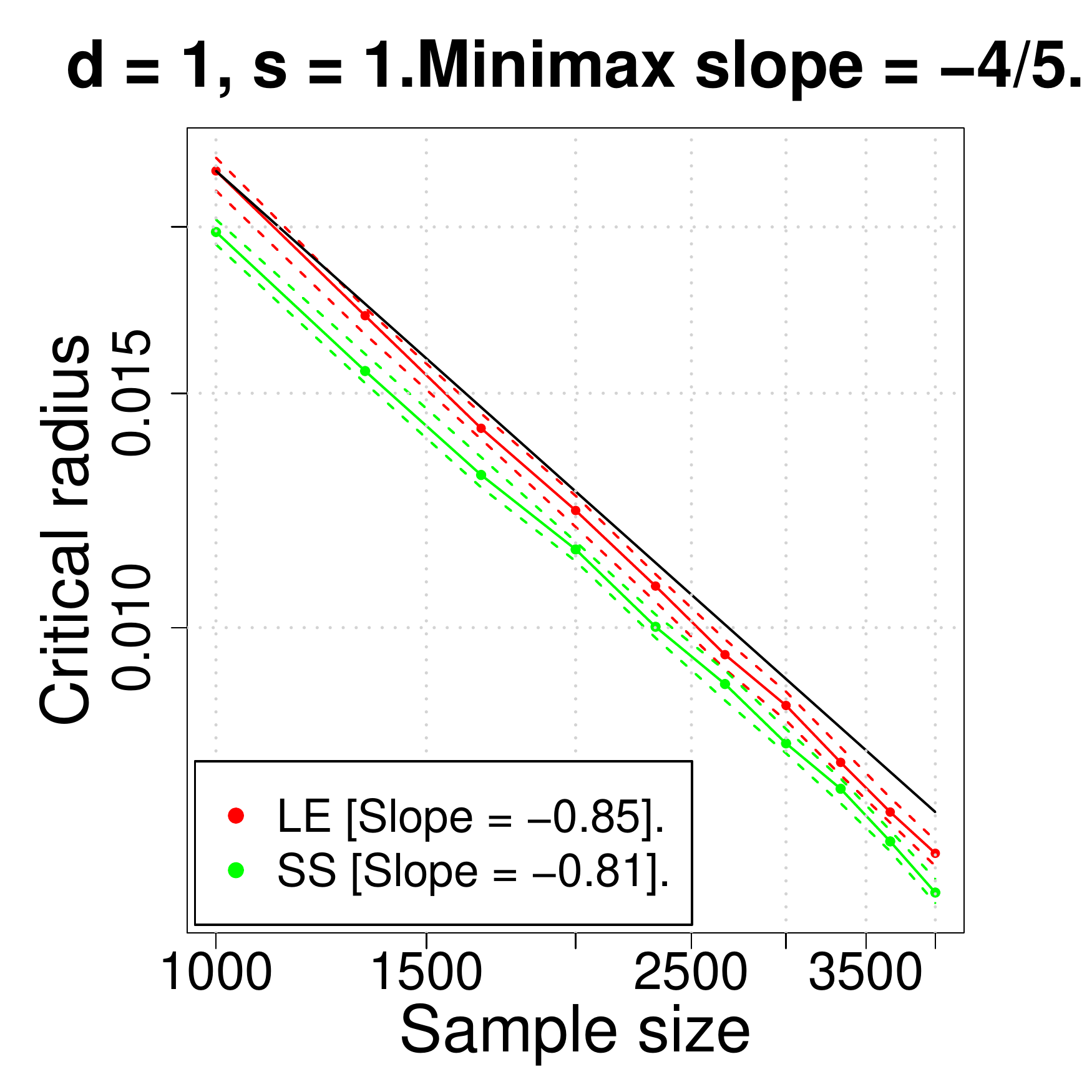}
	\includegraphics[width=.245\textwidth]{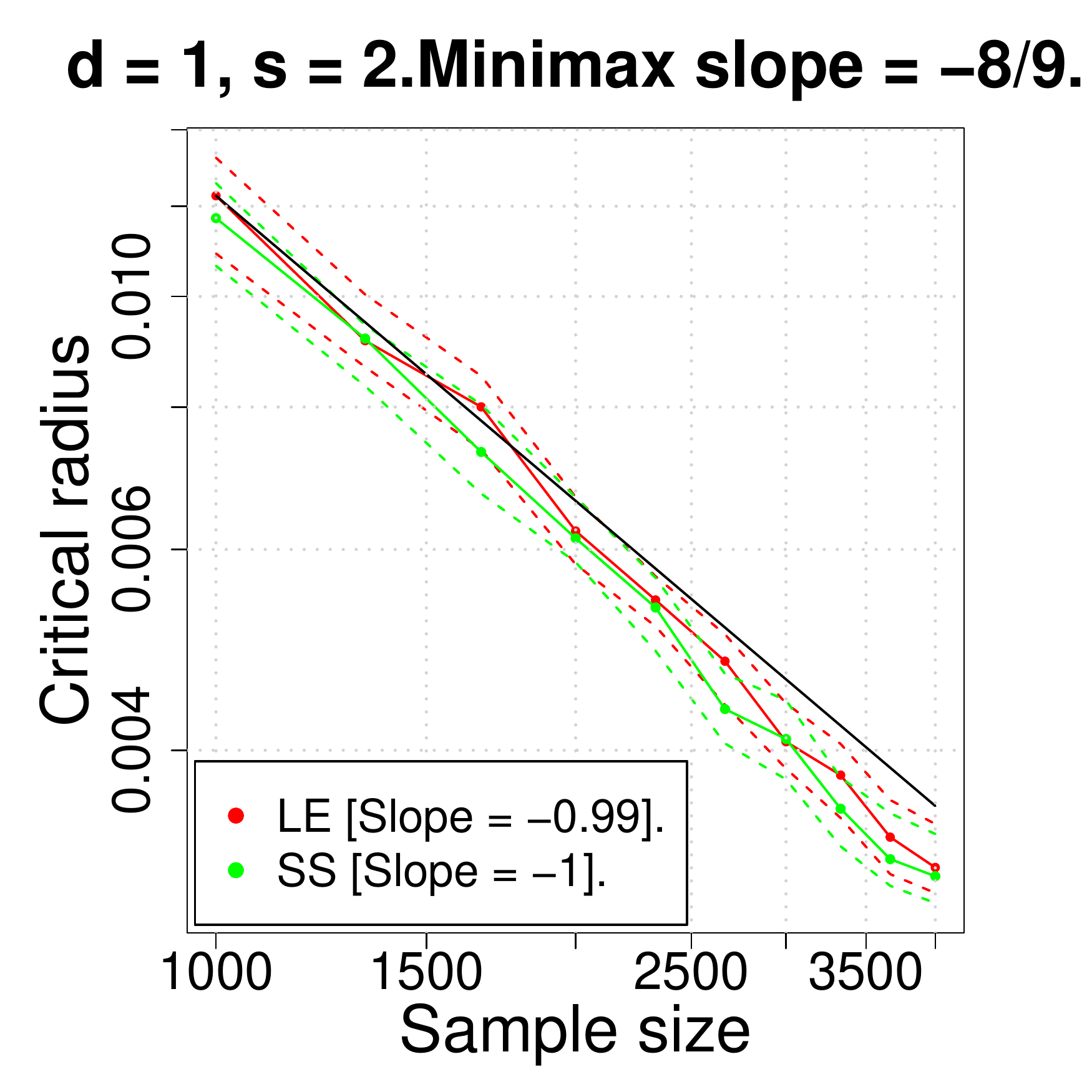}
	\includegraphics[width=.245\textwidth]{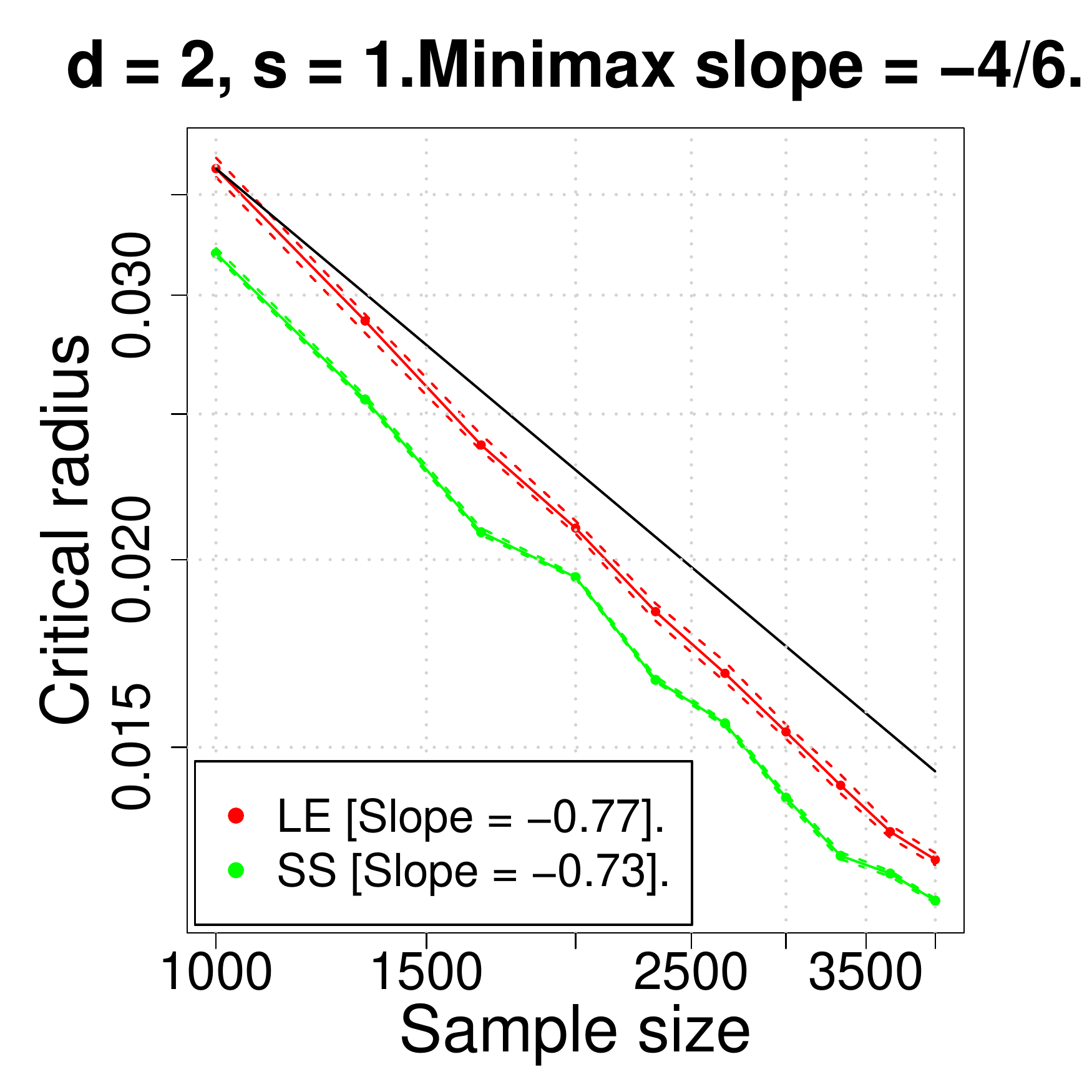}
	\includegraphics[width=.245\textwidth]{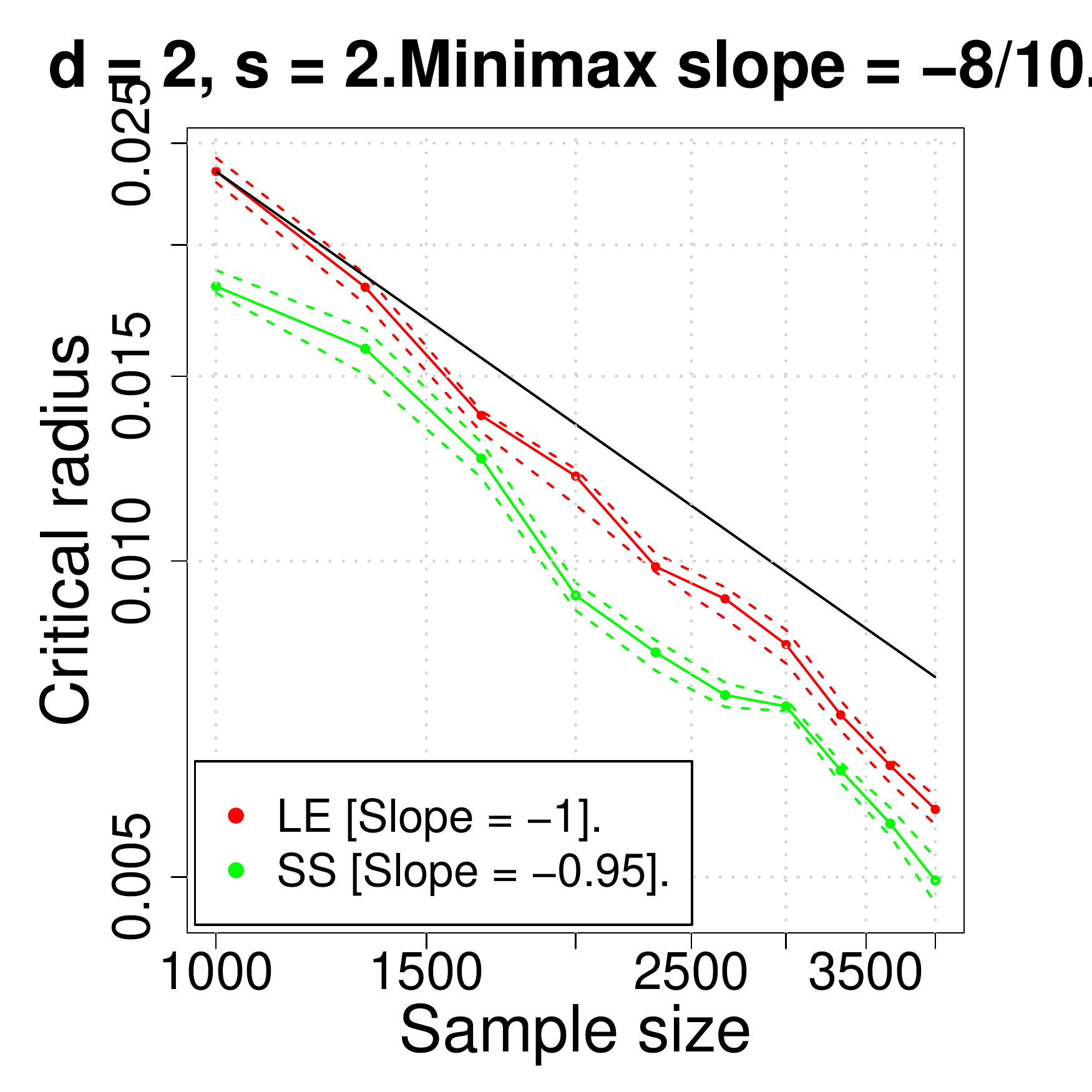}
	\caption{Worst-case testing risk for PCR-LE (\texttt{LE}) and spectral series (\texttt{SP}) tests, as a function of sample size $n$. Plots are on the same scale as Figure~\ref{fig:fig1}, and black line shows the minimax rate. All tests are set to have $.05$ Type I error, and are calibrated by simulation under the null.}
	\label{fig:fig2}
\end{figure*}

In Figure~\ref{fig:fig2}, we see that the estimated critical radii of both the PCR-LE and population-level spectral series tests are quite close to each other, and converge at roughly the minimax rate.

\begin{figure*}[tb]
	\includegraphics[width=.245\textwidth]{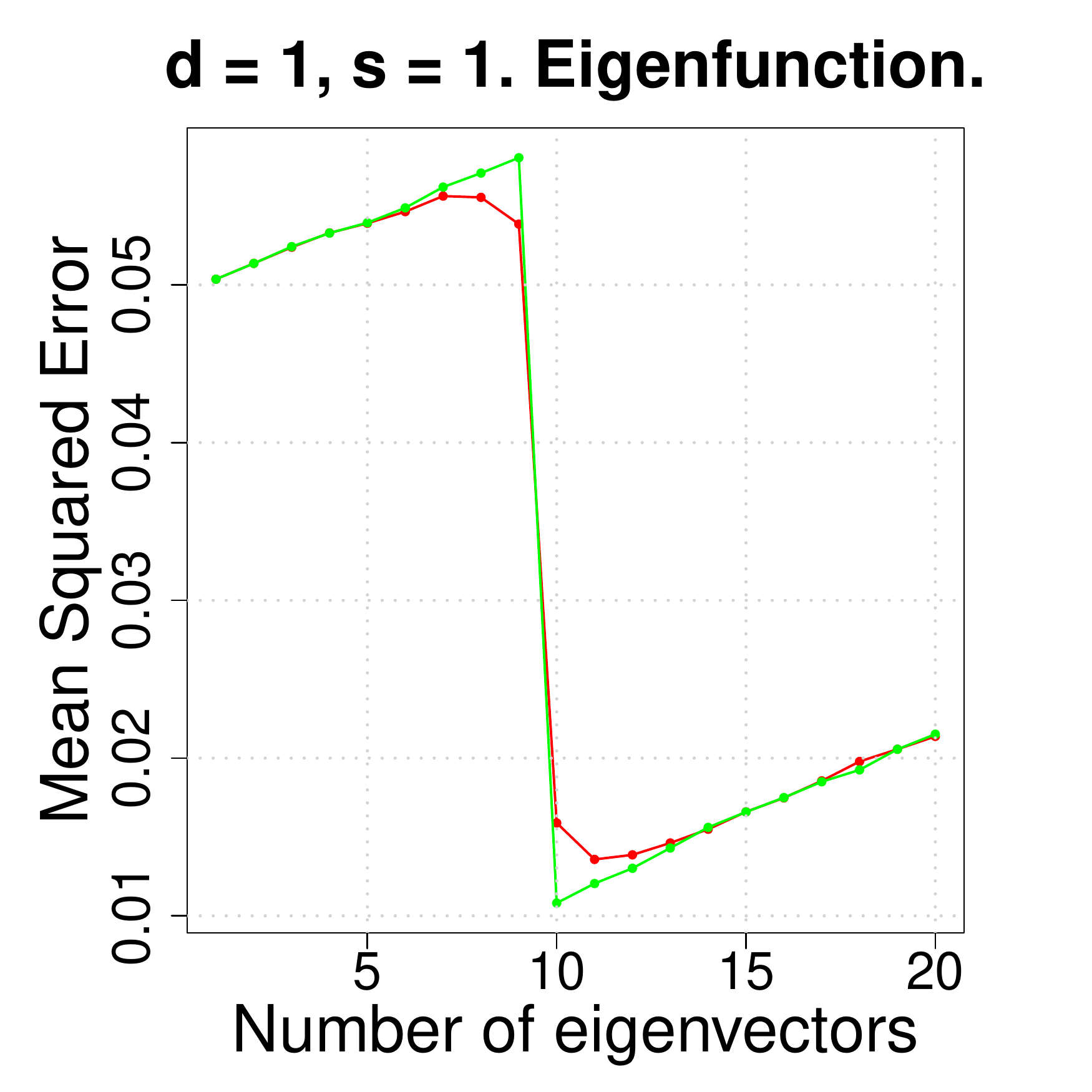}
	\includegraphics[width=.245\textwidth]{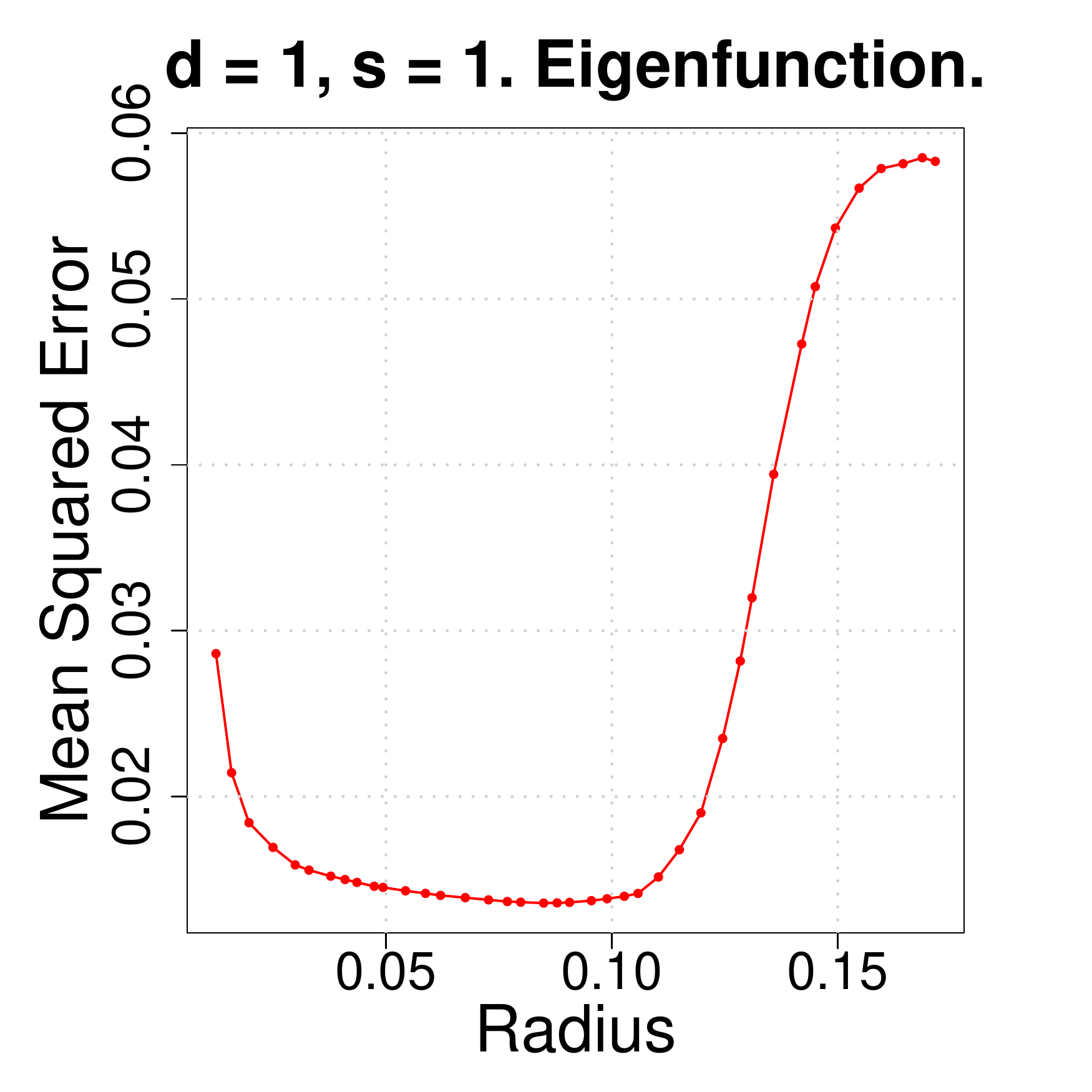} 
	\includegraphics[width=.245\textwidth]{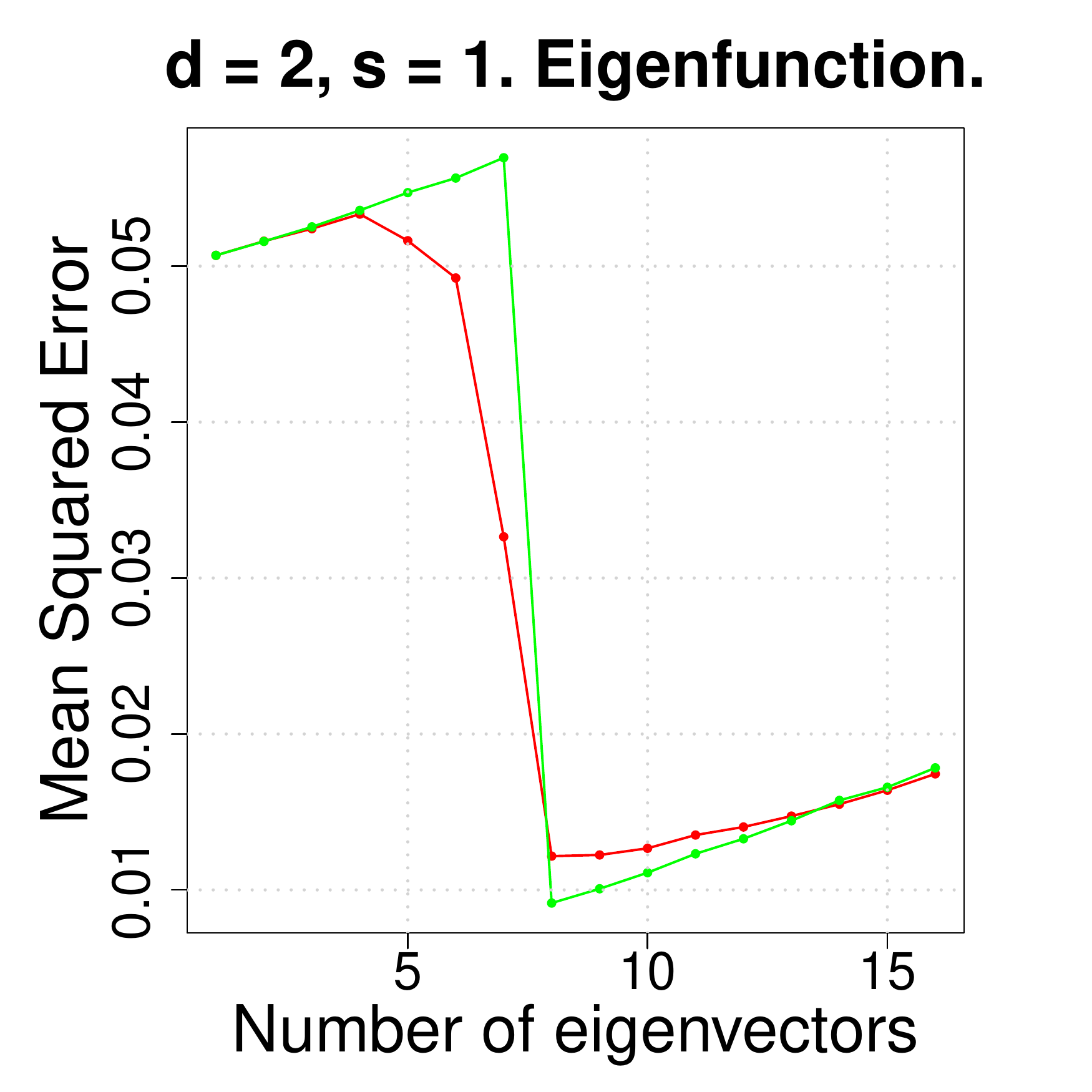}
	\includegraphics[width=.245\textwidth]{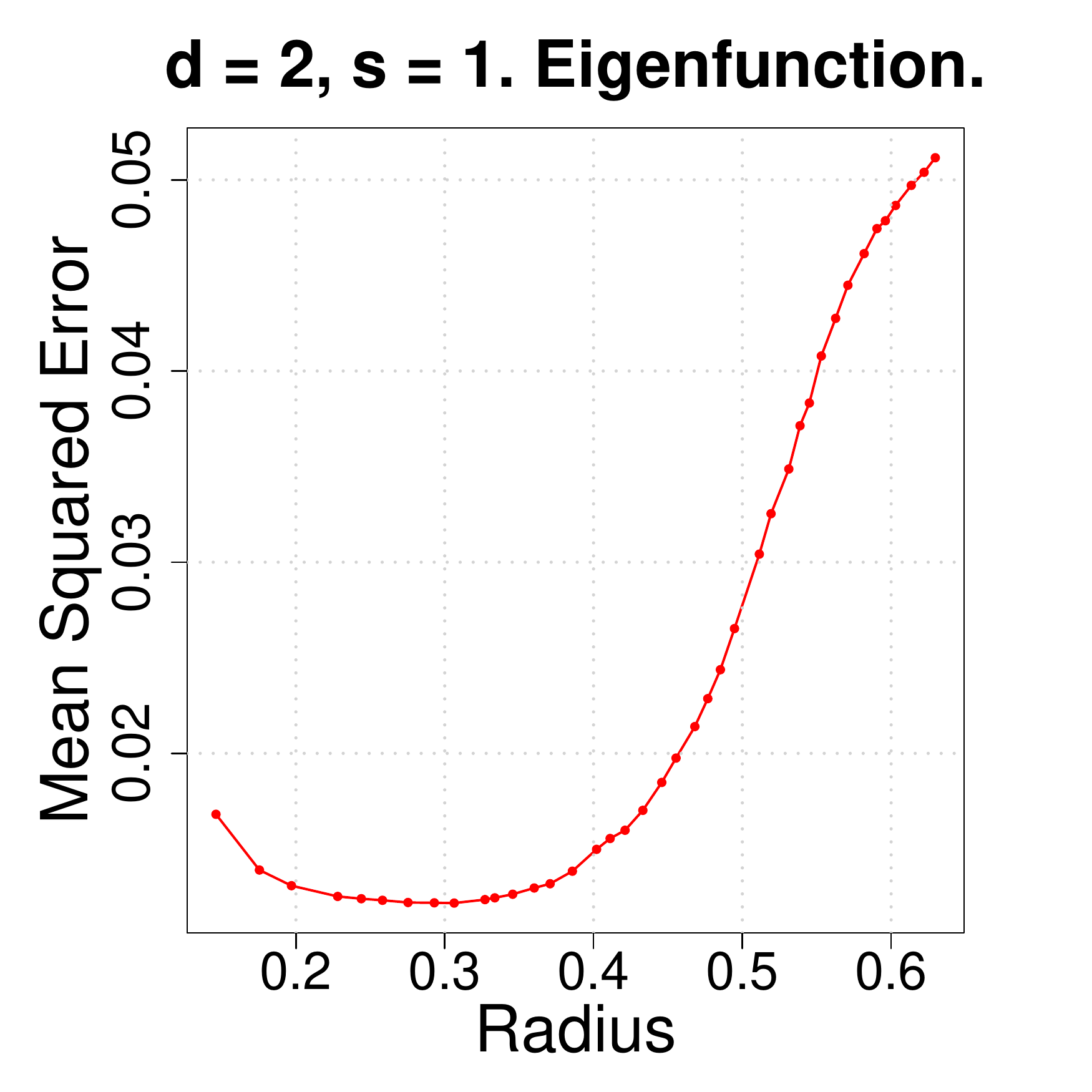} 
	\includegraphics[width=.245\textwidth]{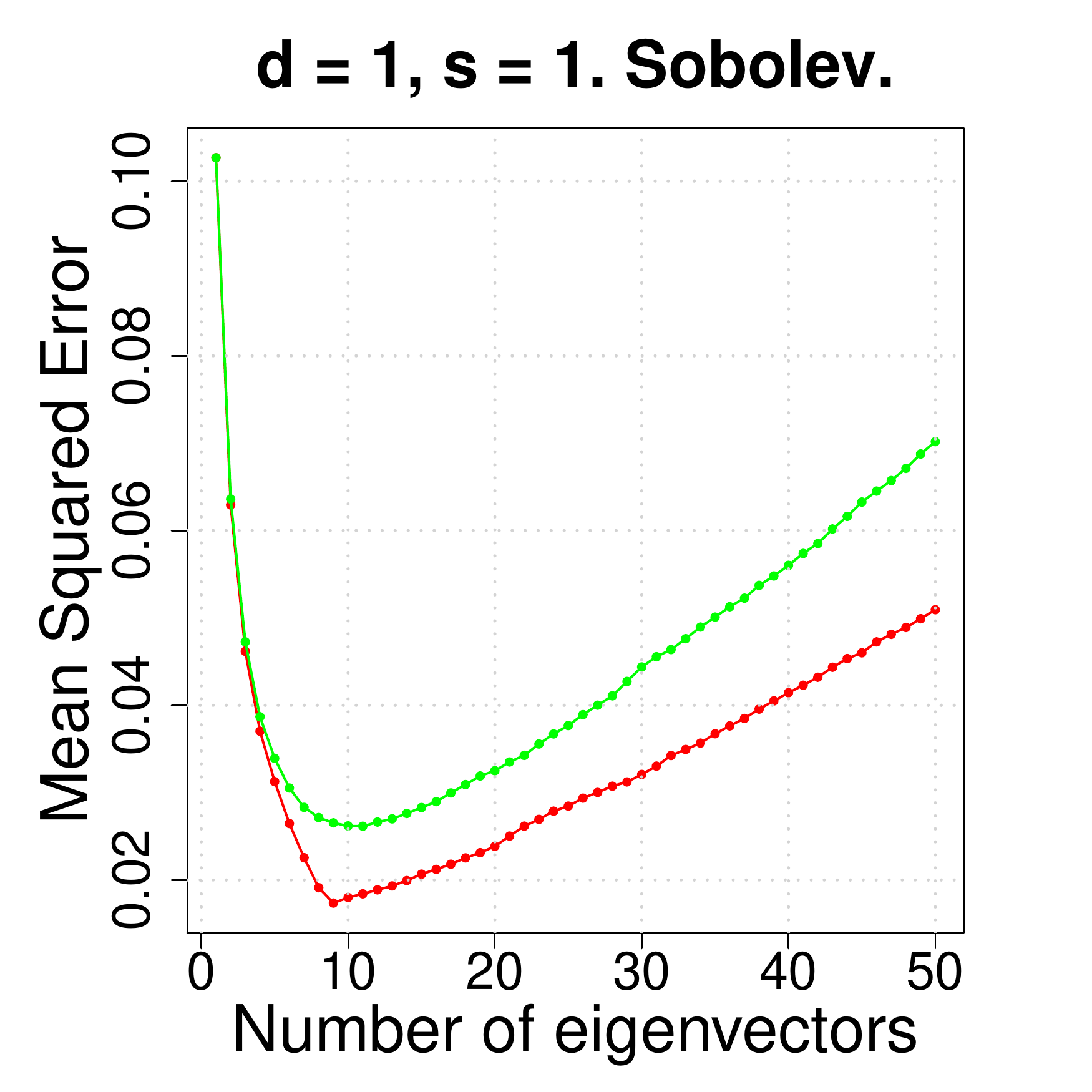}
	\includegraphics[width=.245\textwidth]{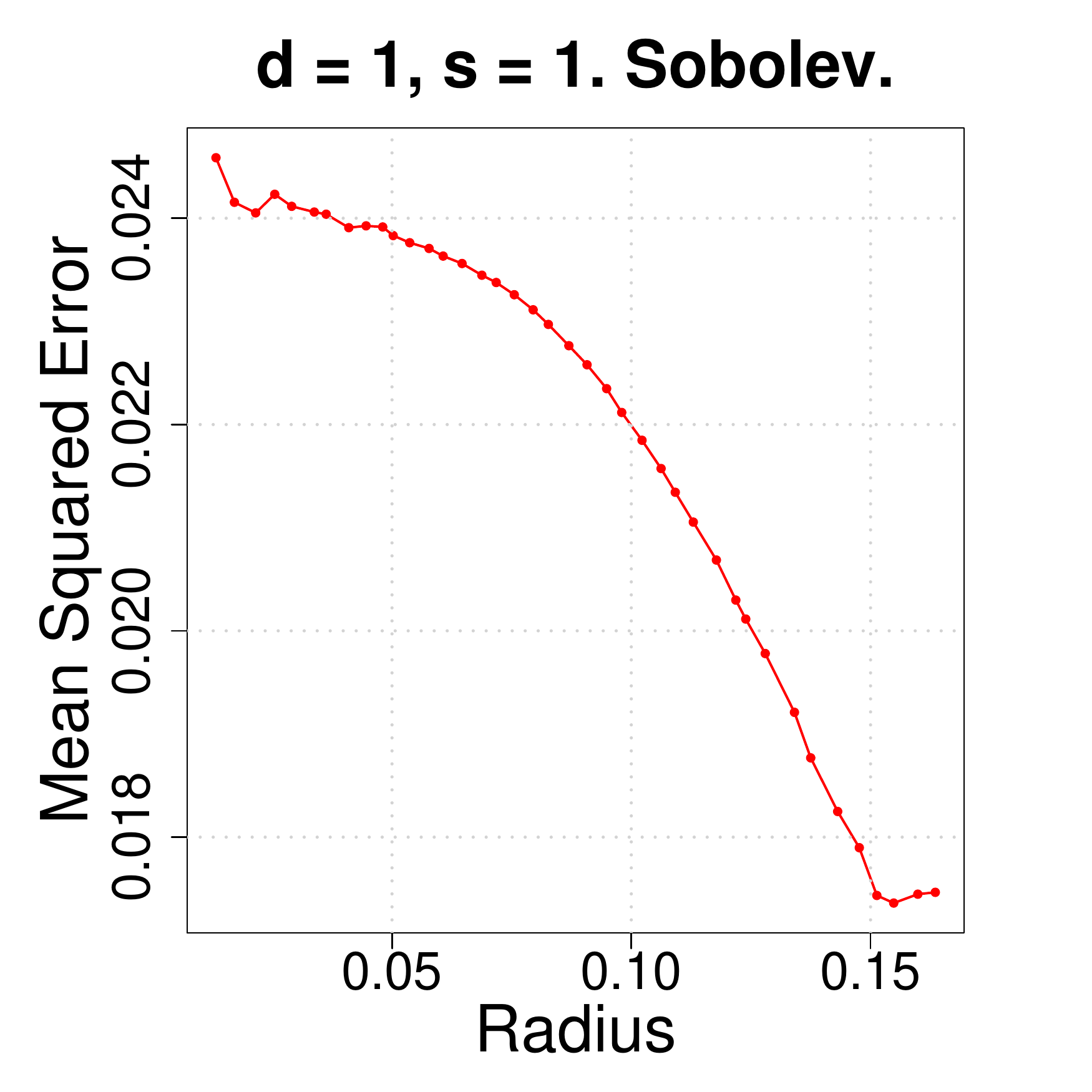}
	\includegraphics[width=.245\textwidth]{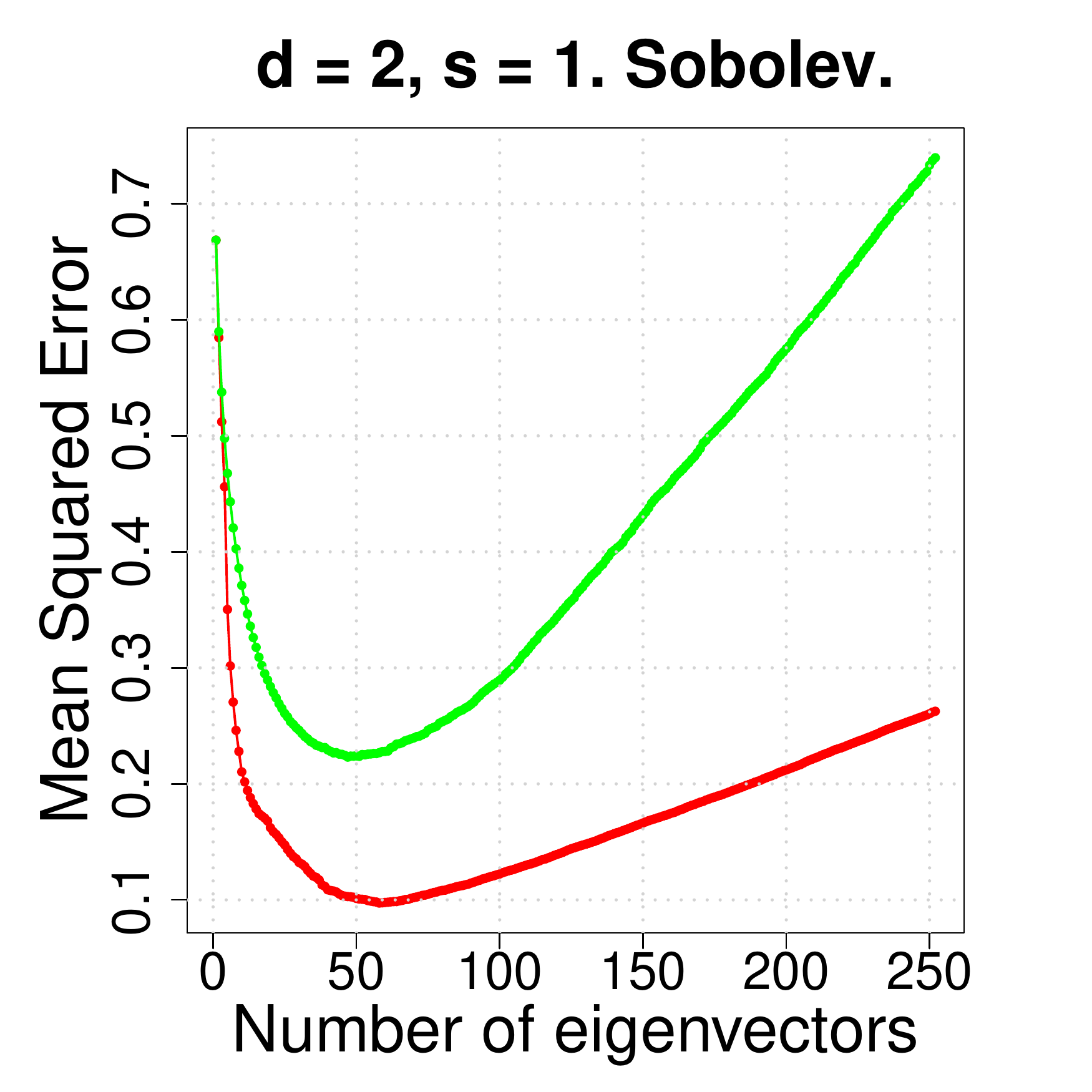}
	\includegraphics[width=.245\textwidth]{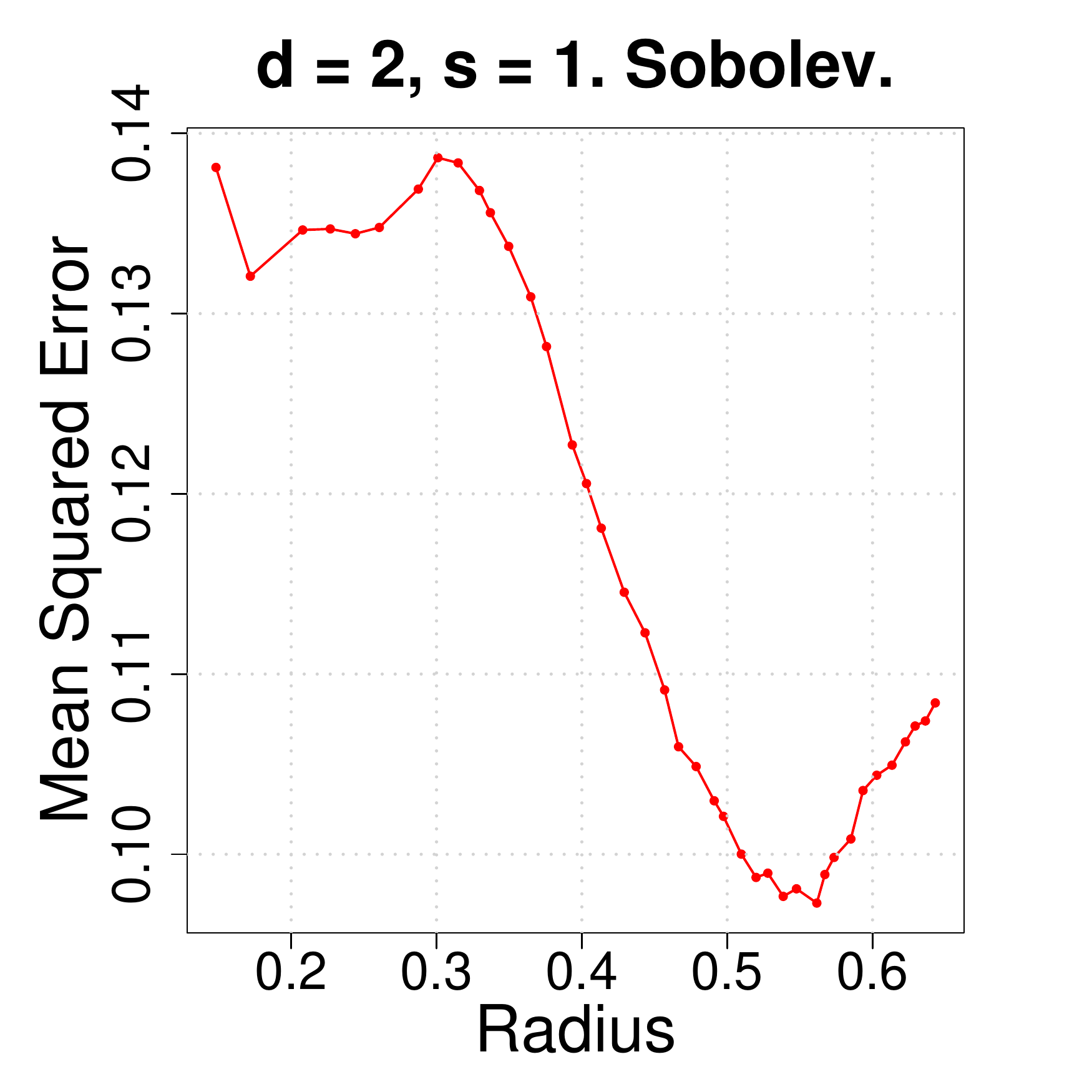}  
	\caption{Mean squared error of PCR-LE (\textcolor{red}{red}), and population-level spectral series (\textcolor{green}{green}) estimators as a function of tuning parameters. Top row: the same regression function $f_0$ as used in Figure~\ref{fig:fig1}. Bottom row: the regression function $f_0 \propto \sum_{k} 1/\rho_k^{1/2} \psi_k$. For all experiments, the sample size $n = 1000$, and the results are averaged over $200$ repetitions. In each panel, all tuning parameters except the one being varied are set to their optimal values.}
	\label{fig:fig3}
\end{figure*}

\paragraph{Tuning parameters.}
Our first two experiments demonstrate that PCR-LE methods have comparable statistical performance to population-level spectral series methods. PCR-LE depends on two tuning parameters, and in our final experiment we investigate the importance of both, focusing now on estimation. In Figure~\ref{fig:fig3}, we see how the mean-squared error of PCR-LE changes as each tuning parameter is varied. As suggested by our theory, properly choosing the number of eigenvectors $K$ is crucial: the mean-squared error curves, as a function of $K$, always have a sharply defined minimum. On the other hand, as a function of the graph radius parameter $\varepsilon$ the mean-squared error curve is much closer to flat. This squares completely with our theory, which requires that the number of eigenvectors $K$ be much more carefully tuned that the graph radius $\varepsilon$.

\section{Discussion}
\label{sec:discussion}

In this work, we have derived upper bounds on the rates of convergence for regression with PCR-LE, which imply that in various settings the PCR-LE estimator and test are minimax rate-optimal over Sobolev classes. Importantly, these upper bounds hold under nonparametric conditions on the design density $p$, and allow for $p$ to be unknown and, potentially, supported on a low-dimensional manifold. Our results help explain the practical success of methods which leverage graph Laplacian eigenvectors for regression. They also distinguish such methods from more traditional spectral series procedures, which rely on a density-dependent basis and thus require the density be known a priori.

Of course, there do exist other methods for nonparametric regression which achieve optimal rates of convergence under similar (or indeed weaker) conditions on $p$. These include other graph-based approaches---Laplacian smoothing---methods besides spectral series methods---e.g. kernel smoothing, local polynomial regression, thin-plate splines---and continuum spectral projection methods which use the eigenfunctions of an operator defined independently of $p$. To be clear, we do not advocate PCR-LE over these alternatives. Rather, we view our results as theoretically justifying a place for regression using Laplacian Eigenmaps in the nonparametric regression toolbox. 

That being said, PCR-LE does have certain advantages over each of the aforementioned approaches. We now conclude by outlining some of these advantages (limiting our discussion to estimation):
\begin{itemize}
	\item \emph{Optimality over high-dimensional Sobolev spaces}. As mentioned in the introduction, Laplacian smoothing (defined via~\eqref{eqn:laplacian_smoothing}) provably achieves minimax optimal rates over $H^1(\mc{X})$ only when $d \in \{1,2,3,4\}$ \citep{sadhanala16, green2021}. In contrast, PCR-LE is optimal over $H^1(\mc{X})$ for all dimensions $d$, and also over the higher-order Sobolev spaces $H^s(\mc{X})$. 
	\item \emph{Manifold adaptivity}. When the design distribution is non-uniform, an oft-recommended alternative to population-level spectral series regression is to run OLS using eigenfunctions of a density-independent differential operator. As a concrete example, let $\Delta$ be the unweighted Laplacian operator on $\Rd$, $\Delta = \sum_{i = 1}^{d} \partial^2f/\partial x_i^2$. Denoting the eigenfunctions of $\Delta$ (under Neumann boundary conditions) by $\phi_1,\phi_2,\ldots$, and letting $\Phi \in \mathbb{R}^{n \times K}$ be the matrix with entries $\Phi_{ik} = \phi_k(X_i)$ and columns $\Phi_1,\ldots,\Phi_K$,  one could compute an estimator by solving the following OLS problem:
	\begin{equation*}
	\minimize_{f \in \mathrm{span}\{\Phi_{1},\ldots,\Phi_{K}\}} \|{\bf Y} - f\|_n^2.
	\end{equation*}
	Unlike with spectral series regression, this approach can produce reasonable estimates even when the sampled eigenfunctions $(\phi_k(X_1),\ldots,\phi_k(X_n)) \in \Reals^n$ are not approximately orthogonal. Indeed, under the conditions of Model~\ref{def:model_flat_euclidean}, such a method will in fact be minimax rate-optimal, though the upper bounds may come with undesirably large constants if $p$ is very non-uniform. However under Model~\ref{def:model_manifold}, this method cannot achieve the faster minimax rates of convergence---which depend only on the intrinsic dimension $m$---and may even be inconsistent. This is because the eigenfunctions $\phi_k$ have no underlying relationship to the Sobolev space $H^s(\mc{X})$ except when $\mc{X}$ is a full-dimension set in $\Rd$. In contrast, PCR-LE uses features which are empirical approximations to eigenfunctions $\psi_k$ of the density-weighted Laplace-Beltrami operator $\Delta_P$. The eigenfunctions of $\Delta_P$ are appropriately adapted to the geometry of the manifold $\mc{X}$, and as a result PCR-LE is consistent and in certain cases minimax optimal, as we have shown. 
	\item \emph{Density adaptivity}. In Appendix~\ref{subsec:eigenmaps_beats_kernel_smoothing}, we give a simple univariate example of a sequence of densities and regression functions $\{(p^{(n)}, f_0^{(n)}: n \in \mathbb{N}\}$ such that the expected in-sample mean squared error of PCR-LE is smaller than that of either kernel smoothing or least squares using eigenfunctions of $\Delta$. This is possible because PCR-LE induces a completely different bias than these latter two methods. In particular, when $f_0$ and $p$ satisfy the so-called \emph{cluster assumption}---meaning $f_0$ is piecewise constant in high-density regions (clusters) of $p$---then the bias of PCR-LE can be much smaller (for equivalent levels of variance) than that of kernel smoothing or least-squares with eigenfunctions of $\Delta$. 
	
	We emphasize that this does not contradict the well-known optimality properties of, for example, kernel smoothing over H\"{o}lder balls. Rather, in the standard nonparametric regression setup---which we adopt in the main part of this paper, and in which $P$ is assumed to be equivalent to Lebesgue measure---the biases of PCR-LE and kernel smoothing happen to be equivalent. But when $P$ is sufficiently non-uniform, this is no longer the case.
\end{itemize}
Grounding each of these three points on a firmer and more complete theoretical basis would be, in our view, a valuable direction for future work.

%% file: appendix.tex
\noindent 

\section{Notation Table}
\label{sec:notation_table}

\begin{table}[h]
	\begin{center}
		\begin{tabular}{p{.14\textwidth} | p{.5\textwidth} }
			Symbol & Definition
			\\
			\hline
			$\mc{X}$ & domain, either an open set in $\Rd$ (Model~\ref{def:model_flat_euclidean}) or a compact manifold embedded in $\Rd$ (Model~\ref{def:model_manifold}) \\
			$\nu$ & Lebesgue measure \\
			$\mu$ & volume form induced by the embedding of $\mc{X}$ into $\Rd$ \\
			$P$ & probability measure associated with the design points \\
			$p$ & density of the probability measure, either with respect to $\nu$ (Model~\ref{def:model_flat_euclidean}) or $\mu$ (Model~\ref{def:model_manifold}). \\
			$L^2(\mc{X})$ & set of square-integrable functions, meaning either $\int_{\mc{X}} f^2 \,d\nu < \infty$ (Model~\ref{def:model_flat_euclidean}) or $\int_{\mc{X}} f^2 \,d\mu < \infty$ (Model~\ref{def:model_manifold}). \\
			$C^k(\mc{X})$ & functions which are $k$-times continuously differentiable in $\mc{X}$ \\
			$C_c^{\infty}(\mc{X})$ & functions in $C^{\infty}(\mc{X})$ which are compactly supported in $\mc{X}$ \\
			$H^s(\mc{X})$ & order-s Sobolev space (Definition~\ref{def:sobolev_space} under Model~\ref{def:model_flat_euclidean}, Definition~\ref{def:sobolev_space_manifold} under Model~\ref{def:model_manifold}.) \\
			$H_0^s(\mc{X})$ & order-s zero-trace Sobolev space (Definition~\ref{def:zero_trace_sobolev_space}) \\
			$\|\cdot\|_2$ & Euclidean distance \\
			$d_{\mc{X}}(\cdot,\cdot)$ & geodesic distance \\
			$B(x,\delta)$ & Ball in Euclidean distance, centered at $x$ with radius $\delta$ \\
			$B_{\mc{X}}(x,\delta)$ & Ball in geodesic distance
		\end{tabular}
	\end{center}
	\caption{Notation.}
\end{table}

\section{Upper bounds on population-level spectral series regression}

In this section we first give the proof of Proposition~\ref{prop:spectral_series_estimation}, then of Proposition~\ref{prop:spectral_series_testing}. In both cases the structure of the analysis, which is fairly classical and straightforward, can be usefully compared to our analysis of PCR-LE (see Section~\ref{subsec:analysis}).

\paragraph{Proof of Proposition~\ref{prop:spectral_series_estimation}.}
We decompose risk into squared bias and variance,
\begin{equation}
\label{pf:spectral_series_estimation_0}
\Ebb \|\wt{f} - f_0\|_P^2 = \Ebb\| \Ebb[\wt{f}]  - f_0\|_P^2 + \Ebb\| \wt{f} - \Ebb[\wt{f}]\|_P^2.
\end{equation}
Since the eigenfunctions $\{\psi_k\}$ form an orthonormal basis of $L^2(\mc{X})$ (with respect to the inner-product $\dotp{\cdot}{\cdot}_P$) and $f_0 \in \mc{H}^{s}(\mc{X}) \subseteq L^2(\mc{X})$, we can write the squared bias in terms of squared Fourier coefficients of $f_0$, leading to the following upper bound,
\begin{equation*}
\|f_0 - \Ebb[\wt{f}]\|_P^2 = \sum_{k = K + 1}^{\infty}  \dotp{f_0}{\psi_k}_P^2 \leq  \frac{1}{\{\lambda_{K + 1}(\Delta_P)\}^s} \sum_{k = K + 1}^{\infty} \{\lambda_{k + 1}(\Delta_P)\}^s \dotp{f_0}{\psi_k}_P^2 \leq \frac{\|f_0\|_{\mc{H}^s(\mc{X})}}{\{\lambda_{K + 1}(\Delta_P)\}^s}.
\end{equation*}
On the other hand, the variance term can be written as the sum of the variance of each empirical Fourier coefficient, and subsequently by the law of total variance we derive that
\begin{align}
\label{pf:spectral_series_estimation_2}
\Ebb\| \wt{f} - \Ebb[\wt{f}]\|_P^2 = \sum_{k = 1}^{K} \Var\Bigl[\dotp{{\bf Y}}{\psi_k}_n\Bigr] & = \sum_{k = 1}^{K} \Var\Bigl[\Ebb[\dotp{Y}{\psi_k}_n|{\bf X}\Bigr] + \Ebb\Bigl[\Var[\dotp{Y}{\psi_k}_n|{\bf X}\Bigr] \nonumber \\
& = \sum_{k = 1}^{K} \Var\Bigl[\dotp{f_0}{\psi_k}_n\Bigr] + \frac{1}{n}\Ebb\Bigl[\|\psi_k\|_n^2\Bigr] \nonumber \\
& \leq \frac{K}{n} + \frac{1}{n}\sum_{k = 1}^{K}\Ebb\Bigl[\Bigl(f_0(X)\psi_k(X)\Bigr)^2\Bigr].
\end{align}
Consequently,
\begin{equation}
\label{pf:spectral_series_estimation_1}
\Ebb \|\wt{f} - f_0\|_P^2 \leq \frac{\|f_0\|_{\mc{H}^s(\mc{X})}^2}{\bigl[\lambda_{K + 1}(\Delta_P)\bigr]^s} + \frac{K}{n} + \frac{1}{n}\Ebb\Bigl[(f_0(X))^2 \cdot \sum_{k = 1}^{K} (\psi_k(X))^2\Bigr].
\end{equation}
The claim of the proposition then follows from variants of two classical results in spectral geometry. The first is a Weyl's Law asymptotic scaling of the eigenvalues of $\Delta_P$ due to~\cite{dunlop2020}; formally, there exist constants $c$ and $C$ (which will depend on $P$ and $d$) such that
\begin{equation}
\label{eqn:weyl}
ck^{2/d} \leq \rho_k \leq Ck^{2/d}\quad\textrm{for all $k \in \mathbb{N}$, $k \geq 2$}.
\end{equation}
The second is a local analog to Weyl's Law, which says that there exists a constant $C$ (again depending on $P$ and $d$) such that
\begin{equation}
\label{eqn:local_weyl}
\sup_{x \in \mc{X}}\biggl\{\sum_{k = 1}^{K} \bigl(\psi_k(x)\bigr)^2\biggr\} \leq CK \quad\textrm{for all $K \in \mathbb{N}$}.
\end{equation}
Equation~\eqref{eqn:local_weyl} is a direct implication of~\eqref{eqn:weyl} along with Theorem 17.5.3 of~\cite{hormander1973}. Plugging the upper bounds~\eqref{eqn:weyl} and~\eqref{eqn:local_weyl} back into~\eqref{pf:spectral_series_estimation_1}, we conclude that
\begin{equation}
\label{pf:spectral_series_estimation_3}
\Ebb \|\wt{f} - f_0\|_P^2 \leq C\biggl(\frac{\|f_0\|_{\mc{H}^s(\mc{X})}^2}{(K + 1)^{2s/d}} + \frac{K}{n}\biggr).
\end{equation}
If $n^{-1/2} \geq M$, then taking $K = 1$ implies $\Ebb \|\wt{f} - f_0\|_P^2 \leq C(M^2 + 1/n)$. Otherwise, setting $K = \floor{M^2n}^{d/(2s + d)}$ balances squared bias and variance, and yields the claim.
\qed.

\paragraph{Proof of Proposition~\ref{prop:spectral_series_testing}.}
We briefly lay out the main ideas needed to prove Proposition~\ref{prop:spectral_series_testing}, following the lead of~\cite{ingster2009} who prove a similar result in the special case where $M = 1$ and $P$ is the uniform distribution over $\mc{X} = [0,1]^d$, and referring to that work for more details.

We begin by computing the first two moments of the test statistic $\wt{T}$. The expectation is
\begin{align*}
\Ebb[\wt{T}] = \frac{(n - 1)}{n} \sum_{k = 1}^{K} \dotp{f_0}{\psi_K}_P^2 + \frac{K}{n} + \Ebb\Bigl[(f_0(X))^2 \dot \sum_{k = 1}^{K} (\psi_k(X))^2\Bigr],
\end{align*}
and from~\eqref{eqn:weyl} (Weyl's Law) we have that under the alternative $f_0 \neq 0$,
\begin{equation*}
\Ebb_{f_0}[\wt{T}] \geq \frac{\|f_0\|_{\mc{H}^s(\mc{X})}^2}{\lambda_{K + 1}^s(\Delta_{P})} + \frac{K}{n}.
\end{equation*}
To compute the variance, we decompose $\wt{T} = \wt{T}_{1,1} + \wt{T}_{1,2} + \wt{T}_{1,3} + \wt{T}_2$ into the sum of 3 U-statistics and the remaining diagonal terms, defined in terms of the equivalent kernel $\kappa(x,x') = \sum_{k = 1}^{K} \psi_k(x) \psi_k(x')$ as,
\begin{align*}
T_{1,1} & := \frac{1}{n^2} \sum_{1 \leq i \neq j \leq n} w_i w_j \kappa(X_i,X_j),\quad && T_{1,2} := \frac{1}{n^2} \sum_{1 \leq i \neq j \leq n} \bigl(w_i f_0(X_j) + w_jf_0(X_i)\bigr)\kappa(X_i,X_j) \\
T_{1,3} & := \frac{1}{n^2} \sum_{1 \leq i \neq j \leq n} f_0(X_i) f_0(X_j) \kappa(X_i,X_j), \quad && T_2 := \frac{1}{n^2} \sum_{i = 1}^{n} Y_i^2 \kappa(X_i,X_i).
\end{align*}
The variances of each statistic can be found by routine computation (see~\cite{ingster2009}), and in particular satisfy the upper bounds
\begin{align*}
\Var(T_{1,1}) & \leq \frac{2K}{n^2}, \quad && \Var(T_{1,2}) \overset{\mathrm{(i)}}{\leq} \frac{C}{n}\|f_0\|_P^2\\
\Var(T_{1,3}) & \overset{\mathrm{(ii)}}{\leq} C\biggl(\frac{K}{n}\|f_0\|_P^4 + \frac{K}{n^2}\|f_0\|_{L^4(\mc{X})}^4\biggr), \quad && \Var(T_{2}) \overset{\mathrm{(iii)}}{\leq} \frac{CK^2}{n^3}\biggl(1 + \|f_0\|_{L^4(\mc{X})}^4\biggr)
\end{align*}
where $\mathrm{(i)}-\mathrm{(iii)}$ hold due to local Weyl's law, i.e.~\eqref{eqn:local_weyl}. Upper bounds on Type I and Type II error,
\begin{equation*}
\Ebb_0[\wt{\varphi}] \leq \biggl(1 + CK/n^2\biggr)a, \quad \Ebb_{f_0}[1 - \wt{\varphi}] \leq \frac{C(K/n^2 + \|f_0\|_P^2 + K/n \|f_0\|_P^4 + K/n^2 \|f_0\|_{L^4(\mc{X})}^4)}{(\sum_{k = 1}^{K} \dotp{f_0}{\psi_k}_P^2 - \sqrt{2K/an})^2},
\end{equation*}
follow from Chebyshev's inequality. It can be verified that so long as
\begin{equation}
\label{pf:spectral_series_test}
\|f_0\|_P^2 \geq C\Biggl(\frac{\|f_0\|_{\mc{H}^s(\mc{X})}^2}{\lambda_{K + 1}^s(\Delta_{P})} + \frac{\sqrt{K}}{n}\biggl(\sqrt{\frac{1}{a}} + \sqrt{\frac{1}{b}}\biggr)\Biggr)
\end{equation}
for a sufficiently large constant $C$, then $\Ebb_{f_0}[1 - \wt{\varphi}] \leq b$. The two summands in~\eqref{pf:spectral_series_test} are bias and standard deviation terms, respectively. When $M^2 \leq n^{-1}$, setting $K = 1$ gives the desired result. Otherwise, choosing $K = \floor{M^2n}^{2d/(4s + d)}$ balances these two terms, and leads to~\eqref{eqn:spectral_series_testing}. \qed

\section{Graph-dependent error bounds}
\label{sec:fixed_graph_error_bounds}
In this section, we adopt the fixed design perspective; or equivalently, condition on $X_i = x_i$ for $i = 1,\ldots,n$. Let $G = \bigl([n],W\bigr)$ be a fixed graph on $\{1,\ldots,n\}$ with Laplacian matrix $L = \sum_{k = 1}^{n}\lambda_k v_k v_k^{\top}$; the eigenvectors have unit empirical norm, $\|v_k\|_n^2 = 1$. The randomness thus all comes from the responses 
\begin{equation}
\label{eqn:fixed_graph_regression_model}
Y_i = f_{0}(x_i) + w_i
\end{equation}
where the noise variables $w_i$ are independent $N(0,1)$. In the rest of this section, we will mildly abuse notation and write $f_0 = (f_0(x_1),\ldots,f_0(x_n)) \in \Reals^n$. We will also write ${\bf Y} = (Y_1,\ldots,Y_n)$.

\subsection{Upper bound on Estimation Error of PCR-LE}

\begin{lemma}
	\label{lem:fixed_graph_estimation}
	For any integer $s > 0$, and any integer $0 \leq K \leq n$, the PCR-LE estimator $\wh{f}$ of~\eqref{eqn:laplacian_eigenmaps_estimator} satisfies
	\begin{equation}
	\label{eqn:fixed_graph_estimation}
	\|\wh{f} - f_0\|_n^2 \leq \frac{\dotp{L^sf_0}{f_0}_n}{\lambda_{K + 1}^s} + \frac{5K}{n};
	\end{equation}
	this is guaranteed if $K = 0$, and otherwise holds with probability at least $1 - \exp(-K)$ if $1 \leq K \leq n$. 
\end{lemma}
\paragraph{Proof (of Lemma~\ref{lem:fixed_graph_estimation}).}
	By the triangle inequality,
	\begin{equation}
	\label{pf:fixed_graph_estimation_1}
	\|\wh{f} - f_0\|_n^2 \leq 2\Bigl(\|\mathbb{E}\wh{f} - f_0\|_n^2 + \|\wh{f} - \mathbb{E}\wh{f}\|_n^2\Bigr).
	\end{equation}
	The first term in~\eqref{pf:fixed_graph_estimation_1} (approximation error) is non-random, since the design is fixed. The expectation $\mathbb{E}\wh{f} = \sum_{k = 1}^{K} \dotp{v_k}{f_0}_n v_k$, so that
	\begin{equation*}
	\|\mathbb{E}\wh{f} - f_0\|_n^2 = \Bigl\|\sum_{k = K + 1}^{n} \dotp{v_k}{f_0}_n v_k\Bigr\|_n^2 = \sum_{k = K + 1}^n \dotp{v_k}{f_0}_n^2.
	\end{equation*}
	In the above, the last equality relies on the fact that $v_k$ are orthonormal with respect to $\dotp{\cdot}{\cdot}_n$. Using the fact that the eigenvalues are in increasing order, we obtain
	\begin{equation*}
	\sum_{k = K + 1}^n \dotp{v_k}{f_0}_n^2 \leq \frac{1}{\lambda_{K + 1}^s} \sum_{k = K + 1}^n \lambda_k^s \dotp{v_k}{f_0}_n^2 \leq \frac{\dotp{L^sf_0}{f_0}_n}{\lambda_{K + 1}^s}.
	\end{equation*}
	
	If $K = 0$, $\wh{f} = \Ebb{\wh{f}} = 0$, and the second term in~\eqref{pf:fixed_graph_estimation_1} is $0$. Otherwise the second   in~\eqref{pf:fixed_graph_estimation_1} (estimation error) is random. Observe that $\dotp{v_k}{\varepsilon}_n \overset{d}{=} Z_k/\sqrt{n}$, where $(Z_1,\ldots,Z_n) \sim N(0,I_{n \times n})$. Again using the orthonormality of the eigenvectors $v_k$, we have
	\begin{equation*}
	\|\wh{f} - \mathbb{E}\wh{f}\|_n^2 = \sum_{k = 1}^{K} \dotp{v_k}{\varepsilon}_n^2 \overset{d}{=} \frac{1}{n}\sum_{k = 1}^{K} Z_k^2.
	\end{equation*}
	Thus $\|\wh{f} - \mathbb{E}\wh{f}\|_n^2$ is equal to $1/n$ times a $\chi^2$ distribution with $K$ degrees of freedom. Consequently, it follows from a result of \citep{laurent00} that
	\begin{equation*}
	\Pbb\biggl(\|\wh{f} - \mathbb{E}\wh{f}\|_n^2 \geq \frac{K}{n} + 2\frac{\sqrt{K}}{n}\sqrt{t} + \frac{2t}{n}\biggr) \leq \exp(-t).
	\end{equation*}
	Setting $t = K$ completes the proof of the lemma.

\subsection{Upper bound on Testing Error of PCR-LE}

Let $\wh{T} = \sum_{k = 1}^{K} \dotp{{\bf Y}}{v_k}_n^2$, and let $\varphi = \1\{\wh{T} \geq t_a\}$. In the following Lemma, we upper bound the Type I and Type II error of the test $\varphi$.

\begin{lemma}
	\label{lem:fixed_graph_testing}
	Suppose we observe $(Y_1,x_1),\ldots,(Y_n,x_n)$ according to~\eqref{eqn:fixed_graph_regression_model}.
	\begin{itemize}
		\item If $f_0 = 0$, then $\Ebb_0[\varphi] \leq a$.
		\item Suppose $f_0 \neq 0$ satisfies
		\begin{equation}
		\label{eqn:fixed_graph_testing_critical_radius}
		\|f_0\|_n^2 \geq \frac{\dotp{L^sf_0}{f_0}_n}{\lambda_{K + 1}^s} + \frac{\sqrt{2K}}{n}\biggl[2\sqrt{\frac{1}{a}} + \sqrt{\frac{2}{b}} + \frac{32}{bn}\biggr],
		\end{equation}
		for some $s \in \mathbb{N}\setminus \{0\}$. Then $\Ebb_{f_0}[1 - \phi] \leq b$.
	\end{itemize}
\end{lemma}
\paragraph{Proof (of Lemma~\ref{lem:fixed_graph_testing}).}
We first compute the expectation and variance of $\wh{T}$, then apply Chebyshev's inequality to upper bound the Type I and Type II error.

\underline{\emph{Expectation}.}
Recall that $\wh{T} = \sum_{k = 1}^{K} \dotp{Y}{v_k}_n^2$. Expanding the square gives
\begin{equation*}
\Ebb[\wh{T}] = \sum_{k = 1}^{K} \Ebb[\dotp{Y}{v_k}_n^2] = \sum_{k = 1}^{K} \dotp{f_0}{v_k}_n^2 + \Ebb[2\dotp{f_0}{v_k}_n\dotp{\varepsilon}{v_k}_n + \dotp{\varepsilon}{v_k}_n^2] = \frac{K}{n} + \sum_{k = 1}^{K} \dotp{f_0}{v_k}_n^2.
\end{equation*}
Thus $\Ebb[\wh{T}] - t_a = \sum_{k = 1}^{K} \dotp{f_0}{v_k}_n^2 - \sqrt{2K}/n \cdot \sqrt{1/a}$. Furthermore, it is a consequence of~\eqref{eqn:fixed_graph_testing_critical_radius} that 
\begin{equation}
\label{pf:fixed_graph_testing_1}
\sum_{k = 1}^{K} \dotp{f_0}{v_k}_n^2 - \frac{\sqrt{2K}}{n}\sqrt{1/a} \geq \|f_0\|_n^2 - \frac{\dotp{L^sf_0}{f_0}_n}{\lambda_{K + 1}^s} - \frac{\sqrt{2K}}{n}\sqrt{1/a} \geq \frac{\sqrt{2K}}{n}\biggl[\sqrt{\frac{1}{a}} + \sqrt{\frac{2}{b}} + \frac{32}{bn}\biggr].
\end{equation} 

\underline{\emph{Variance}.}
Recall from the proof of Lemma~\ref{lem:fixed_graph_estimation} that $\dotp{\varepsilon}{v_k}_n \overset{d}{=} Z_k/\sqrt{n}$ for $(Z_1,\ldots,Z_n) \sim N(0,I_{n \times n})$. Expanding the square, and recalling that $\Cov[Z,Z^2] = 0$ for Gaussian random variables, we have that
\begin{equation*}
\Var\bigl[\dotp{{\bf Y}}{v_k}_n^2\bigr] = \Var\biggl[\frac{2}{n}\dotp{f_0}{v_k}_nZ_k + \frac{2}{n^2}Z_k^2\biggr] = \frac{4\dotp{f_0}{v_k}_n^2}{n} + \frac{2}{n^2}.
\end{equation*}
Moreover, since $\Cov[Z_k^2,Z_{\ell}^2] = 0$ for each $k = 1,\ldots,K$, we see that
\begin{equation*}
\Var\bigl[\wh{T}\bigr] = \sum_{k = 1}^{K} \Var\bigl[\dotp{{\bf Y}}{v_k}_n^2\bigr] = \frac{2K}{n^2} + \sum_{k = 1}^{K}\frac{4\dotp{f_0}{v_k}_n^2}{n}.
\end{equation*}

\underline{\emph{Bounds on Type I and Type II error}.}
The upper bound on Type I error follows immediately from Chebyshev's inequality. 

The upper bound on Type II error also follows from Chebyshev's inequality. We observe that~\eqref{eqn:fixed_graph_testing_critical_radius} implies $\Ebb_{f_0}[\wh{T}] = t_a$, and apply Chebyshev's inequality to deduce
\begin{equation*}
\Pbb_{f_0}\bigl(\wh{T} < t_a\bigr) \leq \Pbb_{f_0}\Bigl(|\wh{T} - \Ebb_{f_0}[\wh{T}]|^2 > |\Ebb_{f_0}[\wh{T}] - t_a|^2\Bigr) \leq \frac{\Var\bigl[\wh{T}\bigr]}{\bigl[\Ebb_{f_0}[\wh{T}] - t_a\bigr]^2} = \frac{2K/n^2 + 4/n\sum_{k = 1}^{K}\dotp{f_0}{v_k}_n^2}{\bigl[\Ebb_{f_0}[\wh{T}] - t_a\bigr]^2}.
\end{equation*}
Thus we have upper bounded the Type II error by the sum of two terms, each of which are no more than $1/(2b)$, as we now show. For the first term, after noting that~\eqref{pf:fixed_graph_testing_1} implies $\Ebb_{f_0}[\wh{T}] - t_a \geq \sqrt{2K}/n \cdot \sqrt{2/b}$, the upper bound follows:
\begin{equation*}
\frac{2K/n^2}{\bigl[\Ebb_{f_0}[\wh{T}] - t_a\bigr]^2} \leq \frac{b}{2}.
\end{equation*}
On the other hand, for the second term we use~\eqref{pf:fixed_graph_testing_1} in two ways: first to conclude that $\Ebb_{f_0}[\wh{T}] - t_a \geq 1/2 \cdot \sum_{k = 1}^{K}\dotp{f_0}{v_k}_n^2$, and second to obtain
\begin{equation*}
\frac{4\sum_{k = 1}^{K}\dotp{f_0}{v_k}_n^2}{n\bigl[\Ebb_{f_0}[\wh{T}] - t_a\bigr]^2} \leq \frac{4\sum_{k = 1}^{K}\dotp{f_0}{v_k}_n^2}{n\bigl(\sum_{k = 1}^{K}\dotp{f_0}{v_k}_n^2/2\bigr)^2} \leq \frac{16}{n\sum_{k = 1}^{K}\dotp{f_0}{v_k}_n^2} \leq \frac{b}{2}.
\end{equation*}

\section{Graph Sobolev semi-norm, flat Euclidean domain}
\label{sec:graph_quadratic_form_euclidean}
In this section we prove Proposition~\ref{prop:graph_seminorm_ho}. The proposition will follow from several intermediate results.
\begin{enumerate}
	\item~In Section~\ref{subsec:decomposition_graph_seminorm}, we show that
	\begin{equation}
	\label{pf:graph_seminorm_ho_1}
	\dotp{L_{n,\varepsilon}^sf}{f}_n \leq \frac{1}{\delta} \dotp{L_{P,\varepsilon}^sf}{f}_{P} + \frac{C\varepsilon^2}{\delta n\varepsilon^{2 + d}}M^2.
	\end{equation}
	with probability at least $1 - 2\delta$. 
	
	We term the first term on the right hand side the \emph{non-local Sobolev semi-norm}, as it is a kernelized approximation to the Sobolev semi-norm $\dotp{\Delta_P^sf}{f}_{P}$. The second term on the right hand side is a pure bias term, which as we will see is negligible compared to the non-local Sobolev semi-norm as long as $\varepsilon \ll n^{-1/(2(s -1 + d))}$. 
	\item~In Section~\ref{subsec:approximation_error_nonlocal_laplacian}, we show that when $x$ is sufficiently in the interior of $\mc{X}$, then $L_{P,\varepsilon}^kf(x)$ is a good approximation to $\Delta_P^kf(x)$, as long as $f \in H^{s}(\mc{X})$ and $p \in C^{s - 1}(\mc{X})$ for some $s \geq 2k + 1$. 
	\item~In Section~\ref{subsec:boundary_behavior_nonlocal_laplacian}, we show that when $x$ is sufficiently near the boundary of $\mc{X}$, then $L_{P,\varepsilon}^kf(x)$ is close to $0$, as long as $f \in H_0^{s}(\mc{X})$ for some $s > 2k$.
	\item~In Section~\ref{subsec:estimate_nonlocal_seminorm}, we use the results of the preceding two sections to show that if $f \in H_0^s(\mc{X};M)$ and $p \in C^{s - 1}(\mc{X})$, there exists a constant $C$ which does not depend on $f$ such that
	\begin{equation}
	\label{pf:graph_seminorm_ho_2}
	\dotp{L_{P,\varepsilon}^sf}{f}_{P} \leq CM^2.
	\end{equation}
\end{enumerate}
Finally, in Section~\ref{subsec:integrals} we provide some assorted estimates used in Sections~\ref{subsec:decomposition_graph_seminorm}. 

\paragraph{Proof (of Proposition~\ref{prop:graph_seminorm_ho}).}
Proposition~\ref{prop:graph_seminorm_ho} follows immediately from~\eqref{pf:graph_seminorm_ho_1} and~\eqref{pf:graph_seminorm_ho_2}. \qed

One note regarding notation: suppose a function $g \in H^{\ell}(U)$, where $\ell \in \mathbb{N}$ and $U$ is an open set. Let $V$ be another open set, compactly contained within $U$. Then we will use the notation $g \in H^{\ell}(V)$ to mean that the restriction $\restr{g}{V}$ of $g$ to $V$ belongs to $H^{\ell}(V)$.

\subsection{Decomposition of graph Sobolev semi-norm}
\label{subsec:decomposition_graph_seminorm}

In Lemma~\ref{lem:graph_seminorm_bias}, we decompose the graph Sobolev semi-norm (a V-statistic) into an unbiased estimate of the non-local Sobolev semi-norm (a U-statistic), and a pure bias term. We establish that the pure bias term will be small (in expectation) relative to the U-statistic whenever $\varepsilon$ is sufficiently small.
\begin{lemma}
	\label{lem:graph_seminorm_bias}
	For any $f \in L^2(\mc{X})$, the graph Sobolev semi-norm satisfies
	\begin{equation}
	\label{eqn:graph_seminorm_bias_1}
	\dotp{L_{n,\varepsilon}^sf}{f}_{n} = U_{n,\varepsilon}^{(s)}(f) + B_{n,\varepsilon}^{(s)}(f),
	\end{equation}
	such that $\mathbb{E}[U_{n,\varepsilon}^{(s)}(f)] = (n - s - 1)!/n! \cdot \dotp{L_{P,\varepsilon}^sf}{f}_P$. If additionally $f \in H^1(\mc{X};M)$ and $\varepsilon \geq n^{-1/d}$, then the bias term $B_{n,\varepsilon}^{(s)}(f)$ satisfies
	\begin{equation}
	\label{eqn:graph_seminorm_bias_2}
	\mathbb{E}\bigl[|B_{n,\varepsilon}^{(s)}(f)|\bigr] \leq \frac{C\varepsilon^2}{\delta n\varepsilon^{2 + d}}M^2.
	\end{equation}
\end{lemma}
Then~\ref{pf:graph_seminorm_ho_1} follows immediately from Lemma~\ref{lem:graph_seminorm_bias}, by Markov's inequality.
\paragraph{Proof (of Lemma~\ref{lem:graph_seminorm_bias}).}
We begin by introducing some notation. We will use bold notation $\bj = (j_1,\ldots,j_s)$ for a vector of indices where $j_i \in [n]$ for each $i$. We write $[n]^s$ for the collection of all such vectors, and $(n)^s$ for the subset of such vectors with no repeated indices. Finally, we write $D_if$ for a kernelized difference operator,
\begin{equation*}
D_if(x) := \bigl(f(x) - f(X_i)\bigr) \eta\biggl(\frac{\|X_i - x\|}{\varepsilon}\biggr),
\end{equation*}
and we let $D_{\bj}f(x) := \bigl(D_{j_1}\circ \cdots \circ D_{j_s}f\bigr)(x)$.

With this notation in hand, it is easy to represent $\dotp{L_{n,\varepsilon}^sf}{f}_{n}$ as the sum of a U-statistic and a bias term,
\begin{align*}
\dotp{L_{n,\varepsilon}^sf}{f}_{n} & = \frac{1}{n} \sum_{i = 1}^{n} L_{n,\varepsilon}^sf(X_i) \cdot f(X_i) \\
& = \underbrace{\frac{1}{n^{s + 1}\varepsilon^{s(d + 2)}} \sum_{{i\bf j} \in (n)^{s + 1}}D_{\bj}f(X_i) \cdot f(X_i)}_{=:U_{n,\varepsilon}^{(s)}(f)} + \underbrace{\frac{1}{n^{s + 1}\varepsilon^{s(d + 2)}} \sum_{\substack{i\bj \in \\ [n]^{s + 1}\setminus (n)^{s + 1}}} D_{\bj}f(X_i) \cdot f(X_i)}_{=:B_{n,\varepsilon}^{(s)}(f)}
\end{align*}
When the indices of $i\bj$ are all distinct, it follows straightforwardly from the law of iterated expectation that
\begin{equation*}
\mathbb{E}[ D_{\bj}f(X_i) \cdot f(X_i)] = \varepsilon^{s(d + 2)}\mathbb{E}[L_{P,\varepsilon}^sf(X_i) \cdot f(X_i)] = \dotp{L_{P,\varepsilon}^sf}{f}_{P}, 
\end{equation*}
which in turn implies $\mathbb{E}[U_{n,\varepsilon}^{(s)}(f)] = (n - s - 1)!/n! \cdot \dotp{L_{P,\varepsilon}^sf}{f}_P$. 

It remains to show~\eqref{eqn:graph_seminorm_bias_2}. By adding and subtracting $f(X_{\bj_1})$, we obtain by symmetry that
\begin{equation*}
\sum_{\substack{i\bj \in \\ [n]^{s + 1}\setminus (n)^{s + 1}}} D_{\bj}f(X_i) \cdot f(X_i) = \frac{1}{2} \cdot \sum_{\substack{i\bj \in \\ [n]^{s + 1}\setminus (n)^{s + 1}}} D_{\bj}f(X_i) \cdot \bigl(f(X_i) - f(X_{\bj_1})\bigr),
\end{equation*}
and consequently
\begin{equation*}
\Ebb\Bigl[\sum_{\substack{i\bj \in \\ [n]^{s + 1}\setminus (n)^{s + 1}}} D_{\bj}f(X_i) \cdot f(X_i)\Bigr] \leq \frac{1}{2} \cdot \sum_{\substack{i\bj \in \\ [n]^{s + 1}\setminus (n)^{s + 1}}} \Ebb\Bigl[\bigl|D_{\bj}f(X_i)\bigr| \cdot \bigl|f(X_i) - f(X_{\bj_1})\bigr|\Bigr].
\end{equation*}
In Lemma~\ref{lem:graph_seminorm_bias2}, we show that if $f \in H^1(\mc{X};M)$, then for any $i\bj \in [n]^{s + 1}$ which contains a total of $k + 1$ distinct indices, 
\begin{equation*}
\Ebb\Bigl[\bigl|D_{\bj}f(X_i)\bigr| \cdot \bigl|f(X_i) - f(X_{\bj_1})\bigr|\Bigr] \leq C_1 \varepsilon^{2 + kd} M^2.
\end{equation*}
This shows us that the expectation of $|B_{n,\varepsilon}^s(f)|$ can bounded from above by the sum over several different terms, as follows:
\begin{align*}
\Ebb\Bigl[|B_{n,\varepsilon}^s(f)|\Bigr] & \leq C_1\frac{\varepsilon^2}{n\varepsilon^{2s}}M^2 \sum_{\substack{i\bj \in \\ [n]^{s + 1}\setminus (n)^{s + 1}}} \frac{1}{(n\varepsilon^d)^s}  \varepsilon^{(|i\bj| - 1)d} \\
& \leq C_1\frac{\varepsilon^2}{n\varepsilon^{2s}}M^2  \sum_{k = 1}^{s - 1} \frac{(n\varepsilon^d)^k}{(n\varepsilon^d)^s}n.
\end{align*}
Finally, we note that by assumption $n\varepsilon^d \geq 1$, so that in the above sum the factor of $(n\varepsilon^d)^k$ is largest when $k = s- 1$. We conclude that
\begin{equation*}
\Ebb\Bigl[|B_{n,\varepsilon}^s(f)|\Bigr] \leq C_1 (s - 1) \frac{\varepsilon^2}{n\varepsilon^{2s + d}}M^2,
\end{equation*}
which is the desired result.

\subsection{Approximation error of non-local Laplacian}
\label{subsec:approximation_error_nonlocal_laplacian}

In this section, we establish the convergence $L_{P,\varepsilon}^kf \to \sigma_{\eta}^k\Delta_P^kf$ as $\varepsilon \to 0$. More precisely, we give an upper bound on the squared difference between $L_{P,\varepsilon}^kf$ and  $\sigma_{\eta}^k\Delta_P^kf$ as a function of $\varepsilon$. The bound holds for all $x \in \mc{X}_{k\varepsilon}$, and $f \in H^{s}(\mc{X})$, as long as $s \geq 2k + 1$. 
\begin{lemma}
	\label{lem:approximation_error_nonlocal_laplacian}
	Assume Model~\ref{def:model_flat_euclidean}. Let $s \in \mathbb{N} \setminus \{0,1\}$, suppose that $f \in H^s(\mc{X};M)$, and if $s > 1$ suppose that $p \in C^{s - 1}(\mc{X})$. Let $L_{P,\varepsilon}$ be define with respect to a kernel $\eta$ that satisfies~\ref{asmp:kernel_flat_euclidean}. Then there exist constants $C_1$ and $C_2$ that do not depend on $f$, such that each of the following statements hold.
	\begin{itemize}
		\item If $s$ is odd and $k = (s - 1)/2$, then
		\begin{equation}
		\label{eqn:approximation_error_nonlocal_laplacian_1}
		\|L_{P,\varepsilon}^kf - \Delta_P^kf\|_{L^2(\mc{X}_{k\varepsilon})} \leq C_1 M \varepsilon
		\end{equation}
		\item If $s$ is even and $k = (s - 2)/2$, then
		\begin{equation}
		\label{eqn:approximation_error_nonlocal_laplacian_2}
		\|L_{P,\varepsilon}^kf - \Delta_P^kf\|_{L^2(\mc{X}_{k\varepsilon})} \leq C_2 M \varepsilon^2.
		\end{equation}
	\end{itemize}
\end{lemma}
We remark that when $k = 1$ and $f \in C^3(\mc{X})$ or $C^4(\mc{X})$, statements of this kind are well known, and indeed stronger results---with $L^{\infty}(\mc{X})$ norm replacing $L^2(\mc{X})$ norm---hold. When dealing with the iterated Laplacian, and functions $f$ which are regular only in the Sobolev sense, the proof is somewhat more lengthy, but in result is similar in spirit.
 
\paragraph{Proof (of Lemma~\ref{lem:approximation_error_nonlocal_laplacian}).}
Throughout this proof, we shall assume that $f$ and $p$ are smooth functions, meaning they belong to $C^{\infty}(\mc{X})$. This is without loss of generality, since $C^{\infty}(\mc{X})$ is dense in both $H^s(\mc{X})$ and $C^{s - 1}(\mc{X})$, and since both sides of the inequalities~\eqref{eqn:approximation_error_nonlocal_laplacian_1} and~\eqref{eqn:approximation_error_nonlocal_laplacian_2} are continuous with respect to $\|\cdot\|_{H^s(\mc{X})}$ and $\|\cdot\|_{C^{s - 1}(\mc{X})}$ norms.

We will actually prove a more general set of statements than contained in Lemma~\ref{lem:approximation_error_nonlocal_laplacian}, more general in the sense that they give estimates for all $k$, rather than simply the particular choices of $k$ given above. In particular, we will prove that the following two statements hold for any $s \in \mathbb{N}$ and any $k \in \mathbb{N} \setminus \{0\}$. 
\begin{itemize}
	\item If $k \geq s/2$, then for every $x \in \mc{X}_{k\varepsilon}$, 
	\begin{equation}
	\label{pf:approximation_error_nonlocal_laplacian_0}
	L_{P,\varepsilon}^kf(x) = g_s(x) \varepsilon^{s - 2k}
	\end{equation}
	for a function $g_s$ that satisfies
	\begin{equation}
	\label{pf:approximation_error_nonlocal_laplacian_0.5}
	\|g_s\|_{L^2(\mc{X}_{k\varepsilon})} \leq C \|p\|_{C^{q}(\mc{X})}^k M 
	\end{equation}
	where $q = 1$ if $s =0$ or $s = 1$, and otherwise $q = s - 1$. 
	\item If $k < s/2$, then for every $x \in \mc{X}_{k\varepsilon}$,
	\begin{equation}
	\label{pf:approximation_error_nonlocal_laplacian_1}
	L_{P,\varepsilon}^kf(x) = \sigma_{\eta}^k \cdot \Delta_{P}^kf(x) + \sum_{j = 1}^{\floor{(s - 1)/2} - k} g_{2(j + k)}(x)\varepsilon^{2j} + g_{s}(x) \varepsilon^{s - 2k}.
	\end{equation}
	for functions $g_j$ that satisfy
	\begin{equation}
	\label{pf:approximation_error_nonlocal_laplacian_1.5}
	\|g_j\|_{H^{s - j}(\mc{X}_{k\varepsilon})} \leq C \|p\|_{C^{s - 1}(\mc{X})}^k M.
	\end{equation}
\end{itemize}
In the statement above, recall that $H^0(\mc{X}_{k\varepsilon}) = L^2(\mc{X}_{k\varepsilon})$. Additionally, note that we may speak of the pointwise behavior of derivatives of $f$ because we have assumed that $f$ is a smooth function. Observe that~\eqref{eqn:approximation_error_nonlocal_laplacian_1} follows upon taking $k = \floor{(s - 1)/2}$ in~\eqref{pf:approximation_error_nonlocal_laplacian_1}, whence we have
\begin{equation*}
\bigl(L_{P,\varepsilon}^kf(x) - \sigma_{\eta}^k \Delta_{P}^kf(x)\bigr)^2 = \varepsilon^2 \bigl(g_s(x)\bigr)^2
\end{equation*}
for some $g_s \in L^2(\mc{X}_{k\varepsilon},C \cdot M \cdot \|p\|_{C^{s - 1}(\mc{X})})$, and integrating over $\mc{X}_{k\varepsilon}$ gives the desired result. \eqref{eqn:approximation_error_nonlocal_laplacian_2} follows from~\eqref{pf:approximation_error_nonlocal_laplacian_1} in an identical fashion. 

It thus remains establish~\eqref{pf:approximation_error_nonlocal_laplacian_1}, and~\eqref{pf:approximation_error_nonlocal_laplacian_0} which is an important part of proving~\eqref{pf:approximation_error_nonlocal_laplacian_1}. We will do so by induction on $k$. Note that throughout, we will let $g_j$ refer to functions which may change from line to line, but which always satisfy~\eqref{pf:approximation_error_nonlocal_laplacian_1.5}. 

\underline{\textit{Proof of~\eqref{pf:approximation_error_nonlocal_laplacian_0} and~\eqref{pf:approximation_error_nonlocal_laplacian_1}, base case.}}

We begin with the base case, where $k = 1$. Again, we point out that although desired result is known when $s = 3$ or $s = 4$, and $f$ is regular in the H\"{o}lder sense, we require estimates for all $s \in \mathbb{N}$ when $f$ is regular in the Sobolev sense.

When $s = 0$, the inequality~\eqref{pf:approximation_error_nonlocal_laplacian_0} is implied by Lemma~\ref{lem:l2estimate_nonlocal_laplacian}.  When $s \geq 1$, we proceed using Taylor expansion. For any $x \in \mc{X}_{\varepsilon}$, we have that $B(x,\varepsilon) \subseteq \mc{X}$. Thus for any $x' \in B(x,\varepsilon)$, we may take an order $s$ Taylor expansion of $f$ around $x' = x$, and an order $q$ Taylor expansion of $p$ around $x' = x$, where $q = 1$ if $s = 1$, and otherwise $q = s - 1$. (See Section~\ref{subsec:taylor_expansion} for a review of the notation we use for Taylor expansions, as well as some properties that we make use of shortly.) This allows us to express $L_{P,\varepsilon}f(x)$ as the sum of three terms,
\begin{align*}
L_{P,\varepsilon}f(x) & = \frac{1}{\varepsilon^{d + 2}}\sum_{j_1 = 1}^{s - 1} \sum_{j_2 = 0}^{q - 1}\frac{1}{j_1!j_2!}  \int_{\mc{X}} \bigl(d_x^{j_1}f\bigr)(x' - x) \bigl(d_x^{j_2}p\bigr)(x' - x) \eta\biggl(\frac{\|x' - x\|}{\varepsilon}\biggr) \,dx' \quad + \\
& \quad \frac{1}{\varepsilon^{d + 2}}\sum_{j = 1}^{s - 1} \frac{1}{j!} \int_{\mc{X}} \bigl(d_x^jf\bigr)(x' - x)  r_{x'}^{q}(x;p) \eta\biggl(\frac{\|x' - x\|}{\varepsilon}\biggr) \,dx' \quad  + \\
& \quad \frac{1}{\varepsilon^{d + 2}} \int_{\mc{X}} r_{x'}^j(x;f) \eta\biggl(\frac{\|x' - x\|}{\varepsilon}\biggr) \,dP(x').
\end{align*}
Here we have adopted the convention that $\sum_{j = 1}^{0} = 0$. 

Changing variables to $z = (x' - x)/\varepsilon$, we can rewrite the above expression as 
\begin{align*}
L_{P,\varepsilon}f(x) & = \frac{1}{\varepsilon^{2}}\sum_{j_1 = 1}^{s - 1} \sum_{j_2 = 0}^{q - 1}\frac{\varepsilon^{j_1 + j_2}}{j_1!j_2!}  \int d_x^{j_1}f(z) d_x^{j_2}p(z) \eta\bigl(\|z\|\bigr) \,dz \quad + \\
& \quad \frac{1}{\varepsilon^{2}} \sum_{j = 1}^{s - 1} \frac{\varepsilon^j}{j!} \int d_x^jf(z)  r_{zh + x}^{q}(x;p) \eta\bigl(\|z\|\bigr) \,dz \quad  + \\
& \quad \frac{1}{\varepsilon^{2}} \int r_{zh + x}^j(x;f) \eta\bigl(\|z\|\bigr) p(zh + x)\,dz \\
& := G_1(x) + G_2(x) + G_3(x).
\end{align*}
We now separately consider each of $G_1(x),G_2(x)$ and $G_3(x)$. We will establish that if $s = 1$ or $s = 2$, then $G_1(x) = 0$, and otherwise if $s \geq 3$ that
\begin{equation*}
G_1(x) = \sigma_{\eta}\Delta_Pf(x) + \sum_{j = 1}^{\floor{(s - 1)/2} - 1}g_{2(j + 1)}(x)\varepsilon^{2j} + g_{s}(x)\varepsilon^{s - 2}.
\end{equation*}
On the other hand, we will establish that if $s = 1$ then $G_2(x) = 0$, and otherwise for $s \geq 2$
\begin{equation}
\label{pf:approximation_error_nonlocal_laplacian_2}
\|G_2\|_{L^2(\mc{X}_{\varepsilon})} \leq C \varepsilon^{s - 2} M \|p\|_{C^{s - 1}(\mc{X})};
\end{equation}
this same estimate will hold for $G_3$ for all $s \geq 1$. Together these will imply~\eqref{pf:approximation_error_nonlocal_laplacian_0} and~\eqref{pf:approximation_error_nonlocal_laplacian_1}. 

\emph{Estimate on $G_1(x)$.}
If $s = 1$, then $s - 1 = 0$, and so $G_1(x) = 0$. We may therefore suppose $s \geq 2$. Recall that
\begin{equation}
G_1(x) = \sum_{j_1 = 1}^{s - 1} \sum_{j_2 = 0}^{q - 1} \frac{\varepsilon^{j_1 + j_2 - 2}}{j_1!j_2!}  \underbrace{\int_{B(0,1)} d_x^{j_1}f(z) d_x^{j_2}p(z) \eta(\|z\|) \,dz}_{:= g_{j_1,j_2}(x)} \label{pf:approximation_error_nonlocal_laplacian_3}
\end{equation}
The nature of $g_{j_1,j_2}(x)$ depends on the sum $j_1 + j_2$. Since $d_x^{j_1}f d_x^{j_2}$ is an order $j_1 + j_2$ (multivariate) monomial, we have (see Section~\ref{subsec:taylor_expansion}) that whenever $j_1 + j_2$ is odd,
\begin{equation*}
g_{j_1,j_2}(x) = \int_{\mc{X}} d_x^{j_1}f(z) d_x^{j_2}p(z) \eta(\|z\|) \,dz = 0.
\end{equation*}
In particular this is the case when $j_1 = 1$ and $j_2 = 0$. Thus when $s = 2$,  $G_1(x) = g_{1,0}(x) = 0$. On the other hand if $s \geq 3$, then the lowest order terms in~\eqref{pf:approximation_error_nonlocal_laplacian_3} are those where $j_1 + j_2 = 2$, so that either $j_1 = 1$ and $j_2 = 1$, or $j_1 = 2$ and $j_2 = 0$. We have that
\begin{align*}
g_{1,1}(x) + \frac{1}{2}g_{2,0}(x) & = \int_{\mc{X}} d_x^{1}f(z) d_x^{1}p(z) \eta(\|z\|) \,dz + \frac{p(x)}{2} \int_{\mc{X}} d_x^{2}f(z) \eta(\|z\|) \,dz \\
& = \sum_{i_1 = 1}^{d} \sum_{i_2 = 1}^{d} D^{e_{i_1}}f(x) D^{e_{i_2}}p(x) \int_{\mc{X}} z^{e_{i_1} + e_{i_2}} \eta(\|z\|) \,dz + \frac{p(x)}{2} \sum_{i_1 = 1}^{d} \sum_{i_2 = 1}^{d}D^{e_{i_2}+e_{i_2}}f(x)\int_{\mc{X}} z^{e_{i_1} + e_{i_2}} \eta(\|z\|) \,dz\\
& = \sum_{i = 1}^{d} D^{e_{i}}f(x) D^{e_{i}}p(x) \int_{\mc{X}} z^2 \eta(\|z\|) \,dz + \frac{p(x)}{2} \sum_{i = 1}^{d} D^{2e_{i}}f(x)\int_{\mc{X}} z^2 \eta(\|z\|) \,dz\\ 
& = \sigma_{\eta}\Delta_Pf(x),
\end{align*}
which is the leading term order term. Now it remains only to deal with the higher-order terms, where $j_1 + j_2 > 2$, and where it suffices to show that each function $g_{j_1,j_2}$ satisfies~\eqref{pf:approximation_error_nonlocal_laplacian_1.5} for $j = \min\{j_1 + j_2 - 2,s - 2\}$. It is helpful to write $g_{j_1,j_2}$ using multi-index notation, 
\begin{align*}
g_{j_1,j_2}(x) = \sum_{|\alpha_1| = j_1} \sum_{|\alpha_2| = j_2} D^{\alpha_1}f(x) D^{\alpha_2}p(x) \int_{B(0,1)} z^{\alpha_1 + \alpha_2} \eta(\|z\|) \,dz,
\end{align*}
where we note that $|\int_{B(0,1)} z^{\alpha_1 + \alpha_2} \eta(\|z\|) \,dz| < \infty$ for all $\alpha_1, \alpha_2$, by the assumption that $\eta$ is Lipschitz on its support. Finally, by H\"{o}lder's inequality we have that
\begin{align*}
\|D^{\alpha_1}f D^{\alpha_2}p\|_{H^{s - (j + 2)}(\mc{X})} & \leq \|D^{\alpha_1}f\|_{H^{s - (j + 2)}(\mc{X})} \|D^{\alpha_2}p\|_{C^{s - (j + 2)}(\mc{X})} \\
& \leq \|D^{\alpha_1}f\|_{H^{s - j_1}(\mc{X})} \|D^{\alpha_2}p\|_{C^{s - (j_2 + 1)}(\mc{X})} \\
& \leq M \cdot \|p\|_{C^{s - 1}(\mc{X})},
\end{align*}
and summing over all $|\alpha_1| = j_1$ and $|\alpha_2| = j_2$ establishes that $g_{j_1,j_2}$ satisfies~\eqref{pf:approximation_error_nonlocal_laplacian_1.5}.

\emph{Estimate on $G_2(x)$.}
Note immediately that $G_2(x) = 0$ if $s = 1$. Otherwise if $s \geq 2$, then $q = s - 1$. Recalling that $|r_{x + z\varepsilon}^{s - 1}(x; p)| \leq C\varepsilon^{s - 1}\|p\|_{C^{s - 1}(\mc{X})}$ for any $z \in B(0,1)$, and that $d_x^jf(\cdot)$ is a $j$-homogeneous function, we have that
\begin{align}
|G_2(x)| & \leq \sum_{j = 1}^{s - 1} \frac{\varepsilon^{j - 2}}{j!}\int_{B(0,1)} \Bigl|\bigl(d_x^{j}f\bigr)(z)\Bigr| \cdot |r_{x + z\varepsilon}^{s - 1}(x;p)| \cdot \eta(\|z\|) \,dz \nonumber \\
& \leq C\varepsilon^{s - 2}\|p\|_{C^{s - 1}(\mc{X})} \sum_{j = 1}^{s - 1} \frac{1}{j!} \int_{B(0,1)} \Bigl|\bigl(d_x^{j}f\bigr)(z)\Bigr| \cdot \eta(\|z\|) \,dz \label{pf:approximation_error_nonlocal_laplacian_4}.
\end{align}
Furthermore, for each $j = 1,\ldots,s - 1$ convolution of $d_x^jf$ with $\eta$ only decreases the $L^2(\mc{X}_{\varepsilon})$ norm, meaning
\begin{equation}
\label{pf:approximation_error_nonlocal_laplacian_5}
\begin{aligned}
\int_{\mc{X}_{\varepsilon}} \biggl(\int_{B(0,1)} \Bigl|\bigl(d_x^{j}f\bigr)(z)\Bigr| \cdot \eta(\|z\|) \,dz\biggr)^2 \,dx & \leq \int_{\mc{X}_{\varepsilon}} \biggl(\int_{B(0,1)} \Bigl|\bigl(d_x^jf\bigr)(z)\Bigr|^2 \eta(\|z\|)\,dz \biggr) \cdot \biggl(\int_{B(0,1)} \eta(\|z\|) \,dz \biggr) \,dx \\
& \leq \int_{B(0,1)} \int_{\mc{X}_{\varepsilon}} \Bigl[\bigl(d^jf\bigr)(x)\Bigr]^2 \eta(\|z\|) \,dx  \,dz \\
& \leq \|d^jf\|_{L^2(\mc{X_{\varepsilon}})}^2.
\end{aligned}
\end{equation}
In the above, we have used both that $|d_x^jf(z)| \leq |d^jf(x)|$ for all $z \in B(0,1)$, and that the kernel is normalized so that $\int \eta(\|z\|) \,dz = 1$. 
Combining this with~\eqref{pf:approximation_error_nonlocal_laplacian_4}, we conclude that
\begin{align*}
\int_{\mc{X}_{\varepsilon}} |G_2(x)|^2 \,dx & \leq C \Bigl(\varepsilon^{s - 2}\|p\|_{C^{s - 1}(\mc{X})}\Bigr)^2 \sum_{j = 1}^{s - 1} \int_{\mc{X}_{\varepsilon}}\biggl(\frac{1}{j!} \int_{B(0,1)} \Bigl|\bigl(d_x^{j}f\bigr)(z)\Bigr| \cdot \Bigl|\eta(\|z\|)\Bigr| \,dz\biggr)^2 \,dx \\
& \leq C \Bigl(\varepsilon^{s - 2}\|p\|_{C^{s - 1}(\mc{X})}\Bigr)^2 \sum_{j = 1}^{s - 1} \|d^ju\|_{L^2(\mc{X_{\varepsilon}})}^2,
\end{align*}
establishing the desired estimate.

\emph{Estimate on $G_3(x)$.}
Applying the Cauchy-Schwarz inequality, we deduce a pointwise upper bound on $|G_3(x)|^2$,
\begin{align*}
|G_3(x)|^2 & \leq \biggl(\frac{p_{\max}}{\varepsilon^2}\biggr)^2 \cdot \biggl(\int_{B(0,1)} \bigl|r_{x + \varepsilon z}^s(x;u)\bigr|^2 \eta(\|z\|)\,dz\biggr) \cdot \biggl(\int_{B(0,1)} \eta(\|z\|) \,dz\biggr) \\
& \leq \biggl(\frac{p_{\max}}{\varepsilon^2}\biggr)^2 \int_{B(0,1)} \bigl|r_{x + \varepsilon z}^s(x;u)\bigr|^2 \eta(\|z\|) \,dz.
\end{align*}
Applying this pointwise over all $x \in \mc{X}_{\varepsilon}$ and integrating, we obtain
\begin{align*}
\int_{\mc{X}_{\varepsilon}} |G_3(x)|^2 \,dx & \leq \biggl(\frac{p_{\max}}{\varepsilon^2}\biggr)^2 \int_{\mc{X}_{\varepsilon}} \int_{B(0,1)} \bigl|r_{x + \varepsilon z}^s(x;f)\bigr|^2 \eta(\|z\|) \,dz \,dx \\
& = \biggl(\frac{p_{\max}}{\varepsilon^2}\biggr)^2 \int_{B(0,1)} \int_{\mc{X}_{\varepsilon}} \bigl|r_{x + \varepsilon z}^s(x;f)\bigr|^2 \eta(\|z\|) \,dx \,dz \\
& \leq \biggl(\frac{p_{\max}\varepsilon^s}{\varepsilon^2}\biggr)^2  \|d^sf\|_{L^2(\mc{X}_{\varepsilon})}^2,
\end{align*}
with the last inequality following from~\eqref{eqn:sobolev_remainder_term}. Noting that $p_{\max} = \|p\|_{C^0(\mc{X})} \leq \|p\|_{C^{s - 1}(\mc{X})}$, we see that this is a sufficient bound on $\|G_3\|_{L^2(\mc{X}_{\varepsilon})}$.

\underline{\textit{Proof of~\eqref{pf:approximation_error_nonlocal_laplacian_0} and~\eqref{pf:approximation_error_nonlocal_laplacian_1}, induction step.}}
We now assume that~\eqref{pf:approximation_error_nonlocal_laplacian_0} and~\eqref{pf:approximation_error_nonlocal_laplacian_1} hold for all order up to some $k$, and show that they then hold for order $k + 1$ as well. The proof is relatively straightforward, once we introduce a bit of notation. Namely, for any $\ell,j \in \mathbb{N}$ such that $1 \leq j \leq \ell \leq$, we will use $g_j^{\ell}$ to refer to a function satisfying
\begin{equation}
\label{pf:approximation_error_nonlocal_laplacian_6}
\|g_j^{\ell}\|_{H^{\ell - j}(\mc{X}_{(k + 1)\varepsilon})} \leq C \|p\|_{C^{q}(\mc{X})}^{k + 1} M.
\end{equation}
Note that $g_j^{\ell}(x) = g_{(s - \ell) + j}(x)$, so that $g_j^{s}(x) = g_j(x)$. As before, the functions $g_j^{\ell}$ may change from line to line, but will always satisfy~\eqref{pf:approximation_error_nonlocal_laplacian_6}. We immediately illustrate the purpose of this notation. Suppose $g \in H^{\ell}(\mc{X}_{k\varepsilon}; C \|p\|_{C^{q}(\mc{X})}^k M)$ for some $\ell \leq s$. If $\ell \leq 2$, then by the inductive hypothesis, it follows that for any $x \in \mc{X}_{(k + 1)\varepsilon}$
\begin{equation}
\label{pf:approximation_error_nonlocal_laplacian_7}
L_{P,\varepsilon}g(x) = g_{\ell}^{\ell}(x) \varepsilon^{\ell - 2}.
\end{equation} 
On the other hand if $2 < \ell \leq s$, then by the inductive hypothesis, it follows that for any $x \in \mc{X}_{(k + 1)\varepsilon}$,
\begin{equation}
\label{pf:approximation_error_nonlocal_laplacian_8}
L_{P,\varepsilon}g(x) = \sigma_{\eta} \Delta_Pg(x) + \sum_{j = 1}^{\floor{(\ell - 1)/2} - 1} g_{2j + 2}^{\ell}(x) \varepsilon^{2j} + g_{\ell}^{\ell}(x) \varepsilon^{\ell - 2}.
\end{equation}

\emph{Proof of \eqref{pf:approximation_error_nonlocal_laplacian_0}.} If $s \leq 2(k + 1)$, then by the inductive hypothesis it follows that for all $x \in \mc{X}_{k\varepsilon}$, we have $L_{P,\varepsilon}^kf(x) = g_{s}(x) \cdot \varepsilon^{s - 2k}$, for some $g_s \in L^2(\mc{X}_{k\varepsilon}, C\|p\|_{C^{s - 1}(\mc{X})}^k M)$. Note that we may know more about $L_P^kf(x)$ than simply that it is bounded in $L^2$-norm, but a bound in $L^2$-norm suffices. In particular, from such a bound along with~\eqref{pf:approximation_error_nonlocal_laplacian_7} we deduce that for any $x \in \mc{X}_{(k + 1)\varepsilon}$,
\begin{equation}
\label{pf:approximation_error_nonlocal_laplacian_8.5}
L_{P,\varepsilon}^{k + 1}f(x) = \bigl(L_{P,\varepsilon} \circ L_{P,\varepsilon}^k f)(x)= L_{P,\varepsilon} g_s(x)\varepsilon^{s - 2k} = g_{s}^{s}(x) \varepsilon^{s - 2(k + 1)},
\end{equation}
establishing~\eqref{pf:approximation_error_nonlocal_laplacian_0}. 

\emph{Proof of \eqref{pf:approximation_error_nonlocal_laplacian_1}.} If $s > 2(k + 1)$, then by the inductive hypothesis we have that for all $x \in \mc{X}_{k\varepsilon}$, 
\begin{equation*}
L_{P,\varepsilon}^kf(x) = \sigma_{\eta}^k \Delta_P^kf(x) + \sum_{j = 1}^{\floor{(s - 1)/2} - k} g_{2(j + k)}(x) \varepsilon^{2j} + g_s(x) \varepsilon^{s - 2k}.
\end{equation*}
Thus for any $x \in \mc{X}_{(k + 1)\varepsilon}$, 
\begin{equation*}
L_{P,\varepsilon}^{k + 1}f(x) = \bigl(L_{P,\varepsilon} \circ L_{P,\varepsilon}^k f\bigr)(x) = \sigma_{\eta}^k L_{P,\varepsilon}\Delta_P^kf(x) + \sum_{j = 1}^{\floor{(s - 1)/2} - k} L_{P,\varepsilon}g_{2(j + k)}(x) \varepsilon^{2j} + L_{P,\varepsilon}g_s(x) \varepsilon^{s - 2k}
\end{equation*}
There are three terms on the right hand side of this equality, and we now analyze each separately.
\begin{enumerate}
	\item Noting that $\Delta_P^kf \in H^{s - 2k}(\mc{X}; C\|p\|_{C^{s - 1}(\mc{X})}^kM)$, we use~\eqref{pf:approximation_error_nonlocal_laplacian_8} to derive that
	\begin{align}
	L_{P,\varepsilon}\Delta_P^kf(x) & = \sigma_{\eta} \Delta_P^{k + 1}f(x) + \sum_{j = 1}^{(s - 2k - 1)/2 - } g_{2j + 2}^{s - 2k}(x)\varepsilon^{2j} + g_{s - 2k}^{s - 2k}(x) \varepsilon^{s - 2k - 2} \nonumber \\
	& = \sigma_{\eta} \Delta_P^{k + 1}f(x) + \sum_{j = 1}^{(s - 1)/2 - (k + 1)} g_{2(k + 1 + j)}(x)\varepsilon^{2j} + g_{s}(x) \varepsilon^{s - 2(k + 1)}, \label{pf:approximation_error_nonlocal_laplacian_9}
	\end{align}
	where in the second equality we have simply used the fact $g_j^{\ell}(x) = g_{(s - \ell) + j}(x)$ to rewrite the equation.
	\item Suppose $j < \floor{(s - 1)/2} - k$. Then we use~\eqref{pf:approximation_error_nonlocal_laplacian_8} to derive that
	\begin{align*}
	L_{P,\varepsilon}g_{2(j + k)}(x) & = \sigma_{\eta}\Delta_P g_{2(j + k)}(x) + \sum_{i = 1}^{\floor{(s - 2j - 2k - 1)/2} - 1} g_{2(i + 1)}^{s - 2(j + k)}(x)\varepsilon^{2i} + g_{s - 2(j + k)}^{s - 2(j + k)}(x) \varepsilon^{s - 2(j + k + 1)} \\
	& = g_{2(j + k + 1)}(x) + \sum_{i = 1}^{\floor{(s - 1)/2} - (j + k + 1)} g_{2(i + j + k + 1)}(x)\varepsilon^{2i} + g_{s}(x) \varepsilon^{s - 2(j + k + 1)},
	\end{align*}
	where in the second equality we have again used $g_j^{\ell}(x) = g_{(s - \ell) + j}(x)$, and also written $\sigma_{\eta} \Delta_Pf = g_{2}^{s - 2(j + k)} = g_{2(j + k + 1)}$, since the particular dependence on the Laplacian $\Delta_P$ will not matter. From here, multiplying by $\varepsilon^{2j}$, we conclude that
	\begin{align}
	\varepsilon^{2j} L_{P,\varepsilon}g_{2(j + k)}(x) & = g_{2(j + k + 1)}(x) \varepsilon^{2j} + \sum_{i = 1}^{\floor{(s - 1)/2} - (j + k + 1)} g_{2(i + j + k + 1)}(x)\varepsilon^{2(i + j)} + g_{s}(x) \varepsilon^{s - 2(k + 1)} \nonumber \\ 
	& = g_{2(j + k + 1)}(x) \varepsilon^{2j} + \sum_{m = 1}^{\floor{(s - 1)/2} - (k + 1)}  g_{2(m + k + 1)}(x)\varepsilon^{2m} + g_{s}(x) \varepsilon^{s - 2(k + 1)} \label{pf:approximation_error_nonlocal_laplacian_10},
	\end{align}
	with the second equality following upon changing variables to $m = i + j$. 
	
	On the other hand if $j = \floor{(s - 1)/2} - k$, then the calculation is much simpler,
	\begin{equation}
	\label{pf:approximation_error_nonlocal_laplacian_11}
	\varepsilon^{2j} L_{P,\varepsilon}g_{2(j + k)}(x) = g_{s - 2(j + k)}^{s - 2(j + k)}(x) \varepsilon^{2j} \varepsilon^{s - 2(j + k) - 2} = g_s(x) \varepsilon^{s - 2(k + 1)}.
	\end{equation}
	\item Finally, it follows immediately from~\eqref{pf:approximation_error_nonlocal_laplacian_8} that
	\begin{equation}
	\label{pf:approximation_error_nonlocal_laplacian_12}
	L_{P,\varepsilon}g_s(x) \varepsilon^{s - 2k} = g_{s}(x) \varepsilon^{s - 2(k + 1)}.
	\end{equation}
\end{enumerate}
Plugging~\eqref{pf:approximation_error_nonlocal_laplacian_9}-\eqref{pf:approximation_error_nonlocal_laplacian_12} back into~\eqref{pf:approximation_error_nonlocal_laplacian_8.5} proves the claim.

\subsection{Boundary behavior of non-local Laplacian}
\label{subsec:boundary_behavior_nonlocal_laplacian}

In Lemma~\ref{lem:approximation_error_nonlocal_laplacian_boundary}, we establish that if $f$ is Sobolev smooth of order $s > 2k$ and zero-trace, then near the boundary of $\mc{X}$ the non-local Laplacian $L_{P,\varepsilon}^kf$ is close to $0$ in the $L^2$-sense.
\begin{lemma}
	\label{lem:approximation_error_nonlocal_laplacian_boundary}
	Assume Model~\ref{def:model_flat_euclidean}. Let $s,k \in \mathbb{N}$. Suppose that $f \in H_0^{s}(\mc{X};M)$. Then there exist numbers $c,C > 0$ that do not depend on $M$, such that for all $\varepsilon < c$, 
	\begin{equation*}
	\|L_{P,\varepsilon}^kf\|_{L^2(\partial_{k\varepsilon}\mc{X})}^2 \leq C \varepsilon^{2(s - 2k)}M^2.
	\end{equation*}
\end{lemma}

\paragraph{Proof (of Lemma~\ref{lem:approximation_error_nonlocal_laplacian_boundary})}
Applying Lemma~\ref{lem:l2estimate_nonlocal_laplacian}, we have that
\begin{equation*}
\|L_{P,\varepsilon}^kf\|_{L^2(\partial_{k\varepsilon}(\mc{X}))}^2 \leq \frac{(Cp_{\max})^{2}}{\varepsilon^4} \|L_{P,\varepsilon}^{k - 1}f\|_{L^2(\partial_{k\varepsilon}(\mc{X}))}^2 \leq \cdots \leq \frac{(Cp_{\max})^{2}}{\varepsilon^{4k}} \|f\|_{L^2(\partial_{k\varepsilon}(\mc{X}))}^2   
\end{equation*}
Thus it remains to show that for all $\varepsilon < c$,
\begin{equation}
\label{pf:approximation_error_nonlocal_laplacian_boundary_0}
\|f\|_{L^2(\partial_{k\varepsilon}(\mc{X}))}^2 = \int_{\partial_{k\varepsilon}(\mc{X})} \bigl(f(x)\bigr)^2 \,dx \leq C_1 \varepsilon^{2s} \|f\|_{H^s(\mc{X})}^2.
\end{equation}
We will build to~\eqref{pf:approximation_error_nonlocal_laplacian_boundary_0} by a series of intermediate steps, following the same rough structure as the proof of Theorem 18.1 in \citet{leoni2017}. For simplicity, we will take $k = 1$; the exact same proof applies to the general case upon assuming $\varepsilon < c/k$.

\underline{\textit{Step 1: Local Patch.}}
To begin, we assume that for some $c_0 > 0$ and a Lipschitz mapping $\phi: \Reals^{d - 1} \to [-c_0,c_0]$, we have that $f \in C_c^{\infty}(U_{\phi}(c_0))$, where 
\begin{equation*}
U_{\phi}(c_0) = \Bigl\{y \in Q(0,c_0): \phi(y_{-d}) \leq y_d\Bigr\}, 
\end{equation*}
and here $Q(0,c_0)$ is the $d$-dimensional cube of side length $c_0$, centered at $0$. We will show that for all $0 < \varepsilon < c_0$, and for the tubular neighborhood $V_{\phi}(\varepsilon) = \{y \in Q(0,c_0): \phi(y_{-d}) \leq y_d \leq \phi(y_{-d}) + \varepsilon\}$, we have that
\begin{equation*}
\int_{V_{\phi}(\varepsilon)} |f(x)|^2 \,dx \leq C\varepsilon^{2s} \|f\|_{H^s(U_{\phi}(c_0))}^2.
\end{equation*}
For a given $y = (y',y_d) \in V_{\phi}(\varepsilon)$, let $y_0 = (y',\phi(y'))$. Taking the Taylor expansion of $f(y)$ around $y = y_0$ because $u$ is compactly supported in $V_{\phi}$ it follows that,
\begin{align*}
f(y) & = f(y_0) + \sum_{j = 1}^{s - 1} \frac{1}{j!} D^{je_d}f(y_0) \bigl(y_d - \phi(y')\bigr)^j + \frac{1}{(s - 1)!}\int_{\phi(y')}^{y_d} (1 - t)^{s - 1} D^{se_d}f(y',z) \bigl(y_d - z\bigr)^{s - 1} \,dz \Longrightarrow\\
|f(y)| & \leq C\varepsilon^{s - 1}\int_{\phi(y')}^{y_d} \bigl|D^{se_d}f(y',z)\bigr| \,dz. 
\end{align*}
Consequently, by squaring both sides and applying Cauchy-Schwarz, we have that
\begin{equation*}
|f(y)|^2 \leq C\varepsilon^{2(s - 1)} \biggl(\int_{\phi(y')}^{y_d} \bigl|D^{se_d}f(y',z)\bigr| \,dz\biggr)^2 \leq C\varepsilon^{2s - 1} \int_{\phi(y')}^{y_d} \bigl|D^{se_d}f(y',z)\bigr|^2 \,dz.
\end{equation*}
Applying this bound for each $y \in V_{\phi}(\varepsilon)$, and then integrating, we obtain
\begin{align}
\int_{V_{\phi}(\varepsilon)} |f(y)|^2 \,dy & \leq \int_{Q_{d - 1}(c_0)} \int_{\phi(y')}^{\phi(y') + \varepsilon} |f(y',y_d)|^2 \,dy_d \,dy' \nonumber \\
& \leq C\varepsilon^{2s - 1}\int_{Q_{d - 1}(c_0)}  \int_{\phi(y')}^{\phi(y') + \varepsilon} \int_{\phi(y')}^{y_d} \bigl|D^{se_d}f(y',z)\bigr|^2 \,dz \,dy_d \,dy' \label{pf:approximation_error_nonlocal_laplacian_boundary_1}
\end{align}
where we have written $Q_{d - 1}(0,c_0)$ for the $d - 1$ dimensional cube of side length $c_0$, centered at $0$. Exchanging the order of the inner two integrals then gives
\begin{align*}
\int_{\phi(y')}^{\phi(y') + \varepsilon} \int_{\phi(y')}^{y_d} \bigl|D^{se_d}f(y',z)\bigr|^2 \,dz \,dy_d & = \int_{\phi(y')}^{\phi(y') + \varepsilon} \int_{z}^{\varepsilon} \bigl|D^{se_d}f(y',z)\bigr|^2 \,dy_d \,dz \\
& \leq C \varepsilon \int_{\phi(y')}^{\phi(y') + \varepsilon} \bigl|D^{se_d}f(y',z)\bigr|^2 \,dz \\
& \leq C \varepsilon \int_{\phi(y')}^{c_0} \bigl|D^{se_d}f(y',z)\bigr|^2 \,dz.
\end{align*}
Finally, plugging back into~\eqref{pf:approximation_error_nonlocal_laplacian_boundary_1}, we conclude that
\begin{equation*}
\int_{V_{\phi}(\varepsilon)} |f(y)|^2 \,dy \leq C \varepsilon^{2s} \int_{Q_{d - 1}(0,c_0)} \int_{\phi(y')}^{c_0} \bigl|D^{se_d}f(y',z)\bigr|^2 \,dz \,dy' \leq C \varepsilon^{2s} |u|_{H^s(U_{\phi}(c_0))}^2.
\end{equation*}

\underline{\textit{Step 2: Rigid motion of local patch.}} Now, suppose that at a point $x_0 \in \partial \mc{X}$, there exists a rigid motion $T: \Rd \to \Rd$ for which $T(x_0) = 0$, and a number $C_0$ such that for all $\varepsilon \cdot C_0 \leq c_0$, 
\begin{equation*}
T\bigl(Q_{T}(x_0,c_0) \cap \partial_{\varepsilon}\mc{X}\bigr) \subseteq V_{\phi}\bigl(C_0\varepsilon\bigr) \quad\textrm{and}\quad T\bigl(Q_T(x_0,c_0) \cap \mc{X}\bigr) = U_{\phi}(c_0).
\end{equation*}
Here $Q_{T}(x_0,c_0))$ is a (not necessarily coordinate-axis-aligned) cube of side length $c_0)$, centered at $x_0$. Define $v(y) := f(T^{-1}(y))$ for $y \in U_{\phi}(c_0)$. If $u \in C_c^{\infty}(\mc{X})$, then $v \in C_c^{\infty}(U_{\phi}(c_0))$, and moreover $\|v\|_{H^s(U_{\phi}(c_0))}^2 = \|f\|_{H^s(Q_{T}(x_0,c_0) \cap \mc{X})}^2$. Therefore, using the upper bound that we derived in Step 1,
\begin{equation*}
\int_{V_{\phi}(C_0 \cdot \varepsilon)} |v(y)|^2 \,dy \leq C \varepsilon^{2s} \|v\|_{H^s(U_{\phi}(c_0))}^2,
\end{equation*}
we conclude that
\begin{align*}
\int_{Q_{T}(x_0,c_0) \cap \partial_{\varepsilon}\mc{X}} |f(x)|^2 \,dx & = \int_{T(Q_T(x_0,c_0)) \cap \partial_{\varepsilon}\mc{X})} |v(y)|^2 \,dy \\
& \leq \int_{V_{\phi}(C_0 \cdot \varepsilon)} |v(y)|^2 \,dy \\
& \leq C \varepsilon^{2s} \|v\|_{H^s(U_{\phi}(c_0))}^2 = C \varepsilon^{2s} \|f\|_{H^s(Q_{T}(x_0,c_0)) \cap \mc{X})}^2 \leq C \varepsilon^{2s} \|f\|_{H^s(\mc{X})}^2.
\end{align*}

\underline{\textit{Step 3: Lipschitz domain}}.
Finally, we deal with the case where $\mc{X}$ is assumed to be an open, bounded subset of $\Rd$, with Lipschitz boundary. In this case, at every $x_0 \in \partial \mc{X}$, there exists a rigid motion $T_{x_0}: \Rd \to \Rd$ such that $T_{x_0}(x_0) = 0$, a number $c_0(x_0)$, a Lipschitz function $\phi_{x_0}:\Reals^{d - 1} \to [-c_0,c_0]$, and a number $C_0(x_0)$, such that for all $\varepsilon \cdot C_0(x_0) \leq c_0(x_0)$,
\begin{equation*}
T\bigl(Q_{T}(x_0,c_0(x_0)) \cap \partial_{\varepsilon}\mc{X}\bigr) \subseteq V_{\phi}\bigl(C_0(x_0) \cdot \varepsilon\bigr) \quad\textrm{and}\quad T\bigl(Q_T(x_0,c_0(x_0)) \cap \mc{X}\bigr) = U_{\phi}(c_0(x_0)).
\end{equation*}
Therefore for every $x_0 \in \partial \mc{X}$, it follows from the previous step that
\begin{equation*}
\int_{Q_{T_{x_0}}(x_0,c_0(x_0)) \cap \partial_{\varepsilon}\mc{X}} |f(x)|^2 \,dx \leq C(x_0) \varepsilon^{2s} \|f\|_{H^s(\mc{X})}^2,
\end{equation*}
where on the right hand side $C(x_0)$ is a constant that may depend on $x_0$, but not on $u$ or $\varepsilon$.

We conclude by taking a collection of cubes that covers $\partial_{\varepsilon}\mc{X}$ for all $\epsilon$ sufficiently small. First, we note that by a compactness argument there exists a finite subset of the collection of cubes $\{Q_{T_{x_0}}(x_0,c_0(x_0)/2): x_0 \in \partial\mc{X} \}$ which covers $\partial \mc{X}$, say $Q_{T_{x_1}}(x_1,c_0(x_1)/2),\ldots, Q_{T_{x_N}}(x_N,c_0(x_N)/2)$. Then, for any $\varepsilon \leq \min_{i = 1,\ldots,N} c_0(x_i)/2$, it follows from the triangle inequality that
\begin{equation*}
\partial_{\varepsilon}\mc{X} \subseteq \bigcup_{i = 1}^{N} Q_{T_{x_i}}(x_i, c_0(x_i)).
\end{equation*}
As a result,
\begin{equation*}
\int_{\partial_{\varepsilon}\mc{X}} |f(x)|^2 \leq \sum_{i = 1}^{N} \int_{Q_{T_{x_i}}(x_i, c_0(x_i)) \cap \partial_{\varepsilon}(\mc{X})} |f(x)|^2 \leq  \varepsilon^{2s} \|f\|_{H^s(\mc{X})}^2 \sum_{i = 1}^{N}C_0(x_i),
\end{equation*}
which proves the claim of~\eqref{pf:approximation_error_nonlocal_laplacian_boundary_0}.

\subsection{Estimate of non-local Sobolev seminorm}
\label{subsec:estimate_nonlocal_seminorm}

Now, we use the results of the preceding two sections to prove~\eqref{pf:graph_seminorm_ho_2}. We will divide our analysis in two cases, depending on whether $s$ is odd or even, but before we do this we state some facts that will be applicable to both cases. First, we  recall that $L_{P,\varepsilon}$ is self-adjoint in $L^2(P)$, meaning $\dotp{L_{P,\varepsilon}f}{g}_{P} = \dotp{f}{L_{P,\varepsilon}g}_{P}$ for all $f, g \in L^2(\mc{X})$. We also recall the definition of the Dirichlet energy $E_{P,\varepsilon}(f;\mc{X})$,
\begin{equation}
\label{eqn:dirichlet_energy}
\dotp{L_{P,\varepsilon}f}{f}_{P} = \frac{1}{\varepsilon^{d + 2}}\int_{\mc{X}} \int_{\mc{X}} \bigl(f(x) - f(x')\bigr)^2 \eta\biggl(\frac{\|x' - x\|}{\varepsilon}\biggr) \,dP(x') \,dP(x) =: E_{P,\varepsilon}(f;\mc{X}).
\end{equation}
Finally, we recall a result of~\cite{green2021}: there exist constants $c_0$ and $C_0$ which do not depend on $M$, such that for all $\varepsilon < c_0$ and for any $f \in H^1(\mc{X};M)$,
\begin{equation}
\label{pf:estimate_nonlocal_seminorm_1}
E_{P,\varepsilon}(f;\mc{X}) \leq C_0 M^2.
\end{equation}
\paragraph{Case 1: $s$ odd.}
Suppose $s$ is odd, so that $s \geq 3$. Taking $k = (s - 1)/2$, we use the self-adjointness of $L_{P,\varepsilon}$ to relate the non-local semi-norm $\dotp{L_{P,\varepsilon}^sf}{f}_{P}$ to a non-local Dirichlet energy,
\begin{equation*}
\dotp{L_{P,\varepsilon}^sf}{f}_P = \dotp{L_{P,\varepsilon}^{k + 1}f}{L_{P,\varepsilon}^{k}f}_P = E_{P,\varepsilon}(L_{P,\varepsilon}^{k}f;\mc{X}).
\end{equation*}
We now separate this energy into integrals over $\mc{X}_{k\varepsilon}$ and $\partial_{k\varepsilon}(\mc{X})$,
\begin{align}
E_{P,\varepsilon}(L_{P,\varepsilon}^{k}f;\mc{X}) & = \frac{1}{\varepsilon^{d + 2}}\Biggl\{\int_{\mc{X}_{k\varepsilon}} \int_{\mc{X}_{k\varepsilon}} \bigl(L_{P,\varepsilon}^kf(x) - L_{P,\varepsilon}^kf(x')\bigr)^2 \eta\biggl(\frac{\|x' - x\|}{\varepsilon}\biggr) \,dP(x') \,dP(x) \nonumber \\
& \quad + \int_{\partial_{k\varepsilon}\mc{X}} \int_{\partial_{k\varepsilon}\mc{X}} \bigl(L_{P,\varepsilon}^kf(x) - L_{P,\varepsilon}^kf(x')\bigr)^2 \eta\biggl(\frac{\|x' - x\|}{\varepsilon}\biggr) \,dP(x') \,dP(x)\Biggr\} \nonumber \\
& := E_{P,\varepsilon}(L_{P,\varepsilon}^{k}f;\mc{X}_{k\varepsilon}) + E_{P,\varepsilon}(L_{P,\varepsilon}^{k}f;\partial_{k\varepsilon}\mc{X}) \label{pf:estimate_nonlocal_seminorm_1.5}
\end{align}
and upper bound each energy separately. For the first term, we add and substract $\sigma_{\eta}^k\Delta_P^kf(x)$ and $\sigma_{\eta}^k\Delta_P^kf(x')$ within the integrand, then use the triangle inequality and the symmetry between $x$ and $x'$ to deduce that
\begin{equation}
\label{pf:estimate_nonlocal_seminorm_2}
E_{P,\varepsilon}(L_{P,\varepsilon}^{k}f;\mc{X}_{k\varepsilon}) \leq 3 \sigma_{\eta}^{2k} E_{P,\varepsilon}(\Delta_P^kf;\mc{X}_{k\varepsilon}) + \frac{2}{\varepsilon^{d + 2}}\int_{\mc{X}_{k\varepsilon}} \int_{\mc{X}_{k\varepsilon}} \bigl(L_{P,\varepsilon}^kf(x) - \sigma_{\eta}^k \Delta_P^kf(x)\bigr)^2 \eta\biggl(\frac{\|x' - x\|}{\varepsilon}\biggr) \,dP(x') \,dP(x).
\end{equation}
Noticing that $\Delta_P^kf \in H^1(\mc{X};\|p\|_{C^{s - 1}(\mc{X})}^kM)$, we use~\eqref{pf:estimate_nonlocal_seminorm_1} to conclude that $E_{P,\varepsilon}(\Delta_P^kf;\mc{X}_{k\varepsilon}) \leq C_0M^2$. On the other hand, it follows from Assumption~\ref{asmp:kernel_flat_euclidean} and~\eqref{eqn:approximation_error_nonlocal_laplacian_1} that
\begin{align*}
\frac{2}{\varepsilon^{d + 2}}\int_{\mc{X}_{k\varepsilon}} \int_{\mc{X}_{k\varepsilon}} \bigl(L_{P,\varepsilon}^kf(x) - \sigma_{\eta}^k \Delta_P^kf(x)\bigr)^2 \eta\biggl(\frac{\|x' - x\|}{\varepsilon}\biggr) \,dP(x') \,dP(x) & \leq \frac{2p_{\max}}{\varepsilon^{2}}\int_{\mc{X}_{k\varepsilon}} \bigl(L_{P,\varepsilon}^kf(x) - \sigma_{\eta}^k \Delta_P^kf(x)\bigr)^2 \,dP(x) \\
& \leq C_1M^2.
\end{align*}
Plugging these two bounds into~\eqref{pf:estimate_nonlocal_seminorm_2} gives the desired upper bound on $E_{P,\varepsilon}(L_{P,\varepsilon}^{k};\mc{X}_{k\varepsilon})$. 

For the second term in~\eqref{pf:estimate_nonlocal_seminorm_1.5}, we apply Lemmas~\ref{lem:dirichlet_estimate_nonlocal_laplacian} and~\ref{lem:approximation_error_nonlocal_laplacian_boundary} and conclude that,
\begin{equation*}
E_{P,\varepsilon}(L_{P,\varepsilon}^{k}f;\partial_{k\varepsilon}\mc{X}) \leq \frac{4p_{\max}^2}{\varepsilon^2}\|L_{P,\varepsilon}^kf\|_{L^2(\partial_{k\varepsilon}\mc{X})} \leq C M^2.
\end{equation*}
\paragraph{Case 2: $s$ even.}
If $s \in \mathbb{N}$ is even, $s \geq 2$, then letting $k = (s - 2)/2$, the  self-adjointness of $L_{P,\varepsilon}$ implies
\begin{equation*}
\dotp{L_{P,\varepsilon}^sf}{f}_P = \|L_{P,\varepsilon}^{k + 1}f\|_P^2.
\end{equation*}
As in the first case, we divide the integral up into the interior region $\mc{X}_{k\varepsilon}$ and the boundary region $\partial_{k\varepsilon}\mc{X}$,
\begin{equation}
\label{pf:estimate_nonlocal_seminorm_3}
\|L_{P,\varepsilon}^{k + 1}f\|_P^2 \leq p_{\max} \|L_{P,\varepsilon}^{k + 1}f\|_{L^2(\mc{X})}^2 \leq p_{\max}\biggl\{\int_{\mc{X}_{k\varepsilon}} \bigl(L_{P,\varepsilon}^{k + 1}f(x)\bigr)^2 \,dP(x) + \int_{\partial_{k\varepsilon}\mc{X}} \bigl(L_{P,\varepsilon}^{k + 1}f(x)\bigr)^2 \,dP(x)\biggr\},
\end{equation}
and upper bound each term separately. For the first term, adding and subtracting $\sigma_{\eta}^k \Delta_P^kf(x)$ gives
\begin{align*}
\int_{\mc{X}_{k\varepsilon}} \bigl(L_{P,\varepsilon}^{k + 1}f(x)\bigr)^2 \,dP(x) & \leq 2\int_{\mc{X}_{k\varepsilon}} \bigl(L_{P,\varepsilon}\Delta_P^kf(x)\bigr)^2 \,dP(x)  + 2 \int_{\mc{X}_{k\varepsilon}} \Bigl(L_{P,\varepsilon}\bigl(L_{P,\varepsilon}^kf- \sigma_{\eta}\Delta_P^kf\bigr)(x)\Bigr)^2 \,dP(x) \\
& \overset{(i)}{\leq} CM^2  + 2 \int_{\mc{X}_{k\varepsilon}} \Bigl(L_{P,\varepsilon}\bigl(L_{P,\varepsilon}^kf- \sigma_{\eta}\Delta_P^kf\bigr)(x)\Bigr)^2 \,dP(x) \\
& \overset{(ii)}{\leq} CM^2  + \frac{Cp_{\max}^2}{\varepsilon^{2}} \|L_{P,\varepsilon}^kf- \sigma_{\eta}\Delta_P^kf\|_{L^2(\mc{X}_{k\varepsilon})}^2 \\
& \overset{(iii)}{\leq} CM^2,
\end{align*}
with $(i)$ following from~\eqref{pf:approximation_error_nonlocal_laplacian_0} since $\Delta_P^kf \in H^2(\mc{X};M\|p\|_{C^{s - 1}(\mc{X})}^l)$, $(ii)$ following from Lemma~\ref{lem:l2estimate_nonlocal_laplacian}, and $(iii)$ following from~\eqref{eqn:approximation_error_nonlocal_laplacian_2}.

Then Lemma~\ref{lem:approximation_error_nonlocal_laplacian_boundary} shows that the second term in~\eqref{pf:estimate_nonlocal_seminorm_3} satisfies
\begin{equation*}
\int_{\partial_{k\varepsilon}\mc{X}} \bigl(L_{P,\varepsilon}^{k + 1}f(x)\bigr)^2 \,dP(x) \leq CM^2.
\end{equation*}

\subsection{Assorted integrals}
\label{subsec:integrals}

\begin{lemma}
	\label{lem:l2estimate_nonlocal_laplacian}
	Assume Model~\ref{def:model_flat_euclidean}. Suppose $f \in L^2(U;M)$ for a Borel set $U \subseteq \mc{X}$, and let $L_{P,\varepsilon}$ be defined with respect to a kernel $\eta$ that satisfies~\ref{asmp:kernel_flat_euclidean}. Then there exists a constant $C$ which does not depend on $f$ or $M$ such that
	\begin{equation}
	\label{eqn:l2estimate_nonlocal_laplacian}
	\|L_{P,\varepsilon}f\|_{L^2(U)} \leq \frac{2 p_{\max}}{\varepsilon^2} \|f\|_{L^2(U)}
	\end{equation}
\end{lemma}

\begin{lemma}
	\label{lem:dirichlet_estimate_nonlocal_laplacian}
	Assume Model~\ref{def:model_flat_euclidean}. Suppose $f \in L^2(U;M)$ for a Borel set $U \subseteq \mc{X}$, and let $L_{P,\varepsilon}$ be defined with respect to a kernel $\eta$ that satisfies~\ref{asmp:kernel_flat_euclidean}. Then there exists a constant $C$ which does not depend on $f$ or $M$ such that
	\begin{equation}
	\label{eqn:dirichlet_estimate_nonlocal_laplacian}
	E_{P,\varepsilon}(f;U) \leq \frac{4 p_{\max}^2}{\varepsilon^2} \|f\|_{L^2(U)}^2
	\end{equation}
\end{lemma}

\begin{lemma}
	\label{lem:graph_seminorm_bias2}
	Assume Model~\ref{def:model_flat_euclidean}. Suppose $f \in H^1(\mc{X};M)$, and let $D_if$ be defined with respect to a kernel $\eta$ that satisfies~\ref{asmp:kernel_flat_euclidean}. Then there exists a constant $C$ which does not depend on $f$ or $M$, such that for any $i \in [n]$ and $\bj \in [n]^s$,
	\begin{equation*}
	\Ebb\Bigl[|D_{\bj}f(X_i)| \cdot |f(X_i) - f(X_{\bj_1})|\Bigr] \leq C \varepsilon^{2 + dk}M^2,
	\end{equation*}
	where $k + 1$ is the number of distinct indices in $i\bj$. 
\end{lemma}

\paragraph{Proof (of Lemma~\ref{lem:l2estimate_nonlocal_laplacian}).}
We fix a version of $f \in L^2(U)$, so that we may speak of its pointwise values.

At a given point $x \in U$, we can upper bound $|L_{P,\varepsilon}f(x)|^2$ using the Cauchy-Schwarz inequality as follows,
\begin{align*}
|L_{P,\varepsilon}f(x)|^2 & \leq \biggl(\frac{p_{\max}}{\varepsilon^{2 + d}}\biggr)^2 \Biggl(\int_U \bigl(|f(x')| + |f(x)|\bigr)^2 \eta\biggl(\frac{\|x' - x\|}{\varepsilon}\biggr) \,dx'\Biggr)^2 \\
& \leq \biggl(\frac{p_{\max}}{\varepsilon^{2 + d}}\biggr)^2 \Biggl(\int_U \bigl(|f(x')| + |f(x)|\bigr)^2 \eta\biggl(\frac{\|x' - x\|}{\varepsilon}\biggr) \,dx' \cdot \int \eta\biggl(\frac{\|x' - x\|}{\varepsilon}\biggr) \,dx' \Biggr) \\
& = \frac{p_{\max}^2}{\varepsilon^{4 + d}} \int_U \bigl(|f(x')| + |f(x)|\bigr)^2 \eta\biggl(\frac{\|x' - x\|}{\varepsilon}\biggr) \,dx'.
\end{align*}
The equality follows by the assumption $\int_{\Rd} \eta(\|z\|)\,dx = 1$ in~\ref{asmp:kernel_flat_euclidean}. Integrating over all $x \in U$, it follows from the triangle inequality that
\begin{align}
\|L_{P,\varepsilon}\|_{L^2(U)}^2 & \leq \frac{2 p_{\max}^2}{\varepsilon^{4 + d}} \int_{U} \int_{U} \bigl(|f(x')|^2 + |f(x)|^2\bigr) \eta\biggl(\frac{\|x' - x\|}{\varepsilon}\biggr) \,dx' \,dx \nonumber \\
& \leq \frac{2 p_{\max}^2}{\varepsilon^{4 + d}} \int_{U} \int_{U} \bigl(|f(x')|^2 + |f(x)|^2\bigr) \eta\biggl(\frac{\|x' - x\|}{\varepsilon}\biggr) \,dx' \,dx. \label{pf:l2estimate_nonlocal_laplacian_1}
\end{align}
Finally, using Fubini's Theorem we determine that
\begin{equation}
\int_{U} \int_{U} \bigl(|f(x')|^2 + |f(x)|^2\bigr) \eta\biggl(\frac{\|x' - x\|}{\varepsilon}\biggr) \,dx' \,dx = 2 \int_{U} \int_{U} |f(x)|^2 \eta\biggl(\frac{\|x' - x\|}{\varepsilon}\biggr) \,dx \leq 2 \varepsilon^d \int_{U} |f(x)|^2 \,dx = 2\varepsilon^d \|f\|_{L^2(U)}^2,
\label{pf:l2estimate_nonlocal_laplacian_2}
\end{equation}
and by combining~\eqref{pf:l2estimate_nonlocal_laplacian_1} and~\eqref{pf:l2estimate_nonlocal_laplacian_2} we conclude that
\begin{equation*}
\|L_{P,\varepsilon}\|_{L^2(U)}^2 \leq \frac{4p_{\max}^2}{\varepsilon^4} \|f\|_{L^2(U)}^2. 
\end{equation*}

\paragraph{Proof (of Lemma~\ref{lem:dirichlet_estimate_nonlocal_laplacian}).}
We have
\begin{equation*}
E_{P,\varepsilon}(f) = \frac{1}{\varepsilon^{2 + d}} \int_{U} \int_{U} \bigl(f(x) - f(x')\bigr)^2 \eta\biggl(\frac{\|x' - x\|}{\varepsilon}\biggr) \,dP(x') \,dP(x) \leq \frac{2p_{\max}^2}{\varepsilon^{2 + d}} \int_{U} \int_{U} \bigl(|f(x)|^2 + |f(x')|^2\bigr) \eta\biggl(\frac{\|x' - x\|}{\varepsilon}\biggr) \,dx' \,dx,
\end{equation*}
and the claim follows from~\eqref{pf:l2estimate_nonlocal_laplacian_2}.

\paragraph{Proof (of Lemma~\ref{lem:graph_seminorm_bias2}).}
Let $G_{n,\varepsilon}[X_{i\bj}]$ be the subgraph induced by vertices $X_i, X_{\bj_1},\ldots,X_{\bj_s}$. We make two observations. First, in order for $|D_{\bj}f(X_i)| \cdot |f(X_i)  - f(X_j)|$ to be non-zero, it must be the case that the subgraph $G_{n,\varepsilon}[X_{i\bj}]$ is connected. Second, noting that for any indices $i$ and $j$,
\begin{equation*}
|D_{ij}f(x)| \leq \Bigl(|D_jf(X_i)| + |D_jf(x)|\Bigr)\|\eta\|_{\infty}, 
\end{equation*}
a straightforward inductive argument implies that 
\begin{equation*}
|D_{\bj}f(X_i)| \leq s\|\eta\|_{\infty}^s \sum_{j \in i\bj} |D_{\bj_s}f(X_j)|.
\end{equation*}
Combining these two observations, we reduce the task to upper bounding the product of two (first-order) differences,
\begin{align*}
\mathbb{E}\Bigl[ |D_{\bj}f(X_i)| |f(X_i) - f(X_{\bj_1})| \Bigr] & = \mathbb{E}\Bigl[ |D_{\bj}f(X_i)| |f(X_i) - f(X_{\bj_1})| \cdot \1\bigl\{G_{n,\varepsilon}[X_{i\bj}]~~\textrm{is connected.} \bigr\} \Bigr] \\
& \leq s\|\eta\|_{\infty}^s \sum_{j \in i\bj} \mathbb{E}\Bigl[ |D_{\bj_s}f(X_j)| \cdot |f(X_i) - f(X_{\bj_1})| \cdot \1\bigl\{G_{n,\varepsilon}[X_{i\bj}]~~\textrm{is connected.} \bigr\}   \Bigr] \\
& \leq s\|\eta\|_{\infty}^s \sum_{j \in i\bj} \mathbb{E}\Bigl[ |f(X_j) - f(X_{\bj_s})| \cdot |f(X_i) - f(X_{\bj_1})| \cdot \1\bigl\{G_{n,\varepsilon}[X_{i\bj}]~~\textrm{is connected.} \bigr\}   \Bigr] 
\end{align*}
Next, from the Cauchy-Schwarz inequality we have that for any $j \in \bj$,
\begin{align*}
& \mathbb{E}\Bigl[ |f(X_j) - f(X_{\bj_s})| \cdot |f(X_i) - f(X_{\bj_1})| \cdot \1\bigl\{G_{n,\varepsilon}[X_{i\bj}]~~\textrm{is connected.} \bigr\}   \Bigr] \\
& \quad \leq \sqrt{\mathbb{E}\Bigl[ |f(X_j) - f(X_{\bj_s})|^2 \cdot \1\bigl\{G_{n,\varepsilon}[X_{i\bj}]~~\textrm{is connected.} \bigr\} \Bigr]} \cdot \sqrt{\mathbb{E}\Bigl[ |f(X_j) - f(X_{\bj_s})|^2 \cdot \1\bigl\{G_{n,\varepsilon}[X_{i\bj}]~~\textrm{is connected.} \bigr\} \Bigr]} \\
& \quad = \mathbb{E}\Bigl[ |f(X_j) - f(X_{i})|^2 \cdot \1\bigl\{G_{n,\varepsilon}[X_{i\bj}]~~\textrm{is connected.} \bigr\} \Bigr],
\end{align*}
with the equality following since each $X_i$ are identically distributed. Marginalizing out the contribution of all indices in $\bj$ not equal to $i$ or $j$ gives
\begin{align}
\mathbb{E}\Bigl[ |f(X_j) - f(X_{i})|^2 \cdot \1\bigl\{G_{n,\varepsilon}[X_{i\bj}]~~\textrm{is connected.} \bigr\} \Bigr] & \leq \bigl((s + 1)p_{\max}\nu_d\varepsilon^d\bigr)^{|i\bj\setminus\{j \cup i\}|} \cdot \mathbb{E}\Bigl[ |f(X_j) - f(X_{i})|^2 \1\{\|X_i - X_j\| \leq \varepsilon\}\Bigr] \nonumber \\
& \leq \bigl((s + 1)p_{\max}\nu_d\varepsilon^d\bigr)^{|i\bj\setminus\{j \cup i\}|} \cdot p_{\max}^2 \nu_d \varepsilon^{2 + d} M^2 \label{pf:graph_seminorm_bias2_1}
\end{align}
with the second inequality following from the proof of Lemma~1 in~\cite{green2021}. Finally, we notice that $|i\bj\setminus \{i \cup j\}| + 1 = k$, so that~\eqref{pf:graph_seminorm_bias2_1} gives the desired result.

\section{Graph Sobolev semi-norm, manifold domain}
\label{sec:graph_quadratic_form_manifold}
In this section we prove Proposition~\ref{prop:graph_seminorm_manifold}. Note that when $s = 1$, the upper bound~\eqref{eqn:graph_seminorm_manifold} follows immediately from Lemma~\ref{lem:dirichlet_energy_sobolev} and Markov's inequality. 

On the other hand when $s = 2$ or $s = 3$, we prove Proposition~\ref{prop:graph_seminorm_manifold} by first establishing some intermediate results, many of which are analogous to results we have already shown in the flat Euclidean case. Indeed, in some ways the proof will be simpler in the manifold setting than in the flat Euclidean case: there is no boundary, and we do not need to analyze the iterated nonlocal Laplacian $L_{P,\varepsilon}^j$ for $j > 1$. 

That being said, as mentioned in our main text, in the manifold setting there is some extra error induced by using Euclidean rather than geodesic distance. We upper bound this error by comparing $L_{P,\varepsilon}$ to an alternative nonlocal Laplacian $\wt{L}_{P,\varepsilon}$, which is defined with respect to geodesic distance. Precisely, let $d_{\mc{X}}(x,x')$ denote the geodesic distance between $x,x' \in \mc{X}$, and define
\begin{equation*}
\wt{L}_{P,\varepsilon}f(x) := \int_{\mc{X}} \bigl(f(x') - f(x)\bigr) \eta \biggl(\frac{d_{\mc{X}}(x',x)}{\varepsilon}\biggr) p(x') \,dx'.
\end{equation*}

We show the following results, each of which hold under the same assumptions as Proposition~\ref{prop:graph_seminorm_manifold}.
\begin{itemize}
	\item In Section~\ref{subsec:manifold_decomposition_graph_seminorm} we show that the graph Sobolev seminorm $\dotp{L_{n,\varepsilon}^sf}{f}_n$ is upper bounded by the sum of a nonlocal seminorm and a pure bias term: specifically, with probability at least $1 - 2\delta$,
	\begin{equation}
	\label{pf:graph_seminorm_manifold_1}
	\dotp{L_{n,\varepsilon}^sf}{f}_n \leq \frac{\dotp{L_{P,\varepsilon}^sf}{f}_P}{\delta} + C_1\frac{\varepsilon^2}{n\varepsilon^{2s + m}}M^2.
	\end{equation}
	This upper bound is essentially the same as~\eqref{pf:graph_seminorm_ho_1}, but with the intrinsic dimension $m$ taking the place of the ambient dimension $d$. The pure bias term will be of at most constant order when $\varepsilon \gtrsim n^{-1/(2(s-1) + m)}$. 
	\item In Section~\ref{subsec:error_euclidean_distance}, we show that the error incurred by using the ``wrong'' metric is negligible. Precisely, we find that
	\begin{equation}
	\label{eqn:nonlocal_laplacian_geodesic_error}
	\|L_{P,\varepsilon}f - \wt{L}_{P,\varepsilon}f\|_{L^2(\mc{X})}^2 \leq C_2 \varepsilon^2 |f|_{H^1(\mc{X})}^2.
	\end{equation}
	\item In Section~\ref{subsec:manifold_approximation_error_nonlocal_laplacian}, we analyze the approximation error of $\wt{L}_{P,\varepsilon}$. We show that when $f \in H^2(\mc{X})$ and $p \in C^1(\mc{X})$, 
	\begin{equation}
	\label{eqn:nonlocal_laplacian_approximation_error_manifold_l2}
	\|\wt{L}_{P,\varepsilon}f\|_{L^2(\mc{X})}^2 \leq C_3 \|f\|_{H^2(\mc{X})}^2,
	\end{equation}
	whereas if $f \in H^3(\mc{X})$ and $p \in C^2(\mc{X})$, 
	\begin{equation}
	\label{eqn:nonlocal_laplacian_approximation_error_manifold_sobolev}
	\|\wt{L}_{P,\varepsilon}f - \sigma_{\eta}\Delta_Pf\|_{L^2(\mc{X})}^2 \leq C_3 \varepsilon^2 \|f\|_{H^3(\mc{X})}^2.
	\end{equation}
	\item In Section~\ref{subsec:manifold_estimate_nonlocal_seminorm}, we use the results of the preceding two sections to show that if $f \in H^s(\mc{X})$ and $p \in C^{s - 1}(\mc{X})$, then 
	\begin{equation}
	\label{eqn:manifold_nonlocal_seminorm}
	\dotp{L_{P,\varepsilon}^sf}{f}_P \leq C_4\|f\|_{H^s(\mc{X})}^2.
	\end{equation}
	\item In Section~\ref{subsec:manifold_integrals} we state some technical results used in the previous sections.
\end{itemize}
We point out that when $f$ is H\"{o}lder smooth, results analogous to~\eqref{eqn:nonlocal_laplacian_approximation_error_manifold_sobolev} have been established in \citet{calder2019}. When $f$ is Sobolev smooth, our analysis (which relies heavily on Taylor expansions) is largely similar, except that the remainder term in the relevant Taylor expansion will be bounded in $L^2(\mc{X})$ norm rather than $L^{\infty}(\mc{X})$ norm. This is analogous to the situation in the flat Euclidean model.

In the proof of~\eqref{pf:graph_seminorm_manifold_1}-\eqref{eqn:manifold_nonlocal_seminorm}, we recall the following estimates from differential geometry: (i) letting $K_0$ be an upper bound on the absolute value of the sectional curvatures of $\mc{X}$, $K_0 \leq 2R$, and  letting (ii) $i_0$ be a lower bound on the injectivity radius of $\mc{X}$, $i_0 \geq \pi R$; see Proposition~1 of \cite{aamari2019}. Additionally, recall that for all $\delta < i_0$, the exponential map $\exp_x: B_m(0,\delta) \subset T_x(\mc{X}) \to B_{\mc{X}}(x,\delta) \subset \mc{X}$ is a diffeomorphism for all $x \in \mc{X}$. We shall therefore always assume $\varepsilon < i_0$. 

\paragraph{Proof (of Proposition~\ref{prop:graph_seminorm_manifold}).} Follows immediately from~\eqref{pf:graph_seminorm_manifold_1} and~\eqref{eqn:manifold_nonlocal_seminorm}. \qed

\subsection{Decomposition of graph Sobolev seminorm}
\label{subsec:manifold_decomposition_graph_seminorm}
The proof of~\eqref{pf:graph_seminorm_manifold_1} is identical to the proof of~\eqref{pf:graph_seminorm_ho_1}, except substituting the intrinsic dimension $m$ for ambient dimension $d$, and using Lemma~\ref{lem:manifold_graph_seminorm_bias2} rather than Lemma~\ref{lem:graph_seminorm_bias2}.

\subsection{Error due to Euclidean Distance}
\label{subsec:error_euclidean_distance}

In this section, we prove~\eqref{eqn:nonlocal_laplacian_geodesic_error}. By applying Cauchy-Schwarz we obtain an upper bound on $|L_{P,\varepsilon}f(x) - \wt{L}_{P,\varepsilon}f(x)|^2$:
\begin{align}
\bigl[L_{P,\varepsilon}f(x) - \wt{L}_{P,\varepsilon}f(x)\bigr]^2 & \leq \frac{p_{\max}^2}{\varepsilon^{2(2 + m)}} \int_{\mc{X}} \bigl[f(x') - f(x)\bigr]^2 \biggl|\eta\biggl(\frac{\|x' - x\|}{\varepsilon}\biggr) - \eta\biggl(\frac{d_{\mc{X}}(x',x)}{\varepsilon}\biggr)\biggr| \,d\mu(x') \nonumber \\
& \quad \cdot \int_{\mc{X}} \biggl|\eta\biggl(\frac{\|x' - x\|}{\varepsilon}\biggr) - \eta\biggl(\frac{d_{\mc{X}}(x',x)}{\varepsilon}\biggr)\biggr| \,d\mu(x') \nonumber \\
& = \frac{1}{\varepsilon^{2(2 + m)}} A_1(x) \cdot A_2(x) \label{pf:error_euclidean_distance_0}
\end{align} 
Thus we have upper bounded $|L_{P,\varepsilon}f(x) - \wt{L}_{P,\varepsilon}f(x)|^2$ by the product of two terms, each of which we now suitably bound. To do so, we will use the following estimate, from Proposition 4 of \cite{trillos2019}: for all $\|x' - x\| \leq R/2$,
\begin{equation}
\label{eqn:distance_error}
\|x' - x\| \leq d_{\mc{X}}(x',x) \leq \|x' - x\| + \frac{8}{R^2} \|x' - x\|^3.
\end{equation}
From here forward we will assume $\varepsilon < R/2$. 
\paragraph{Upper bound on $A_1(x)$.}
Consequently $\eta(\|x' - x\|/\varepsilon) \geq \eta(d_{\mc{X}}(x',x)/\varepsilon)$. Furthermore, letting $L_{\eta}$ denote the Lipschitz constant of $\eta$, and setting $\wt{\varepsilon} := (1 + 27\varepsilon^2/R^2)\varepsilon$ we have that
\begin{equation*}
\biggl|\eta\biggl(\frac{\|x' - x\|}{\varepsilon}\biggr) - \eta\biggl(\frac{d_{\mc{X}}(x',x)}{\varepsilon}\biggr)\biggr| \leq \frac{L_{\eta} 8 \varepsilon^2}{R^2} \cdot \1\bigl\{d_{\mc{X}}(x',x) \leq \varepsilon\bigr\} + \|\eta\|_{\infty} \cdot 1\{\varepsilon < d_{\mc{X}}(x',x) \leq \wt{\varepsilon}\}.
\end{equation*}
Thus,
\begin{equation*}
A_1(x) \leq \frac{8L_{\eta}\varepsilon^2}{R^2}\int_{\mc{X}}\bigl[f(x') - f(x)\bigr]^2 \1\{\|x' - x\| \leq \varepsilon\} \,d\mu(x') + \|\eta\|_{\infty} \int_{\mc{X}}\bigl[f(x') - f(x)\bigr]^2 \1\bigl\{\varepsilon < d_{\mc{X}}(x',x) \leq \wt{\varepsilon}\bigr\} \,d\mu(x') \\
\end{equation*}
Integrating over $\mc{X}$, we conclude from Lemma~\ref{lem:dirichlet_energy_remainder} and Lemma~3.3 of \citep{burago2014} and  that
\begin{equation*}
\int_{\mc{X}} A_1(x) \,d\mu(x) \leq \frac{8L_{\eta}\nu_m\varepsilon^2}{R^2(m + 2)} \Bigl(1 + CmK_0R^2\Bigr) \varepsilon^{m + 2} |f|_{H^1(\mc{X})}^2 + C\|\eta\|_{\infty}\varepsilon^{m + 4} |f|_{H^1(\mc{X})}^2 =: C_5 \varepsilon^{m + 4}|f|_{H^1(\mc{X})}.
\end{equation*}

\paragraph{Upper bound on $A_2(x)$.}
Integrating over $x' \in \mc{X}$, we see that
\begin{align}
\int_{\mc{X}} \biggl|\eta\biggl(\frac{\|x' - x\|}{\varepsilon}\biggr) - \eta\biggl(\frac{d_{\mc{X}}(x',x)}{\varepsilon}\biggr)\biggr| \,d\mu(x') & \leq \frac{8L_{\eta}\varepsilon^2}{R^2} \int_{\mc{X}} \1\bigl\{d_{\mc{X}}(x',x)\bigr\} \,d\mu(x') + p_{\max}\|\eta\|_{\infty} \int_{\mc{X}} \1\bigl\{ \varepsilon < d_{\mc{X}}(x',x) \leq \wt{\varepsilon} \bigr\} \,d\mu(x') \nonumber \\
& = \frac{8L_{\eta}\varepsilon^2}{R^2} \cdot \mu\bigl(B(x,\varepsilon)\bigr) +  p_{\max}\|\eta\|_{\infty} \Bigl[\mu\bigl(B(x,\wt{\varepsilon})\bigr) - \mu\bigl(B(x,\varepsilon)\bigr) \Bigr]. \label{pf:error_euclidean_distance_1}
\end{align}
Equation (1.36) in \cite{trillos2019} states that
\begin{equation*}
\bigl|\mu(B_\mc{X}(x,\varepsilon))  - \omega_m \varepsilon^m\bigr|  \leq CmK_0\varepsilon^{m + 2},
\end{equation*}
where we recall $K_0$ is an upper bound on the sectional curvature of $\mc{X}$. Plugging this back into~\eqref{pf:error_euclidean_distance_1}, we conclude that
\begin{align*}
\int_{\mc{X}} \biggl|\eta\biggl(\frac{\|x' - x\|}{\varepsilon}\biggr) - \eta\biggl(\frac{d_{\mc{X}}(x',x)}{\varepsilon}\biggr)\biggr| \,d\mu(x') & \leq  \frac{8L_{\eta}\varepsilon^2}{R^2}\Bigl[\omega_m\varepsilon^{m} + CmK_0\varepsilon^{m + 2}\Bigr] + \|\eta\|_{\infty} \Bigl[\omega_m(\wt{\varepsilon}^m - {\varepsilon}^m) + 2CmK_0\varepsilon^{m + 2}\Bigr] \\
& \leq \frac{8L_{\eta}\varepsilon^2}{R^2}\Bigl[\omega_m\varepsilon^{m} + R^2CmK_0\varepsilon^{m}\Bigr] + \|\eta\|_{\infty} \varepsilon^{m + 2}\Bigl[\frac{27\omega_m}{R^2} + 2CmK_0\Bigr] \\
& =: C_6 \varepsilon^{m + 2}.
\end{align*}

\paragraph{Putting together the pieces.}
Plugging our upper bounds on $A_1(x)$ and $A_2(x)$ back into~\eqref{pf:error_euclidean_distance_0}, we deduce that
\begin{align*}
\|\wt{L}_{P,\varepsilon}f - L_{P,\varepsilon}f\|_{L^2(\mc{X})}^2 & \leq \frac{1}{\varepsilon^{2(2 + m)}} \int_{\mc{X}} A_1(x) \cdot A_2(x) \,d\mu(x) \\
& \leq \frac{C_6}{\varepsilon^{(2 + m)}} \int_{\mc{X}} A_1(x) \,d\mu(x) \\
& \leq C_5 C_6 \varepsilon^2 |f|_{H^1(\mc{X})}^2,
\end{align*}
thus proving the claimed result.

\subsection{Approximation Error of non-local Laplacian}
\label{subsec:manifold_approximation_error_nonlocal_laplacian}
Fix $x \in \mc{X}$. We begin with a pointwise estimate of $\wt{L}_{P,\varepsilon}f$, facilitated by expressing $w(v) = f(\exp_x(v))$ and $q(v) = p(\exp_x(v))$ in normal coordinates, as in \citep{calder2019}. Let $J_x(\cdot)$ be the Jacobian of the exponential map $\exp_x$, we have
\begin{align*}
\wt{L}_{P,\varepsilon}f(x) & = \frac{1}{\varepsilon^{m + 2}} \int_{\mc{X}} \bigl(f(x') - f(x)\bigr) \eta\biggl(\frac{d_{\mc{X}}(x',x)}{\varepsilon}\biggr) \,dP(x') \\
& = \frac{1}{\varepsilon^{m + 2}} \int_{B(0,\varepsilon) \subset T_x(\mc{X})} \bigl(w(v) - w(0)\bigr) \eta\biggl(\frac{\|v\|}{\varepsilon}\biggr) J_x(v) q(v) \,dv \\
& = \frac{1}{\varepsilon^{2}}\biggl\{\int_{B(0,1)} \bigl(w(\varepsilon v) - w(0)\bigr) \eta(\|v\|) q(\varepsilon v) \,dv + \int_{B(0,1)} \bigl(w(\varepsilon v) - w(0)\bigr) \eta(\|v\|) q(\varepsilon v) \bigl(J_x(\varepsilon v) - 1\bigr) \,dv \biggr\} \\
& = A_1(x) + A_2(x)
\end{align*}
Note that $w$ and $q$ have the same smoothness properties as $f$ and $p$. Moreover, arguing exactly as we did in the flat Euclidean case, we can show that when $f \in H^2(\mc{X})$ and $p \in C^1(\mc{X})$, then
\begin{equation*}
\|A_1\|_{L^2(\mc{X})}^2 \leq C \|f\|_{H^2(\mc{X})}^2
\end{equation*}
whereas if $f \in H^3(\mc{X})$ and $p \in C^2(\mc{X})$ then 
\begin{equation*}
\|A_1 - \sigma_{\eta} \Delta_Pf\|_{L^2(\mc{X})}^2 \leq C\|f\|_{H^3(\mc{X})}^2 \varepsilon^2.
\end{equation*}

Therefore it remains only to upper bound $A_2$ in $L^2(\mc{X})$ norm. To do so, we recall (1.34) of \cite{trillos2019}: for any $\varepsilon < i_0$ and all $x \in \mc{X}$, the Jacobian $J_x(v)$ satisfies the upper bound
\begin{equation*}
|J_x(v) - 1| \leq CmK_0\varepsilon^2, \quad  \textrm{for all} ~~ v \in B(0,\varepsilon) \subseteq T_x(\mc{X}).
\end{equation*}
Combining this estimate with the Cauchy-Schwarz inequality, we conclude that
\begin{align*}
\|A_2\|_{L^2(\mc{X})}^2 & \leq Cm^2K_0^2 \biggl[\int_{B(0,1)} \bigl(w(\varepsilon v) - w(0)\bigr)^2 \eta(\|v\|) q(\varepsilon v) \,dv\biggr] \cdot \biggl[\int_{B(0,1)} \eta(\|v\|) q(\varepsilon v) \,dv\biggr] \\
& \leq Cm^2K_0^2 \sigma_{\eta} (1 + L_q\varepsilon)  \int_{B(0,1)} \bigl(w(\varepsilon v) - w(0)\bigr)^2 \eta(\|v\|) q(\varepsilon v) \,dv \\
& \leq Cm^2K_0^2 \sigma_{\eta}^2 (1 + L_q\varepsilon) p_{\max}  \varepsilon^2 |f|_{H^1(\mc{X})}^2,
\end{align*}
with the final inequality following from (3.2) of~\cite{burago2014}. Combining our estimates on $A_1$ and $A_2$ yields the claim.

\subsection{Estimate of non-local Sobolev seminorm}
\label{subsec:manifold_estimate_nonlocal_seminorm}
In this subsection we establish that the upper bound~\eqref{eqn:manifold_nonlocal_seminorm} holds when $f \in H^s(\mc{X})$ and $p \in C^{s - 1}(\mc{X})$. We first consider $s = 2$, and then $s = 3$.

\paragraph{Case 1: $s = 2$.}
When $s = 2$, the triangle inequality implies that
\begin{equation*}
\dotp{L_{P,\varepsilon}^sf}{f}_P \leq 2p_{\max}\Bigl(\|L_{P,\varepsilon}f - \wt{L}_{P,\varepsilon}\|_{L^2(\mc{X})}^2 + \|\wt{L}_{P,\varepsilon}f\|_{L^2(\mc{X})}^2\Bigr)
\end{equation*}
The first term on the right hand side is upper bounded in~\eqref{eqn:nonlocal_laplacian_geodesic_error}, and the second term is upper bounded in~\eqref{eqn:nonlocal_laplacian_approximation_error_manifold_l2}. Together these estimates imply the claim.

\paragraph{Case 2: $s = 3$.}
When $s = 3$, the triangle inequality implies that
\begin{equation*}
\dotp{L_{P,\varepsilon}^sf}{f}_P = E_{P,\varepsilon}(L_{P,\varepsilon}f;\mc{X}) \leq 3\Bigl( E_{P,\varepsilon}(L_{P,\varepsilon}f - \wt{L}_{P,\varepsilon}f;\mc{X}) +  E_{P,\varepsilon}(\wt{L}_{P,\varepsilon}f - \sigma_{\eta} \Delta_Pf;\mc{X}) + \sigma_{\eta}^2 E_{P,\varepsilon}(\Delta_Pf;\mc{X})\Bigr)
\end{equation*}
We now upper bound each of the three terms on the right hand side of the above inequality. First, we note that by Lemma~\ref{lem:dirichlet_energy_l2} and~\eqref{eqn:nonlocal_laplacian_geodesic_error}, 
\begin{equation*}
 E_{P,\varepsilon}(L_{P,\varepsilon}f - \wt{L}_{P,\varepsilon}f;\mc{X}) \leq  \frac{C}{\varepsilon^2}\|L_{P,\varepsilon}f - \wt{L}_{P,\varepsilon}f\|_{L^2(\mc{X})}^2 \leq C |f|_{H^1(\mc{X})}^2.
\end{equation*}
An equivalent upper bound on $E_{P,\varepsilon}(\wt{L}_{P,\varepsilon}f - \sigma_{\eta} \Delta_Pf;\mc{X})$ follows from Lemma~\ref{lem:dirichlet_energy_l2} and~\eqref{eqn:nonlocal_laplacian_approximation_error_manifold_sobolev}. Finally, we notice that $f \in H^3(\mc{X})$ and $p \in C^2(\mc{X})$ implies $\Delta_Pf \in H^1(\mc{X})$, and furthermore $|\Delta_Pf|_{H^1(\mc{X})} \leq \|p\|_{C^2(\mc{X})} \cdot \|f\|_{H^3(\mc{X})}$. We conclude from Lemma~\ref{lem:dirichlet_energy_sobolev} that
\begin{equation*}
E_{P,\varepsilon}(\Delta_Pf;\mc{X}) \leq C |\Delta_Pf|_{H^1(\mc{X})}^2 \leq C \|f\|_{H^3(\mc{X})}^2,
\end{equation*}
where in the final inequality we have absorbed $\|p\|_{C^2(\mc{X})}$ into the constant $C$. Together, these upper bounds prove the claim.

\subsection{Integrals}
\label{subsec:manifold_integrals}
Recall the Dirichlet energy $E_{P,\varepsilon}(f;\mc{X}) = \dotp{L_{P,\varepsilon}f}{f}_P$, defined in~\eqref{eqn:dirichlet_energy}. Now we establish some estimates on $E_{P,\varepsilon}(f;\mc{X})$ under Model~\ref{def:model_manifold}, and under various assumptions regarding the regularity of $f$.
\begin{lemma}
	\label{lem:dirichlet_energy_l2}
	Suppose Model~\ref{def:model_manifold}, and additionally that $f \in L^2(\mc{X})$. Then there exists a constant $C$ such that
	\begin{equation}
	\label{eqn:dirichlet_energy_l2}
	E_{P,\varepsilon}(f;\mc{X}) \leq \frac{C}{\varepsilon^2} \|f\|_{L^2(\mc{X})}^2.
	\end{equation}
\end{lemma}
\begin{lemma}
	\label{lem:dirichlet_energy_sobolev}
	Suppose Model~\ref{def:model_manifold}, and additionally that $f \in H^1(\mc{X})$. Then there exist constants $c$ and $C$ which do not depend on $f$ such that for any $0 < \varepsilon < c$,
	\begin{equation}
	\label{eqn:dirichlet_energy_sobolev}
	E_{P,\varepsilon}(f;\mc{X}) \leq C |f|_{H^1(\mc{X})}^2.
	\end{equation}
\end{lemma}

We use Lemma~\ref{lem:dirichlet_energy_remainder} to help upper bound the error incurred by using $\|\cdot\|$ rather than $d_{\mc{X}}(\cdot,\cdot)$. Recall the notation $\wt{\varepsilon} = (1 + 27\varepsilon^2/R^2)\varepsilon$, where $R$ is the reach of $\mc{X}$.
\begin{lemma}
	\label{lem:dirichlet_energy_remainder}
	Suppose Model~\ref{def:model_manifold}, and additionally that $f \in H^1(\mc{X})$. There exist constants $c$ and $C$ such that for any $\varepsilon < c$,
	\begin{equation}
	\label{eqn:dirichlet_energy_remainder}
	\int_{\mc{X}} \int_{\mc{X}} \bigl(f(x') - f(x)\bigr)^2 \1\{\varepsilon < d_{\mc{X}}(x',x) \leq \wt{\varepsilon}\} \,d\mu(x') \,d\mu(x) \leq C \varepsilon^{4 + m} \|f\|_{H^1(\mc{X})}^2
	\end{equation}
\end{lemma}

Finally, we use Lemma~\ref{lem:manifold_graph_seminorm_bias2} to show that the pure bias component of $\dotp{L_n^sf,f}_n$ is small in expectation. This is analogous to Lemma~\ref{lem:graph_seminorm_bias2}, except assuming Model~\ref{def:model_manifold} rather than Model~\ref{def:model_flat_euclidean}.
\begin{lemma}
	\label{lem:manifold_graph_seminorm_bias2}
	Assume Model~\ref{def:model_manifold}. Suppose $f \in H^1(\mc{X})$, and let $D_if$ be defined with respect to a kernel $\eta$ that satisfies~\ref{asmp:kernel_manifold}. Then there exists a constant $C$ which does not depend on $f$ or $n$, such that for any $i \in [n]$ and $\bj \in [n]^s$,
	\begin{equation*}
	\Ebb\Bigl[|D_{\bj}f(X_i)| \cdot |f(X_i) - f(X_{\bj_1})|\Bigr] \leq C \varepsilon^{2 + mk} \cdot \|f\|_{H^1(\mc{X})}^2,
	\end{equation*}
	where $k + 1$ is the number of distinct indices in $i\bj$. 
\end{lemma}

\paragraph{Proof (of Lemmas~\ref{lem:dirichlet_energy_l2} and~\ref{lem:dirichlet_energy_sobolev}).}
Define the non-local energy $\wt{E}_{P,\varepsilon}$ with respect to geodesic distance,
\begin{equation*}
\wt{E}_{P,\varepsilon}(f;{\mc{X}}) := \dotp{\wt{L}_{P,\varepsilon}f}{f}_{P} = \int_{\mc{X}} \int_{\mc{X}} \bigl(f(x') - f(x)\bigr)^2 \eta\biggl(\frac{d_{\mc{X}}(x',x)}{\varepsilon}\biggr) \,dP(x') \,dP(x).
\end{equation*}
From the lower bound in~\eqref{eqn:distance_error}, it follows that $E_{P,\varepsilon}(f;X) \leq \wt{E}_{P,\varepsilon}(f;{\mc{X}})$, and from the upper bounds $p(x) \leq p_{\max}$ and $\eta(|x|) \leq \|\eta\|_{\infty} \cdot \1\{x \in [-1,1]\}$ we further have
\begin{equation*}
\wt{E}_{P,\varepsilon}(f;{\mc{X}}) \leq p_{\max}^2 \|\eta\|_{\infty} \cdot \int_{\mc{X}} \int_{B_{\mc{X}}(\varepsilon)} \bigl(f(x') - f(x)\bigr)^2 \,d\mu(x') \,d\mu(x).
\end{equation*}
The estimates~\eqref{eqn:dirichlet_energy_l2} and~\eqref{eqn:dirichlet_energy_sobolev} then respectively follow from (3.1) and Lemma~3.3 of \cite{burago2014}.

\paragraph{Proof (of Lemma~\ref{lem:dirichlet_energy_remainder}).}
Following exactly the steps of the proof of Lemma~3.3 of \citet{burago2014}, but replacing all references to a ball of radius $r$ by references to the set difference between balls of radius $\wt{\varepsilon}$ and $\varepsilon$, we obtain that
\begin{equation*}
\int_{\mc{X}} \int_{\mc{X}} \bigl(f(x') - f(x)\bigr)^2 \1\{\varepsilon < d_{\mc{X}}(x',x) \leq \wt{\varepsilon}\} \,d\mu(x') \,d\mu(x) \leq (1 + CmK_0\varepsilon^2) \cdot \int_{\mc{X}} \int_{B_{m}(0,\wt{\varepsilon})} |d_x^{1}f(v)|^2 \,dv \,d\mu(x).
\end{equation*}
From (2.7) of~\citet{burago2014}, we further have
\begin{equation*}
\int_{\mc{X}} \int_{B_{m}(0,\wt{\varepsilon})} |d_x^{1}f(v)|^2 \,dv \,d\mu(x)  = \frac{\nu_m}{2 + m} (\wt{\varepsilon}^{2 + m} - \varepsilon^{2 + m}) \int_{\mc{X}} |d_x^1f|^2 \,d\mu(x) = 27\frac{\nu_m}{(2 + m)R^2} \varepsilon^{4 + m} \|d^1f\|_{L^2(\mc{X})}^2. 
\end{equation*}
Recalling that $\|d^1f\|_{L^2(\mc{X})}^2 \leq \|f\|_{H^1(\mc{X})}^2$, we see that this implies the claim of Lemma~\ref{lem:dirichlet_energy_remainder}.

\paragraph{Proof (of Lemma~\ref{lem:manifold_graph_seminorm_bias2}).}
The proof of Lemma~\ref{lem:manifold_graph_seminorm_bias2} is identical to the proof of Lemma~\ref{lem:graph_seminorm_bias2}, upon substituting the ambient dimension $m$ for the intrinsic dimension $d$, and using Lemma~\ref{lem:dirichlet_energy_sobolev} rather than Lemma~\ref{lem:dirichlet_estimate_nonlocal_laplacian} to establish~\eqref{pf:graph_seminorm_bias2_1}.

\section{Lower bound on empirical norm}
\label{sec:empirical_norm}
In this Section we prove Proposition~\ref{prop:empirical_norm_sobolev} (in Section~\ref{subsec:empirical_norm_sobolev}). We also prove an analogous result when $\mc{X}$ is a manifold as in Model~\ref{def:model_manifold} (in Section~\ref{subsec:empirical_norm_sobolev_manifold}).

\subsection{Proof of Proposition~\ref{prop:empirical_norm_sobolev}}
\label{subsec:empirical_norm_sobolev}
In this section we establish Proposition~\ref{prop:empirical_norm_sobolev}. As mentioned, the proof of this Proposition follows from the Gagliardo-Nirenberg interpolation inequality, and a one-sided Bernstein's inequality (Lemma~\ref{lem:one_sided_bernstein}). 

\begin{lemma}[Gagliardo-Nirenberg interpolation inequality]
	\label{lem:gagliardo_nirenberg}
	Suppose Model~\ref{def:model_flat_euclidean}, and that $f \in H^s(\mc{X})$ for some $s \geq d/4$. Then there exist constants $C_1$ and $C_2$ that do not depend on $f$, such that
	\begin{equation}
	\label{eqn:gagliardo_nirenberg}
	\|f\|_{L^4(\mc{X})} \leq C_1 |f|_{H^s(\mc{X})}^{d/4s} \|f\|_{L^2(\mc{X})}^{1 - d/(4s)} + C_2 \|f\|_{L^2(\mc{X})}
	\end{equation}
\end{lemma}

\paragraph{Proof (of Proposition~\ref{prop:empirical_norm_sobolev}).}
Rearranging~\eqref{eqn:gagliardo_nirenberg} and raising both sides to the $4$th power, we see that
\begin{equation*}
\frac{\Ebb[f^4(X)]}{\|f\|_P^4} \leq C \biggl(\frac{\|f\|_{L^4(\mc{X})}}{\|f\|_{L^2(\mc{X})}}\biggr)^4 \leq C_1\biggl(\frac{|f|_{H^s(\mc{X})}}{\|f\|_{L^2(\mc{X})}}\biggr)^{d/s} + C_2,
\end{equation*}
here the constants $C_1,C_2$ are not the same as in~\eqref{eqn:gagliardo_nirenberg}. Therefore taking the constant $C$ in assumption~\eqref{eqn:empirical_norm_sobolev_1} to be sufficiently large relative to $C_1$ and $C_2$, we have that
\begin{equation*}
C_1\biggl(\frac{|f|_{H^s(\mc{X})}}{\|f\|_{L^2(\mc{X})}}\biggr)^{d/s} \leq \frac{\delta n}{64},
\end{equation*} 
and consequently 
\begin{equation*}
\frac{\Ebb[f^4(X)]}{\|f\|_P^4} \leq \frac{\delta n}{8} + 8C_2^3.
\end{equation*}
The claim then follows from Lemma~\ref{lem:one_sided_bernstein}, upon taking $c = 1/(64C_2^3)$ in the statement of Proposition~\ref{prop:empirical_norm_sobolev}.

\subsection{Proof of Proposition~\ref{prop:empirical_norm_sobolev_manifold}}
\label{subsec:empirical_norm_sobolev_manifold}
The proof of Proposition~\ref{prop:empirical_norm_sobolev_manifold} follows exactly the same steps as the proof of Proposition~\ref{prop:empirical_norm_sobolev}, upon replacing Lemma~\ref{lem:gagliardo_nirenberg} by Lemma~\ref{lem:gagliardo_nirenberg_manifold}.
\begin{lemma}[(c.f Theorem~3.70 of~\citet{aubin2012})]
	\label{lem:gagliardo_nirenberg_manifold}
	Suppose Model~\ref{def:model_manifold}, and that $f \in H^s(\mc{X})$ for some $s \geq m/4$. Then there exist constants $C_1$ and $C_2$ that do not depend on $f$, such that
	\begin{equation}
	\label{eqn:gagliardo_nirenberg_manifold}
	\|f\|_{L^4(\mc{X})} \leq C_1 |f|_{H^s(\mc{X})}^{m/4s} \|f\|_{L^2(\mc{X})}^{1 - m/(4s)} + C_2 \|f\|_{L^2(\mc{X})}.
	\end{equation}
\end{lemma}
	
\section{Proof of Main Results}
\label{sec:proofs_main_results}

\subsection{Estimation Results}

\paragraph{Proof of Theorem~\ref{thm:laplacian_eigenmaps_estimation_fo}.}
We condition on the event that the design points $X_1,\ldots,X_n$ satisfy
\begin{equation}
\label{pf:laplacian_eigenmaps_estimation_fo_1}
\dotp{L_{n,\varepsilon}f_0}{f_0}_n \leq \frac{C}{\delta}M^2 \quad \textrm{and} \quad \lambda_k \geq \min\{\lambda_k(\Delta_P), \varepsilon^{-2}\}~~\textrm{for all $2 \leq k \leq n$.}
\end{equation}
Note that by Propositions~\ref{prop:graph_seminorm_fo} and~\ref{prop:graph_eigenvalue}, these statements are both satisfied with probability at least $1 - \delta - Cn\exp\{-cn\varepsilon^d\}$. 

Conditional on~\eqref{pf:laplacian_eigenmaps_estimation_fo_1}, we have from Lemma~\ref{lem:fixed_graph_estimation} that for any $0 \leq K \leq n$,
\begin{equation*}
\|\wh{f} - f_0\|_n^2 \leq C\biggl\{\frac{M^2}{\delta (\lambda_{K + 1}(\Delta_P) \wedge \varepsilon^{-2})} + \frac{K}{n}\biggr\},
\end{equation*}
either deterministically (when $K = 0$), or with probability at least $1 - \exp(-K)$ (when $K \geq 1$). Further, from the bounds $\varepsilon \leq c_0 K^{-1/d}$ (Assumption~\ref{asmp:parameters_estimation_fo}) and $\lambda_{K + 1}(\Delta_P) \geq c (K + 1)^{2/d}$ (Weyl's Law) we can simply the above expression to the following,
\begin{equation}
\label{pf:laplacian_eigenmaps_estimation_fo_2}
\|\wh{f} - f_0\|_n^2 \leq C\biggl\{\frac{M^2}{\delta}(K + 1)^{-2/d} + \frac{K}{n}\biggr\}.
\end{equation}
We now upper bound the right hand side of~\eqref{pf:laplacian_eigenmaps_estimation_fo_2}, based on the value of $K$ chosen in~\ref{asmp:parameters_estimation_fo}.  When possible we choose $K = \floor{M^2n}^{d/(2 + d)}$ to balance bias and variance, in which case~\eqref{pf:laplacian_eigenmaps_estimation_fo_2} implies
\begin{equation*}
\|\wh{f} - f_0\|_n^2 \leq \frac{C}{\delta} M^2 (M^2n)^{-2/(2 + d)}.
\end{equation*}
If $M^2 < n^{-1}$, then we take $K = 1$, and from~\eqref{pf:laplacian_eigenmaps_estimation_fo_2} we get
\begin{equation*}
\|\wh{f} - f_0\|_n^2 \leq \frac{C}{n\delta}.
\end{equation*}
Finally if $M > n^{1/d}$, we take $K = n$. In this case, we note that $\wh{f}(X_i) = Y_i$ for all $i = 1,\ldots,n$, and it immediately follows that
\begin{equation*}
\|\wh{f} - f_0\|_n^2 = \frac{1}{n}\sum_{i = 1}^{n} w_i^2 \leq 5,
\end{equation*}
with probability at least $1 - \exp(-n)$. Combining these three separate cases yields the conclusion of Theorem~\ref{thm:laplacian_eigenmaps_estimation_fo}.

\paragraph{Proof of Theorem~\ref{thm:laplacian_eigenmaps_estimation_ho}.}
Follows identically to the proof of Theorem~\ref{thm:laplacian_eigenmaps_estimation_fo}, except substituting $L_{n,\varepsilon}^s$ for $L_{n,\varepsilon}$, $\lambda_k^s$ for $\lambda_k$, and using Proposition~\ref{prop:graph_seminorm_ho} rather than Proposition~\ref{prop:graph_seminorm_fo} and Assumption~\ref{asmp:parameters_estimation_ho} rather than Assumption~\ref{asmp:parameters_estimation_fo}.

\paragraph{Proof of Theorem~\ref{thm:laplacian_eigenmaps_estimation_manifold}.}
Follows identically to the proof of Theorem~\ref{thm:laplacian_eigenmaps_estimation_fo}, substituting $L_{n,\varepsilon}^s$ for $L_{n,\varepsilon}$, $\lambda_k^s$ for $\lambda_k$, and using Proposition~\ref{prop:graph_seminorm_manifold} rather than Proposition~\ref{prop:graph_seminorm_fo}, Proposition~\ref{prop:graph_eigenvalue_manifold} rather than Proposition~\ref{prop:graph_eigenvalue}, and Assumption~\ref{asmp:parameters_estimation_manifold} rather than Assumption~\ref{asmp:parameters_testing_fo}.

\subsection{Testing Results}

\paragraph{Proof of Theorem~\ref{thm:laplacian_eigenmaps_testing_fo}.}
We have already upper bounded the Type I error of $\varphi$ in Lemma~\ref{lem:fixed_graph_testing}, and it remains to upper bound the Type II error. To do so, we condition on the event that the design points $X_1,\ldots,X_n$ satisfy,
\begin{equation}
\label{pf:laplacian_eigenmaps_testing_fo_1}
\dotp{L_{n,\varepsilon}f_0}{f_0}_n \leq \frac{C}{\delta}M^2,\quad \textrm{and} \quad \lambda_k \geq \min\{\lambda_k(\Delta_P), \varepsilon^{-2}\}~~\textrm{for all $2 \leq k \leq n$,}
\end{equation}
as well as that
\begin{equation}
\label{pf:laplacian_eigenmaps_testing_fo_2}
\|f_0\|_n^2 \geq \frac{1}{2}\|f_0\|_P^2.
\end{equation}
Note that by Propositions~\ref{prop:graph_seminorm_fo} and~\ref{prop:graph_eigenvalue}, both statements in~\eqref{pf:laplacian_eigenmaps_testing_fo_1} are satisfied with probability at least $1 - \delta - Cn\exp\{-cn\varepsilon^d\}$. Additionally, by Proposition~\ref{prop:empirical_norm_sobolev} and the assumption in~\eqref{eqn:laplacian_eigenmaps_testing_criticalradius_fo} that $\|f_0\|_P^2 \geq CM^2/(bn^{2/d})$, the one-sided inequality~\eqref{pf:laplacian_eigenmaps_testing_fo_2} follows with probability at least $1 - \exp\{-(cn \wedge 1/b)\}$. Setting $\delta = b/3$ and taking $n \geq N$ to be sufficiently large, the bottom line is that both~\eqref{pf:laplacian_eigenmaps_testing_fo_1} and~\eqref{pf:laplacian_eigenmaps_testing_fo_2} are together satisfied with probability at least $1 - b/2$.

Now, to complete the proof of Theorem~\ref{thm:laplacian_eigenmaps_testing_fo}, we would like to invoke Lemma~\ref{lem:fixed_graph_testing}, and conclude that conditional on $X_1,\ldots,X_n$ satisfying~\eqref{pf:laplacian_eigenmaps_testing_fo_1} and~\eqref{pf:laplacian_eigenmaps_testing_fo_2}, our test $\varphi$ will equal $1$ with probability at least $1 - b/2$. To use Lemma~\ref{lem:fixed_graph_testing}, we will need to establish that~\eqref{eqn:fixed_graph_testing_critical_radius} is satisfied, which we now show. 

On the one hand, we have that the right hand side of~\eqref{eqn:fixed_graph_testing_critical_radius} is upper bounded, 
\begin{align*}
\frac{\dotp{L_{n,\varepsilon}f_0}{f_0}_n}{\lambda_{K + 1}} + \frac{\sqrt{2K}}{n}\biggl[2\sqrt{\frac{1}{a}} + \sqrt{\frac{2}{b}} + \frac{32}{bn}\biggr] & \leq C\biggl(\frac{M^2}{b \min\{\lambda_{K+1}(\Delta_P), \varepsilon^{-2}\}} + \frac{\sqrt{2K}}{n}\biggl[\sqrt{\frac{1}{a}} + \frac{1}{b}\biggr]\biggr) \\
& \leq C\biggl(\frac{M^2}{b}K^{-2/d} + \frac{\sqrt{2K}}{n}\biggl[\sqrt{\frac{1}{a}} + \frac{1}{b}\biggr]\biggr)
\end{align*}
with the second inequality following by the assumption $\varepsilon \leq K^{-1/d}$ and Weyl's Law. On the other hand, we have that $\|f_0\|_n^2 \geq \|f_0\|_P^2/2$. Consequently, to prove Theorem~\ref{thm:laplacian_eigenmaps_testing_fo}, it remains only to verify that
\begin{equation}
\label{pf:laplacian_eigenmaps_testing_fo_3}
\|f_0\|_P^2 \geq C\biggl(\frac{M^2}{b}K^{-2/d} + \frac{\sqrt{2K}}{n}\biggl[\sqrt{\frac{1}{a}} + \frac{1}{b}\biggr]\biggr).
\end{equation}
As in the estimation case, we can further upper bound the right hand side of~\eqref{pf:laplacian_eigenmaps_testing_fo_3}, depending on the value of $K$ chosen in~\ref{asmp:parameters_testing_fo}. The classical case is $K = (M^2n)^{d/(2 + d)}$, in which case~\eqref{pf:laplacian_eigenmaps_testing_fo_3} is satisfied as long as
\begin{equation*}
\|f_0\|_P^2 \geq CM^2(M^2n)^{-4/(4 + d)}\biggl[\sqrt{\frac{1}{a}} + \frac{1}{b}\biggr]
\end{equation*}
If $M^2 < n^{-1}$, then we take $K = 1$, and~\eqref{pf:laplacian_eigenmaps_testing_fo_3} is satisfied whenever
\begin{equation*}
\|f_0\|_P^2 \geq \frac{C}{n}\biggl[\sqrt{\frac{1}{a}} + \frac{1}{b}\biggr].
\end{equation*}
Finally if $M > n^{1/d}$, we take $K = n$, and~\eqref{pf:laplacian_eigenmaps_testing_fo_3} is satisfied if
\begin{equation*}
\|f_0\|_P^2 \geq C\biggl(\frac{M^2}{n^{2/d}b} + n^{-1/2}\biggl[\sqrt{\frac{1}{a}} + \frac{1}{b}\biggr]\biggr).
\end{equation*}
We conclude by observing that~\eqref{eqn:laplacian_eigenmaps_testing_criticalradius_fo} implies each of these three inequalities, and thus implies~\eqref{pf:laplacian_eigenmaps_testing_fo_3}.

\paragraph{Proof of Theorem~\ref{thm:laplacian_eigenmaps_testing_ho}.}
Follows identically to the proof of Theorem~\ref{thm:laplacian_eigenmaps_estimation_fo}, except substituting $L_{n,\varepsilon}^s$ for $L_{n,\varepsilon}$, $\lambda_k^s$ for $\lambda_k$, and using Proposition~\ref{prop:graph_seminorm_ho} rather than Proposition~\ref{prop:graph_seminorm_fo} and Assumption~\ref{asmp:parameters_testing_ho} rather than Assumption~\ref{asmp:parameters_testing_fo}.

\paragraph{Proof of Theorem~\ref{thm:laplacian_eigenmaps_testing_manifold}.}
Follows identically to the proof of Theorem~\ref{thm:laplacian_eigenmaps_estimation_fo}, except substituting $L_{n,\varepsilon}^s$ for $L_{n,\varepsilon}$, $\lambda_k^s$ for $\lambda_k$, and using Proposition~\ref{prop:graph_seminorm_manifold} rather than Proposition~\ref{prop:graph_seminorm_fo}, Proposition~\ref{prop:graph_eigenvalue_manifold} rather than Proposition~\ref{prop:graph_eigenvalue}, Proposition~\ref{prop:empirical_norm_sobolev_manifold} rather than Proposition~\ref{prop:empirical_norm_sobolev}, and Assumption~\ref{asmp:parameters_testing_manifold} rather than Assumption~\ref{asmp:parameters_testing_fo}.

\paragraph{Proof of Theorem~\ref{thm:laplacian_eigenmaps_testing_ho_suboptimal}.}
Note that our choices of $K$ and $\varepsilon$ ensure  that~\eqref{pf:laplacian_eigenmaps_testing_fo_1} (with $L_{n,\varepsilon}^s$ replacing $L_{n,\varepsilon}$) and~\eqref{pf:laplacian_eigenmaps_testing_fo_2} are satisfied with probability at least $1 - b/2$. Proceeding as in the proof of Theorem~\ref{thm:laplacian_eigenmaps_testing_fo}, we upper bound the right hand side of~\eqref{eqn:fixed_graph_testing_critical_radius},
\begin{align*}
\frac{\dotp{L_{n,\varepsilon}f_0}{f_0}_n}{\lambda_{K + 1}} + \frac{\sqrt{2K}}{n}\biggl[2\sqrt{\frac{1}{a}} + \sqrt{\frac{2}{b}} + \frac{32}{bn}\biggr] & \leq C\biggl(\frac{M^2}{b \min\{\lambda_{K+1}(\Delta_P), \varepsilon^{-2}\}} + \frac{\sqrt{2K}}{n}\biggl[\sqrt{\frac{1}{a}} + \frac{1}{b}\biggr]\biggr) \\
& \leq C\biggl(\frac{M^2}{b}\varepsilon^2 + \frac{\sqrt{2K}}{n}\biggl[\sqrt{\frac{1}{a}} + \frac{1}{b}\biggr]\biggr).
\end{align*}
Unlike in the proof of Theorem~\ref{thm:laplacian_eigenmaps_testing_fo}, we note that in this case $\varepsilon^2 \leq C\lambda_K(\Delta_P)$ rather than vice versa. From here, proceeding as in the proof of Theorem~\ref{thm:laplacian_eigenmaps_testing_fo} gives the claimed result.

\section{Graph Laplacian methods and the cluster assumption}
\label{subsec:eigenmaps_beats_kernel_smoothing}

A main conclusion of our paper is that PCR-LE is minimax optimal for nonparametric regression over certain Sobolev classes. It is not the only optimal method. For instance, kernel smoothing and least squares using an appropriate set of basis functions as features are two other minimax optimal methods over these Sobolev classes. We now give an example where PCR-LE is better than these two alternatives, in the sense of having (much) smaller risk. This is possible because PCR-LE performs remarkably well when the regression function $f_0$ and design distribution $P$ satisfy a \emph{cluster assumption}: that is, when the regression function is (approximately) piecewise constant over high-density clusters of the design distribution $P$. On the other hand, kernel smoothing (with Euclidean distance) and least squares (using eigenfunctions of an unweighted Laplace operator) cannot take advantage of the cluster assumption. We call this property of PCR-LE \emph{density adaptivity}.

\subsection{Setup}
We begin by specifying a sequence of design densities and regression functions $\{(p^{(n)}, f_0^{(n)}): n \in \mathbb{N}\}$. These distributions will all be chosen to satisfy the cluster assumption. To that end, we define two clusters $Q_1,Q_2 \subset \Reals$ using a cluster separation parameter $r$, as
\begin{equation*}
Q_1 := [0,1/2 - r], \quad Q_2 := [1/2 + r,1],
\end{equation*}
and take the domain $\mc{X}^{(n)} := Q_1 \cup Q_2$. We then take the design density to be uniform over $\mc{X}^{(n)}$ and the regression function to be a piecewise constant function over $Q_1$ and $Q_2$ of height $\theta$,
\begin{equation}
\label{def:model_cluster_assumption}
p^{(n)}(x) := \frac{1}{1 - 2r}\1\bigl\{x \in Q_1 \cup Q_2\bigr\}, \quad f_0^{(n)}(x) := \theta \cdot \Bigl(\1\bigl\{x \in Q_1\bigr\} - \1\bigl\{x \in Q_2\bigr\}\Bigr).
\end{equation}
Thus $p^{(n)}$ and $f_0^{(n)}$ belong to a two-parameter family, where the parameters are the cluster separation $r$ and height $\theta$. Generally speaking, the smaller the separation $r$, and the larger the height $\theta$, the more graph Laplacian methods will outperform both kernel smoothing and linear regression using eigenfunctions of the unweighted Laplace operator as features.

We now define kernel smoothing and least squares using eigenfunction of an unweighted Laplace operator For a kernel function $\psi$ and bandwidth parameter $h$, the kernel smoothing estimator $\wt{f}_{\mathrm{KS}}$ is defined at a point $x \in \mc{X}$ as 
\begin{equation}
\label{eqn:kernel_smoothing}
\wt{f}_{\mathrm{KS}}(x) := 
\begin{dcases*}
0, & \quad \textrm{if $d_{n,h}(x) = 0$,} \\
\frac{1}{d_{n,h}(x)} \sum_{i = 1}^{n} Y_i \psi\biggl(\frac{\|X_i - x\|}{h}\biggr), & \quad \textrm{otherwise.}
\end{dcases*}
\end{equation}

Let $(\lambda_1,\phi_1), (\lambda_2,\phi_2),\ldots$ be eigenpairs of the unweighted Laplace operator $\Delta$ on $[0,1]$, meaning
\begin{equation}
\label{eqn:laplace_beltrami_eigenpairs}
\Delta \phi_k = \lambda_k \phi_k, ~~ \|\phi_k\|_{L^2([0,1])} = 1,~~\frac{d}{dx}\phi_k(0) = \frac{d}{dx}\phi_k(1) = 0. 
\end{equation}
In this case the eigenfunctions $\phi_k$ of $\Delta$ are simply cosine functions, with eigenvalues proportional to their squared frequency. Noting that $\phi_1(x) = 1$ and $\lambda_1 = 0$, for $k = 2,3,\ldots$ we have
\begin{equation*}
\phi_k(x) = \sqrt{2}\cdot\cos(2\pi k x),~~\lambda_k(\Delta) = \pi^2 k^2.
\end{equation*}
The least squares estimator using $\phi_1,\ldots,\phi_K$ ($1 \leq K \leq n$) eigenfunctions as features is simply\footnote{For convenience, we will assume $\Phi \in \Reals^{n \times K}$ is full rank. If this is not the case, the least squares estimator $\wt{f}_K$ is not uniquely defined, but any solution will equal ${\bf Y}$ in-sample, and will satisfy $\|\wt{f}_K - f_0\|_n^2 \geq 1/2$ with high probability.}
\begin{equation}
\label{eqn:least_squares}
\wt{f}_{K} := \argmin_{f \in \mathrm{span}\{\phi_1,\ldots,\phi_K\}} \|Y - f\|_n^2 = \Phi (\Phi^{\top} \Phi)^{-1} \Phi^{\top} Y.
\end{equation}
Hereafter, we will refer to $\wt{f}_K$ as the \emph{uniform least squares} estimator.

\subsection{Upper bounds on risk of PCR-LE}
Now we are in a position to state our results. Both PCR-LE and kernel smoothing depend in part on the choice of kernel. For simplicity, in our analysis we only consider the boxcar kernel,
\begin{equation}
\label{asmp:boxcar_kernel}
\eta(z) = \psi(z) = \1\{z \leq 1\}.
\end{equation}
This is strictly for convenience, and the following results will also hold for any kernel that satisfies~\ref{asmp:kernel_flat_euclidean}.

\begin{proposition}
	\label{prop:cluster_assumption_ub}
	Suppose $(X_1,Y_1),\ldots(X_n,Y_n)$ are sampled according to~\eqref{def:model_cluster_assumption}. Compute the PCR-LE estimator $\wh{f}$ using a kernel $\eta$ which satisfies~\eqref{asmp:boxcar_kernel}, number of eigenvectors $K = 2$, and radius $\varepsilon = r/2$. Then,
	\begin{equation}
	\label{eqn:cluster_assumption_le_ub}
	\Ebb\Bigl[\|\wh{f} - f_0^{(n)}\|_{n}^2\Bigr] \leq \biggl(6\theta^2 + \frac{1}{n}\biggr) \cdot \frac{8}{r} \exp(-nr/8) + \frac{1}{n}.
	\end{equation}
\end{proposition}
\paragraph{Proof of Proposition~\ref{prop:cluster_assumption_ub}.}

We begin by showing that, with high probability, the eigenvectors $v_1,v_2$ respect the cluster structure of $p^{(n)}$. Denote $u_1 = (\1\{X_i \in Q_1\})_{i \in [n]}$, and likewise $u_2 = (\1\{X_i \in Q_2\})_{i \in [n]}$. We make the following two observations:
\begin{enumerate}
	\item Because $\varepsilon < r$ and the kernel $\eta$ is compactly supported on $[0,1]$, for each $X_i \in Q_1$ and $X_j \in Q_2$, it must be the case that $\eta(\|X_i - X_j\|/\varepsilon) = 0$.
	\item Using an elementary concentration argument (stated in Lemma~\ref{lem:balls_in_bins}) and the triangle inequality, we deduce that with probability at least $1 - 4/\varepsilon\exp(-n\varepsilon/4)$ there exists a path in $G_{n,\varepsilon}$ between each $X_i, X_j \in Q_1$, and likewise between each $X_i,X_j \in Q_2$.
\end{enumerate}
Together these observations imply that with high probability the neighborhood graph $G_{n,\varepsilon}$ consists of exactly two connected components: one consisting of all design points $X_i \in Q_1$, and the other consisting of all design points $X_i \in Q_2$. In other words, 
\begin{equation}
\label{pf:eigenmaps_beats_kernel_smoothing_1}
\Pbb\Bigl(\mathrm{span}\{v_1,v_2\} = \mathrm{span}\{u_1,u_2\}\Bigr) \geq 1 - 4/\varepsilon\exp(-n\varepsilon/4).
\end{equation}

Let us condition on the ``good'' event $\mc{E}$ that the design points $X_1,\ldots,X_n$ satisfy~\eqref{eqn:balls_in_bins}, and therefore that $\mathrm{span}\{v_1,v_2\} = \mathrm{span}\{u_1,u_2\}$. Consider the empirical mean $\wb{Y}_Q := \frac{1}{\sharp\{Q \cup {\bf X}\}} \sum_{i: X_i \in Q} Y_i$. Since $\mathrm{span}\{v_1,v_2\} = \mathrm{span}\{u_1,u_2\}$, the estimator $\wh{f} = \wh{f}_{\mathrm{LE}}$ will be piecewise constant on $Q_1$ and $Q_2$, and in fact we have that
\begin{equation}
\label{pf:eigenmaps_beats_kernel_smoothing_-1}
\wh{f} = \wb{Y}_{Q_1}u_1 + \wb{Y}_{Q_2}u_2.
\end{equation}
Therefore conditional on $\mc{E}$,
\begin{equation*}
\|\wh{f} - f_0^{(n)}\|_{n}^2 = P_n(Q_1) \cdot (\wb{Y}_{Q_1} - \theta)^2 + P_n(Q_2) \cdot (\wb{Y}_{Q_2} + \theta)^2
\end{equation*}
and consequently,
\begin{equation}
\label{pf:eigenmaps_beats_kernel_smoothing_0}
\Ebb\Bigl[\|\wh{f} - f_0^{(n)}\|_{n}^2 \Big|\mc{E}\Bigr] = \Ebb\Bigl[\Ebb\bigl[\|\wh{f} - f_0^{(n)}\|_{n}^2\big|X_1,\ldots,X_n\bigr] \Bigr|\mc{E}\Bigr] = \frac{1}{n}.
\end{equation}

Now we derive a crude upper bound on $\|\wh{f} - f_0^{(n)}\|_{n}$ that will suffice to control the error conditional on $\mc{E}^c$. We observe that the empirical norm of $\wh{f}$ is bounded,
\begin{equation*}
\|\wh{f}\|_n^2 \leq \frac{2}{n}\sum_{i = 1}^{n}\dotp{Y}{v_1}_n^2v_{1,i}^2 + \dotp{{\bf Y}}{v_2}_n^2v_{2,i}^2 \leq 2\bigl(\dotp{{\bf Y}}{v_1}_n^2 + \dotp{{\bf Y}}{v_2}_n^2\bigr) \leq 4 \|{\bf Y}\|_n^2. 
\end{equation*}
Noting that $\Ebb[\|{\bf Y}\|_n^2|X_1,\ldots,X_n] = \|f_0\|_n^2 + 1/n  = \theta^2 + 1/n$, we conclude that
\begin{equation*}
\Ebb\Bigl[\|\wh{f} - f_0\|_{n}^2 \cdot \1\{\mc{E}^c\}\Bigr] \leq \Ebb\Bigl[\Bigl(2 \|f_0\|_{n}^2 + 4(\theta^2 + 1/n) \cdot \1\{\mc{E}^c\} \Bigr] \leq (6\theta^2 + n^{-1}) \cdot 4\varepsilon^{-1} \exp(-n\varepsilon/4).
\end{equation*}
Combining this with~\eqref{pf:eigenmaps_beats_kernel_smoothing_0} implies~\eqref{eqn:cluster_assumption_le_ub}.

\subsection{Lower bounds on risk of kernel smoothing and least squares}

\begin{proposition}
	\label{prop:cluster_assumption_lb}
	Suppose $(X_1,Y_1),\ldots,(X_n,Y_n)$ are sampled according to~\eqref{def:model_cluster_assumption}. Suppose $(\log n)^2/n \leq r \leq c$, where $c$ is a universal constant. 
	\begin{itemize}
		\item Compute the kernel smoothing estimator $\wt{f} = \wt{f}_{\mathrm{KS}}$ as in~\eqref{eqn:kernel_smoothing}, using a kernel $\psi$ which satisfies~\eqref{asmp:boxcar_kernel}. Then there exist universal constants $c, N > 0$ such that for all $n > N$,
		\begin{equation}
		\label{eqn:cluster_assumption_ks_lb}
		\inf_{h' > 0} \Ebb\Bigl[\|\wt{f} - f_0^{(n)}\|_{n}^2\Bigr] \geq c\min\biggl\{ \frac{r^{-1}}{n}, \frac{\theta}{\sqrt{n}} \biggr\}.
		\end{equation}
		\item Compute the least squares estimator $\wt{f} = \wt{f}_{\mathrm{SP}}$ as in~\eqref{eqn:least_squares}. Then there exist universal constants $c, N > 0$ such that for all $n > N$,
		\begin{equation}
		\label{eqn:cluster_assumption_sp_lb}
		\inf_{1 \leq K \leq n} \Ebb \Bigl[\|\wt{f} - f_0^{(n)}\|_n^2\Bigr] \geq c\min\biggl\{\frac{r^{-1}}{n}, \frac{1}{\log(n)}, \frac{r^{-2/3}}{n}, \frac{\sqrt{\theta}}{n^{3/4}}\biggr\}.
		\end{equation}
	\end{itemize}
\end{proposition}
The proof of Proposition~\ref{prop:cluster_assumption_lb} is long, and we defer it until after some discussion of the implications of the proposition.

Together, Propositions~\ref{prop:cluster_assumption_ub} and~\ref{prop:cluster_assumption_lb} illustrate that the risk of PCR-LE can be dramatically smaller than that of kernel smoothing or uniform least squares. For instance, taking $\theta = n^{1/2}$ and $r = n^{-3/4}$, when appropriately tuned, $\wh{f}$ satisfies
\begin{equation*}
\Ebb\Bigl[\|\wh{f} - f_0^{(n)}\|_{n}^2\Bigr] \leq C\biggl(n^{7/4}\exp(-n^{1/4}/8)) + \frac{1}{n}\biggr) \leq \frac{C}{n},
\end{equation*}
for a universal constant $C$ and all $n$ larger than some universal constant $N$, whereas for $\wt{f} = \wt{f}_{\mathrm{KS}}$,
\begin{equation*}
\inf_{h' > 0} \Ebb\Bigl[\|\wt{f} - f_0^{(n)}\|_{n}^2\Bigr] \geq \frac{c}{n^{1/4}},
\end{equation*}
and for $\wt{f} = \wt{f}_{\mathrm{SP}}$,
\begin{equation*}
\inf_{1 \leq K \leq n} \Ebb\Bigl[\|\wt{f} - f_0^{(n)}\|_{n}^2\Bigr] \geq \frac{c}{n^{1/2}}.
\end{equation*}
Other choices of $\theta$ and $r$ lead to even more dramatic gaps between the risk of PCR-LE, and the risk of kernel smoothing and least squares. The overall takeaway is that under Model~\ref{def:model_cluster_assumption}, estimators that use the graph Laplacian can converge to the true regression function $f_0^{(n)}$ at fast rates---parametric rates that do not depend on the $L^2$ norm of $f_0^{(n)}$---whereas other estimators, optimal for estimation over Sobolev spaces, converge to $f_0^{(n)}$ at slow rates---nonparametric rates that deteriorate as the $L^2$ norm of $f_0^{(n)}$ grows.

Some remarks:
\begin{itemize}
	\item The lower bound on the in-sample risk of $\wt{f}_{\mathrm{KS}}$  given by \eqref{eqn:cluster_assumption_ks_lb} is larger than that of $\wt{f}_{\mathrm{SP}}$ given by \eqref{eqn:cluster_assumption_sp_lb}. This does not mean that kernel smoothing exhibits less adaptivity to the cluster assumption than uniform least squares. Instead, we suspect it is due to looseness in our lower bounds: we are able to tightly control the bias of kernel smoothing, whereas we must use a potentially loose bound on the bias of uniform least squares. Experimentally, it appears that kernel smoothing usually outperforms uniform least squares, under various instantiations of the cluster assumption.
	\item The cluster assumption---in which the regression function is piecewise constant and $p$ consists of multiple connected components---is a very strong assumption. The \emph{low-density separation} condition is a related but weaker assumption, in which the regression function is assumed to be smoother (but not constant) in regions of higher density. This is a rather general hypothesis which can formalized in a number of different ways. For instance, one could insist that the regression function $f_0$ belong to a normed ball in a \emph{weighted Sobolev space}, with semi-norm given by
	\begin{equation*}
	|f_0|_{H^s(P)} := \dotp{\Delta_P^s f_0}{f_0}_P.
	\end{equation*}
	Intuitively, when $\|f_0\|_{H^s(P)}$ is much smaller than $\|f_0\|_{H^s(\mc{X})}$, density-adaptive learners such as PCR-LE should have the advantage on non-density adaptive linear smoothers, such as kernel smoothing or uniform least squares. Indeed, in the case of Model~\ref{def:model_cluster_assumption} we see that
	\begin{equation*}
	\|f_0^{(n)}\|_{H^s(P^{(n)})} = 0~~\textrm{for all $s \in \mathbb{N}$, and all $r, \theta > 0$,}
	\end{equation*}
	whereas $f_0^{(n)}$ does not even belong to the first-order Sobolev space $H^1([0,1])$. In words, this shows the cluster assumption is an extreme case of the low-density separation condition. Unfortunately, it is quite difficult to analyze graph-based estimators under the general low-density separation condition, without making strong assumptions on $P$. 
	\item Finally, we note that either changing the graph or the normalization of the Laplacian fundamentally alters the type of density adaptivity displayed by graph-Laplacian-based estimators; see~\citet{hoffmann2019} for an extensive discussion.
\end{itemize}

\subsection{Proof of Proposition~\ref{prop:cluster_assumption_lb}}
First we show~\eqref{eqn:cluster_assumption_ks_lb}, then~\eqref{eqn:cluster_assumption_sp_lb}.

\subsubsection{Proof of~\eqref{eqn:cluster_assumption_ks_lb}}
A standard argument using the law of iterated expectation implies the following lower bound on the pointwise risk in terms of squared-bias and variance-like quantities,
\begin{equation*}
\Ebb\Bigl[\Bigl(\wt{f}(X_i) - f_0(X_i)\Bigr)^2|X_i = x\Bigr] \geq \frac{(n - 1)}{n}\Ebb\biggl[\Bigl(f_0(X) - f_0(x)\Bigr)^2|X \in B(x,h')\biggr] + \Ebb\biggl[\frac{1}{d_{n,h'}(x)}\biggr].
\end{equation*}
The variance term can be lower bounded quite simply for any $x \in \mc{X}^{(n)}$; noting that $\sup_{x} p^{(n)}(x) < 2$ and $\nu(B(x,h') \cap \mc{X}^{(n)}) \leq 2h'$, it follows by Jensen's inequality that 
\begin{equation*}
\Ebb\biggl[\frac{1}{d_{n,h'}(x)}\biggr] \geq \frac{1}{\Ebb[d_{n,h'}(x)]} \geq \frac{1}{4nh'}.
\end{equation*}
On the other hand the squared bias term is quite large for $x$ close to $1/2$. Precisely, if $h' \geq 4r$ then a simple calculation implies
\begin{equation*}
\Ebb[(f_0(X) - f_0(x))^2|X \in B(x,h')] \geq \frac{\theta^2}{8} \quad \textrm{for all}~x \in [(1 - h'/2)_{+}, 1/2 - r].
\end{equation*}
Combining these lower bounds on variance and squared bias terms and summing over $X_1,\ldots,X_n$, we arrive at the following: if $h' \leq 4r$, then
\begin{equation*}
\Ebb\Bigl[\|\wt{f} - f_0^{(n)}\|_{n}^2\Bigr] = \frac{1}{n}\sum_{i = 1}^{n} \Ebb\biggl[\Ebb\Bigl[\Bigl(\wt{f}(X_i) - f_0(X_i)\Bigr)^2|X_i\Bigr]\biggr]\geq \frac{1}{16rn},
\end{equation*}
whereas if $h' > 4r$  then
\begin{align*}
\Ebb\Bigl[\|\wt{f} - f_0^{(n)}\|_{n}^2\Bigr] & = \frac{1}{n}\sum_{i = 1}^{n} \Ebb\biggl[\Ebb\Bigl[\Bigl(\wt{f}(X_i) - f_0(X_i)\Bigr)^2|X_i\Bigr]\biggr] \\
& \geq \frac{1}{4nh'} + \frac{\theta^2}{8}\frac{(n - 1)}{n} P^{(n)}\Bigl([(1 - h'/2)_{+}, 1/2 - r]\Bigr) \\
&\geq \frac{1}{4nh'} + \frac{\theta^2h'}{64}.
\end{align*}
In the latter case, setting the derivative equal to $0$ shows that the right hand side is always at least $\theta/\sqrt{64 n}$, and taking the minimum over the two cases then yields~\eqref{eqn:cluster_assumption_ks_lb}.

\subsubsection{Proof of~\eqref{eqn:cluster_assumption_sp_lb}}
We begin by decomposing the risk into conditional bias and variance terms. Let $\mathbb{E}_n = \Ebb[\cdot|X_1,\ldots,X_n]$ denote expectation conditional on the design points $X_1,\ldots,X_n$.  Then by the law of iterated expectation, and the fact that $\Ebb_n[w] = 0$,
\begin{equation*}
\Ebb\bigl[\|\wt{f}_{\mathrm{SP}} - f_0\|_n^2\bigr] = \Ebb\bigl[\|\Ebb_n\wt{f}_{\mathrm{SP}} -  f_0\|_n^2\bigr] + \Ebb\bigl[\|\wt{f}_{\mathrm{SP}} - \Ebb_n\wt{f}_{\mathrm{SP}}\|_n^2\bigr]. 
\end{equation*}
We separately lower bound the expected conditional squared bias and variance terms. To anticipate what is to come: we will show that the expected conditional variance is equal to $K/n$; on the other hand we will show that the expected conditional squared bias is lower bounded,
\begin{equation}
\label{pf:cluster_assumption_sp_lb_1}
\Ebb\bigl[\|\Ebb_n\wt{f}_{\mathrm{SP}} -  f_0\|_n^2\bigr] = \frac{K}{n}\quad\textrm{and}\quad\Ebb\bigl[\|\Ebb_n\wt{f}_{\mathrm{SP}} -  f_0\|_n^2\bigr] \geq \frac{\theta^2}{2601\pi^2 K^{3}},
\end{equation}
with the lower bound holding so long as $K \leq \min\{1/(16r),n/(8 \log(8 n)), (\sqrt{160} \pi/r)^{2/3}\}$. If $K$ is larger than this, then the expected conditional variance is lower bounded,
\begin{equation}
\label{pf:cluster_assumption_sp_lb_2}
\Ebb\bigl[\|\wt{f}_{\mathrm{SP}} - \Ebb_n\wt{f}_{\mathrm{SP}}\|_n^2\bigr] \geq \min\biggl\{\frac{1}{16rn}, \frac{1}{8 \log(8n)}, \frac{(\sqrt{160} \pi)^{2/3}}{r^{2/3}n} \biggr\}
\end{equation}
Otherwise~\eqref{pf:cluster_assumption_sp_lb_1} implies that the in-sample risk is always at least 
\begin{equation*}
\Ebb\bigl[\|\wt{f}_{\mathrm{SP}} - f_0\|_n^2\bigr] \geq \frac{\theta^2}{2601\pi^2 K^{3}} + \frac{K}{n} \geq 2 \frac{\theta^{1/2}}{n^{3/4}} \frac{1}{(2601 \pi^2)^{1/4}}.
\end{equation*}
Along with~\eqref{pf:cluster_assumption_sp_lb_2}, this implies the claim. It remains to show the bounds on conditional bias and variance.

\paragraph{Conditional variance.}
The expected conditional variance is exactly equal to $K/n$, a standard fact that is verified by the following calculations: first,
\begin{align*}
\|\wt{f}_{\mathrm{SP}} - \Ebb_n\wt{f}_{\mathrm{SP}}\|_n^2 = \|\Phi(\Phi^{\top}\Phi)^{-1}\Phi^{\top}w\|_n^2 = \frac{1}{n}w^{\top}\Phi(\Phi^{\top}\Phi)^{-1}\Phi^{\top}w;
\end{align*}
thus standard properties of the Gaussian distribution and the trace trick imply
\begin{equation*}
\Ebb_n\bigl[\|\wt{f}_{\mathrm{SP}} - \Ebb_n\wt{f}_{\mathrm{SP}}\|_n^2 \bigr] = \frac{1}{n} \mathrm{tr} (\Phi(\Phi^{\top}\Phi)^{-1}\Phi^{\top}) = \frac{K}{n};
\end{equation*}
and finally by the law of iterated expectation and the independence of the noise $(w_1,\ldots,w_n)$ and the design points $X_1,\ldots,X_n$,
\begin{equation*}
\Ebb\Bigl[\Ebb_n\bigl[\|\wt{f}_{\mathrm{SP}} - \Ebb_n\wt{f}_{\mathrm{SP}}\|_n^2\Bigr] = K/n.
\end{equation*}

\paragraph{Conditional bias.}
It takes more work to lower bound the conditional bias. We will first upper bound the Lipschitz constant of $\Ebb_n\wt{f}_{\mathrm{SP}}$ in terms of the empirical norm $\|\Ebb_n\wt{f}_{\mathrm{SP}}\|_n$. Then we will use this upper bound to argue that either $\Ebb_n\wt{f}_{\mathrm{SP}}$ has empirical norm much larger than that of $f_0$, or $\Ebb_n\wt{f}_{\mathrm{SP}}$ is a smooth function, in the sense of having a small Lipschitz constant. In the former case, the triangle inequality will then imply that $\|\Ebb_n\wt{f}_{\mathrm{SP}} -  f_0\|_n$ must be large. In the latter case, the smoothness of $\Ebb_n\wt{f}_{\mathrm{SP}}$ will imply that $\Ebb_n\wt{f}_{\mathrm{SP}}$ must be far from $f_0$ at many points $X_i$ close to $x = 1/2$.

The following Lemma gives our upper bound on the Lipschitz constant of $\|\Ebb_n\wt{f}_{\mathrm{SP}}\|_n$. Here we treat $\Ebb_n\wt{f}_{\mathrm{SP}} = \sum_{k = 1}^{K} \wt{\beta}_k \phi_k$ as a function defined at all $x \in [0,1]$ by extending it in the canonical way. As a function over $[0,1]$, clearly $\Ebb_n\wt{f}_{\mathrm{SP}} \in C^\infty([0,1])$. Let $\Sigma \in \Reals^{K \times K}$ be the covariance matrix of $(\phi_1,\ldots,\phi_K)$, i.e. the matrix with entries $\Sigma_{k\ell} = \dotp{\phi_k,\phi_{\ell}}_{P^{(n)}}$. Let $\wh{\Sigma} := (\Phi^{\top}\Phi)/n$ be the empirical covariance matrix. Let $I_K \in \Reals^{K \times K}$ be the identity matrix.
\begin{lemma}[Lipschitz regularity of $\Ebb_n\wt{f}_{\mathrm{SP}}$.]
	\label{lem:lipschitz_regularity_sp}
	Let $\wt{f}_n = \Ebb_n\wt{f}_{\mathrm{SP}}$. Then 
	\begin{equation}
	\label{eqn:lipschitz_regularity_sp1}
	\|\wt{f}_n\|_{C^1(\mc{X})}^2 \leq \pi^2 \frac{K^{3} \cdot \|\Sigma^{1/2} \wh{\Sigma}^{-1} \Sigma^{1/2}\|_{\mathrm{op}}}{(1 - \|I_K - \Sigma\|_F)} \cdot \|\wt{f}_n\|_n^2.
	\end{equation}
	Moreover, suppose $K \leq 1/(16r)$ and $r \leq (1 - 2^{-1/2})/2$.
	\begin{itemize}
		\item {\bf(Matrix perturbation)} Then
		\begin{equation}
		\label{eqn:lipschitz_regularity_matrix_perturbation}
		\|\Sigma - I_K\|_F \leq \frac{1}{2}.
		\end{equation}
		\item {\bf(Matrix concentration, cf. \cite{hsu2012})} If additionally $n \geq 8K \log(K/\delta)$ for some $\delta \in (0,1/2)$, then with probability at least $1 - 2\delta$,
		\begin{equation}
		\label{eqn:lipschitz_regularity_matrix_concentration}
		\|\Sigma^{1/2} \wh{\Sigma}^{-1} \Sigma^{1/2}\|_{\mathrm{op}} \leq 5.
		\end{equation}
	\end{itemize}
	Therefore, if $K \leq \min\{1/(16r),n/(8 \log(K/\delta))\}$, then with probability at least $1 - 2\delta$,
	\begin{equation}
	\label{eqn:lipschitz_regularity_sp2}
	\|\wt{f}_n\|_{C^1(\mc{X})}^2 \leq 10 \pi^2 K^3 \|\wt{f}_n\|_n^2.
	\end{equation}
\end{lemma}
We defer the proof of Lemma~\ref{lem:lipschitz_regularity_sp} until after we complete the proof of~\eqref{eqn:cluster_assumption_sp_lb}.

Now, if $\|\wt{f}_n\|_n^2 \geq \frac{3}{2}\|f_0\|_n^2$, then by the triangle inequality
\begin{equation*}
\|\wt{f}_n - f_0\|_n \geq \|\wt{f}_n\|_n - \|f_0\|_n \geq \sqrt{\frac{3}{2}} \cdot \|f_0\|_n  = \sqrt{\frac{3}{2}} \cdot \theta.
\end{equation*}
Otherwise $\|\wt{f}_n\|_n^2 \geq \frac{3}{2}\|f_0\|_n^2$. In this case, we show that $|\wt{f}_n(X_i) - f_0(X_i)|$ must be large (on the order of $\theta$) for many points $X_i$ which are close to $x = 1/2$. Let us suppose without loss of generality that $\wt{f}_n(1/2) \leq \theta/2$ and consider points $X_i \in Q_1$ close to $x = 1/2$; otherwise if $\wt{f}_n(1/2) > \theta/2$ we could obtain the exact same bound by considering $X_i \in Q_2$. For each point $X_i \in Q_1$, by Lemma~\ref{lem:lipschitz_regularity_sp} we have that with probability at least $1 - 2\delta$,
\begin{equation*}
|\wt{f}_n(X_i) - \wt{f}_n(1/2)| \leq CK^{3/2} \|\wt{f}_n\|_n \cdot |X_i - 1/2| \leq \sqrt{10}\pi K^{3/2} \theta \cdot |X_i - 1/2|.
\end{equation*}
Since $\wt{f}_n(1/2) \leq \theta/2$ and $f_0(X_i) = \theta/2$ for all $X_i \in Q_1$ it follows that
\begin{equation*}
|\wt{f}_n(X_i) - f_0(X_i)| \geq \theta - \sqrt{10}\pi K^{3/2} \theta \cdot |X_i - 1/2|,
\end{equation*}
and consequently
\begin{equation*}
|\wt{f}_n(X_i) - f_0(X_i)| \geq \theta/2,\quad\textrm{for any $X_i \in Q_1$ such that $|X_i - 1/2| \leq 1/(\sqrt{40}\pi K^{3/2})$.}
\end{equation*}
This yields a lower bound on $\|\wt{f}_n - f_0\|_n$; letting $Q_K := \Bigl[\frac{1}{2} - \frac{1}{\sqrt{40}\pi K^{3/2}},\frac{1}{2} - r\Bigr]$, we have that 
\begin{equation*}
\|\wt{f}_n - f_0\|_n \geq \frac{\theta}{2} \cdot P_n\Bigl(Q_k\Bigr).
\end{equation*}
Then as long as $K^{-3/2} \geq \sqrt{160} \pi r$, from the multiplicative form of Hoeffding's inequality (Lemma~\ref{lem:hoeffding_smoothing})
\begin{equation*}
P^{(n)}(Q_K) \geq \frac{1}{\sqrt{160} \pi K^{3/2}} \geq 2r \Longrightarrow \Pbb\biggl(P_n(Q_K) \geq \frac{1}{\sqrt{640} \pi K^{3/2}}\biggr) \geq 1 - \exp(-nr/4) \geq 1 - \frac{4}{n^2}.
\end{equation*}
Putting the pieces together, we conclude that if $K \leq \min\{1/(8r),n/(8 \log(K/\delta)), (\sqrt{160} \pi/r)^{2/3}\}$, then
\begin{equation*}
\|\wt{f}_n - f_0\|_n \geq \frac{\theta}{51 \pi K^{3/2}},
\end{equation*}
with probability at least $1 - 2\delta - {4}{n^2}$. Taking $\delta = 1/8$ then implies the claim.

\paragraph{Proof of Lemma~\ref{lem:lipschitz_regularity_sp}.}

\underline{\it{Proof of~\eqref{eqn:lipschitz_regularity_sp1}}.}
Recall that $\wt{f}_n = \sum_{k = 1}^{K} \wt{\beta}_k \phi_k$. Exchanging sum with derivative, we have that
\begin{equation*}
\frac{d}{dx} \wt{f}_n(x) = -\pi \sum_{k = 1}^{K}\wt{\beta}_k k \sin(k\pi x).
\end{equation*}
Thus taking absolute value and applying the Cauchy-Schwarz inequality gives
\begin{equation*}
|\wt{f}_n'(x)|^2 \leq \pi^2 K^2 \sum_{k = 1}^{K} \bigl(\sin(k\pi x)\bigr)^2 \|\beta\|_2^2 \leq \pi^2 K^3 \|\beta\|_2^2.
\end{equation*}
On the other hand, we can also relate the empirical norm $\|\wt{f}_n\|_n^2$ to the $\ell^2$ norm of $\beta$. Specifically,
\begin{equation*}
\|\wt{f}_n\|_n^2 = \beta^{\top} \wh{\Sigma} \beta \geq \frac{\|\beta\|_2^2}{\|\wh{\Sigma}^{-1}\|_{\mathrm{op}}} \geq \frac{\|\beta\|_2^2}{\|\Sigma^{-1}\|_{\mathrm{op}} \cdot \|\Sigma^{1/2} \wh{\Sigma}^{-1} \Sigma^{1/2}\|_{\mathrm{op}}} = \frac{\|\beta\|_2^2 \|\Sigma\|_{\mathrm{op}} }{\|\Sigma^{1/2} \wh{\Sigma}^{-1} \Sigma^{1/2}\|_{\mathrm{op}}}
\end{equation*}
Rearranging, we see that 
\begin{equation*}
\sup_{x \in [0,1]} |\wt{f}_n'(x)|^2 \leq \frac{\pi^2 K^3}{\|\Sigma\|_{\mathrm{op}}} \|\Sigma^{1/2} \wh{\Sigma}^{-1} \Sigma^{1/2}\|_{\mathrm{op}} \leq \frac{\pi^2 K^3}{1 - \|I_K - \Sigma\|_F} \|\Sigma^{1/2} \wh{\Sigma}^{-1} \Sigma^{1/2}\|_{\mathrm{op}}
\end{equation*}
with the latter inequality following since $\|\Sigma\|_{\mathrm{op}} \geq \|I_K\|_{\mathrm{op}} - \|I_K - \Sigma\|_{\mathrm{op}} \geq 1 - \|I_K - \Sigma\|_{F}$. 

\underline{\it{Proof of~\eqref{eqn:lipschitz_regularity_matrix_perturbation}}.}
We will show that for all $1 \leq k < \ell \leq K$, 
\begin{equation}
\label{pf:lipshcitz_regularity_matrix_perturbation_1}
(1 - \dotp{\phi_k}{\phi_k}_{P^{(n)}})^2 \leq 32r^2, \quad\textrm{and}\quad|\dotp{\phi_k}{\phi_{\ell}}_{P^{(n)}}| \leq 64r^2.
\end{equation}
This implies $\|I - \Sigma\|_F^2 \leq 32K^2r^2$, so that $\|I - \Sigma\|_F \leq 1/2$ so long as $K \leq 1/(16r)$.

The proof of~\eqref{pf:lipshcitz_regularity_matrix_perturbation_1} follows from computing some standard integrals. We separate the computation based on whether $k = 1$ or $k > 1$. \newline

\underline{\it{Case 1: $k = 1$.}}
When $k = 1$, $\dotp{\phi_1}{\phi_1}_{P^{(n)}} = 1$ and $(1 - \dotp{\phi_1}{\phi_1}_{P^{(n)}})^2 = 0$. Additionally, by symbolic integration we find that
\begin{equation*}
\dotp{\phi_k}{\phi_{\ell}}_{P^{(n)}} = \frac{-2\sqrt{2}}{(1 - 2r)} \cdot \frac{\cos(\ell\pi/2) \sin(\ell\pi r)}{\ell \pi},
\end{equation*}
and therefore
\begin{equation*}
\bigl[\dotp{\phi_k}{\phi_{\ell}}_{P^{(n)}}\bigr]^2 \leq \frac{8}{(1 - 2r)^2} \cdot \biggl(\frac{\sin(\ell\pi r)}{\ell \pi}\biggr)^2 \leq \frac{8}{(1 - 2r)^2} r^2 \leq 16r^2,
\end{equation*}
where in the second-to-last inequality follows because $\sin(x)/x \leq 1$, and the last inequality follows by our assumed upper bound on $r$. \newline

\underline{\it{Case 2: $k > 1$.}}
When $k > 1$,
\begin{equation*}
\dotp{\phi_k}{\phi_k}_{P^{(n)}} = 1 - \frac{2}{(1 - 2r)} \frac{\cos(k\pi) \sin(2k\pi r)}{k\pi} \Longrightarrow \bigl[1 - \dotp{\phi_k}{\phi_k}_{P^{(n)}}\bigr]^2 \leq \frac{4}{(1 - 2r)^2} \cdot \biggl(\frac{\sin(2k\pi r)}{k\pi}\biggr)^2 \leq 32r^2.
\end{equation*}
Similarly,
\begin{equation*}
\dotp{\phi_k}{\phi_\ell}_{P^{(n)}} = -\frac{4}{(1 - 2r)}\biggl[\frac{\cos((k + \ell)\pi) \sin((k + \ell)\pi r)}{(k + \ell)\pi} + \frac{\cos((k - \ell)\pi) \sin((k - \ell)\pi r)}{(k - \ell)\pi}\biggr]
\end{equation*}
and therefore
\begin{equation*}
\bigl[\dotp{\phi_k}{\phi_\ell}_{P^{(n)}}\bigr]^2 \leq \frac{16}{(1 - 2r)^2} \biggl(\biggl[\frac{\sin((k + \ell)\pi r)}{(k + \ell)\pi}\biggr]^2 + \biggl[\frac{\sin((k - \ell)\pi r)}{(k - \ell)\pi}\biggr]^2\biggr) \leq 64 r^2.
\end{equation*}

\underline{\it{Proof of~\eqref{eqn:lipschitz_regularity_matrix_concentration}}}
Denote $\Phi(x) = (\phi_1,\ldots,\phi_K(x)) \in \Reals^K$ for any $x \in [0,1]$. Then for any $x \in [0,1]$,
\begin{equation*}
\|\Sigma^{-1/2} \Phi(x)\| \leq \|\Sigma^{-1}\|_{\mathrm{op}}^{1/2} \|\Phi(x)\|_2 \leq \|\Sigma^{-1}\|_{\mathrm{op}}^{1/2} \sqrt{2K} \leq 2 \sqrt{K}
\end{equation*}
with the second-to-last inequality following from~\eqref{eqn:lipschitz_regularity_matrix_perturbation}, and the last inequality following since $|\phi_k(x)| \leq \sqrt{2}$ for all $k$. Thus $\|\Sigma^{-1/2} \Phi(x)\|/\sqrt{K} \leq 2$, and~\eqref{eqn:lipschitz_regularity_matrix_concentration} follows from Theorem~1 of \citet{hsu2012}.

\underline{\it{Proof of~\eqref{eqn:lipschitz_regularity_sp2}}.}
Follows immediately.

\section{Miscellaneous}
Here we give assorted helpful Lemmas used at various points in the above proofs. We also review notation and relevant facts regarding Taylor expansion.

\subsection{Concentration Inequalities}
Lemma~\ref{lem:chi_square_bound} controls the deviation of a chi-squared random variable. It is from~\cite{laurent00}.
\begin{lemma}
	\label{lem:chi_square_bound}
	Let $\xi_1,\ldots,\xi_N$ be independent $N(0,1)$ random variables, and let $U := \sum_{k = 1}^{N} a_k(\xi_k^2 - 1)$.  Then for any $t > 0$,
	\begin{equation*}
	\Pbb\Bigl[U \geq 2 \|a\|_2 \sqrt{t} + 2 \|a\|_{\infty}t\Bigr] \leq \exp(-t).
	\end{equation*}
	In particular if $a_k = 1$ for each $k = 1,\ldots,N$, then
	\begin{equation*}
	\Pbb\Bigl[U\geq 2\sqrt{N t} + 2t\Bigr] \leq \exp(-t).
	\end{equation*}
\end{lemma}

Lemma~\ref{lem:one_sided_bernstein} is an immediate consequence of the one-sided Bernstein's inequality (14.23) in \cite{wainwright2019}.
\begin{lemma}[One-sided Bernstein's inequality]
	\label{lem:one_sided_bernstein}
	Let $X, X_1,\ldots,X_n \sim P$, and $f$ satisfy $\Ebb[f^4(X)] < \infty$. Then
	\begin{equation*}
	\|f\|_n^2 \geq \frac{1}{2}\|f\|_P^2,
	\end{equation*}
	with probability at least $1 - \exp\bigl(-n/8 \cdot \|f\|_P^4 /\Ebb[f^4(X)]\bigr)$.
\end{lemma}

Lemma~\ref{lem:hoeffding_smoothing} is a multiplicative form of Hoeffding's inequality.
\begin{lemma}[Hoeffding's Inequality, multiplicative form]
	\label{lem:hoeffding_smoothing}
	Suppose $Z_i$ are independent random variables, which satisfy $Z_i \in [0,B]$ for $i = 1,\ldots,n$. For any $0 < \delta < 1$, it holds that
	\begin{equation*}
	\Pbb\biggl(\Bigl|S_n - \mu\Bigr| \geq \delta \mu\biggr) \leq 2\exp\biggl(-\frac{\delta^2\mu}{3B^2}\biggr).
	\end{equation*}
\end{lemma}

The following Lemma gives a ``balls-in-bins'' result. More precisely, it gives a lower bound on the probability that every bin
\begin{equation*}
Q_{i1} = [i/m,(i + 1)/m] \cdot (1/2 - r), \quad Q_{i2} = 1/2 + [i/m,(i + 1)/m] \cdot (1/2 - r).
\end{equation*}
will contain at least one ball.
\begin{lemma}
	\label{lem:balls_in_bins}
	Suppose $(X_1,Y_1),\ldots(X_n,Y_n)$ are sampled according to~\eqref{def:model_cluster_assumption}, and suppose $r \leq 1/4$. We have that 
	\begin{equation}
	\label{eqn:balls_in_bins}
	\Pbb\Bigl(\sharp\{Q_{ij} \cup {\bf X}\} > 0 ~~\textrm{for all $i = 1,\ldots,m - 1$ and $j = 1,2$} \Bigr) \geq 1 - 2m\exp\{-n/2m\}.
	\end{equation}
\end{lemma}
\paragraph{Proof (of Lemma~\ref{lem:balls_in_bins}).}
For each $Q_{ij}$, we have that $P(Q_{ij}) = (1/2 - r)/m \geq 1/(2m)$. Therefore
\begin{equation*}
\Pbb\bigl(\sharp\{Q_{ij} \cup {\bf X}\} = 0\bigr) = (1 - 1/(2m))^{n} \leq \exp\{-n/2m\}.
\end{equation*}
By a union bound,
\begin{equation*}
\Pbb\Bigl(\sharp\{Q_{ij} \cup {\bf X}\} = 0 ~~\textrm{for any $i = 1,\ldots,m - 1$ and $j = 1,2$} \Bigr) \leq 2m\exp\{-n/2m\}. \qed
\end{equation*}
Let $\varepsilon = 2/m$. Note that by construction,~\eqref{eqn:balls_in_bins} implies that any points $x$ and $x'$ in adjacent intervals $Q_{ij}$ and $Q_{i'j}$ must be connected in $G_{n,\varepsilon}$. Likewise, it implies that for $h = 1/m$ the degree $d_{n,h}(x) > 0$ for every $x \in Q_1 \cup Q_2$.

\subsection{Taylor expansion}
\label{subsec:taylor_expansion}
We begin with some notation that allows us to concisely derivatives. For a given $z \in \Rd$ and $s$-times differentiable function $f: \mc{X} \to \Reals$, we denote $\bigl(d_x^sf\bigr)(z) := \sum_{|\alpha| = s} D^{\alpha}f(x) z^{\alpha}$. We also write $d^sf := \sum_{|\alpha| = j} D^{\alpha}f$. We point out that in the first-order case $d_x^1f$ is the differential of $f$ at $x \in \mc{X}$, while $d^1f$ is the divergence of $f$.

Let $u$ be a function which is $s$ times continuously differentiable at all $x \in \mc{X}$, for $k \in \mathbb{N}\setminus\{0\}$. Suppose that for some $h > 0$, $x \in \mc{X}_{h}$ and $x' \in B(x,h)$. We write the order-$s$ Taylor expansion of $u(x')$ around $x' = x$ as
\begin{equation*}
u(x') = u(x) + \sum_{j = 1}^{s - 1} \frac{1}{j!}\bigl(d_x^{j}u\bigr)(x' - x) + r_{x'}^{s}(x;u)
\end{equation*}
For notational convenience we have adopted the convention that $\sum_{j = 1}^{0} a_j = 0$.  Thus $\bigl(d_x^{j}f\bigr)(z)$ is a degree-$j$ polynomial---and so a $j$-homogeneous function---in $z$, meaning for any $t \in \Reals$,
\begin{equation*}
\bigl(d_x^{j}f\bigr)(tz) = t^{j} \cdot \bigl(d_x^{j}f\bigr)(z).
\end{equation*}
The remainder term $r_{x'}$ is given by
\begin{equation*}
r_{x'}^s(x;f) = \frac{1}{(j - 1)!} \int_{0}^{1}(1 - t)^{j - 1} \bigl(d_{x + t(x' - x)}^{s}f\bigr)(x' - x) \,dt,
\end{equation*}
where we point out that the integral makes sense because $x + t(x' - x) \in B(x,h) \subseteq \mc{X}$. We now give estimates on the remainder term in both sup-norm and $L^2(\mc{X}_{h})$ norm, each of which hold for any $z \in B(0,1)$. In sup-norm, we have that 
\begin{equation*}
\sup_{x \in \mc{X}_{h}}|r_{x + hz}^j(x;f)| \leq C h^j \|f\|_{C^j(\mc{X})},
\end{equation*}
whereas in $L^2(\mc{X}_{h})$ norm we have,
\begin{equation}
\label{eqn:sobolev_remainder_term}
\int_{\mc{X}_{h}} \bigl|r_{x + thz}^j(x;f)\bigr|^2 \,dx \leq h^{2j} \int_{\mc{X}_{h}} \int_{0}^{1} |d_{x + thz}^jf(z)|^2 \,dt \,dx \leq h^{2j} \|d^jf\|_{L^2(\mc{X})}^2.
\end{equation}
In the last inequality 

Finally, we recall some facts regarding the interaction between smoothing kernels and polynomials.  Let $q_j(z)$ be an arbitrary degree-$j$ (multivariate) polynomial. If $\eta$ is a radially symmetric kernel and $j$ is odd, then by symmetry it follows that
\begin{equation*}
\int_{B(0,1)} q_j(z) \eta(\|z\|) \,dz = 0.
\end{equation*}
On the other hand, if $\psi$ is an order-$s$ kernel for some $s > j$, then by converting to polar coordinates we can verify that
\begin{equation*}
\int_{B(0,1)} q_j(z) \eta(\|z\|) \,dz = 0.
\end{equation*}

\section{Computational considerations}
\label{sec:computational_considerations}
Recall that when $s = 1$, we have shown that PCR-LE is optimal when $\varepsilon \asymp (\log n/n)^{1/d}$ is (up to a constant) as small as possible while still ensuring the graph $G$ is connected. On the other hand, when $s > 1$, we can show PCR-LE is optimal only when $\varepsilon = \omega(n^{-c})$ for some $c < 1/d$. For such a choice of $\varepsilon$, the average degree in $G$ will grow polynomially in $n$ as $n \to \infty$, and computing eigenvectors of the Laplacian of a graph will be more computationally intensive than if the graph were sparse. 
In this dense-graph setting, we now discuss a procedure to more efficiently compute an approximation to the PCR-LE estimate: \emph{edge sparsification}.

By now there exist various methods see (e.g., the seminal papers of \citet{spielman2011,spielman2013,spielman2014}, or the overview by \citet{vishnoi2012} and references therein) to efficiently remove many edges from the graph $G$ while only slightly perturbing the spectrum of the Laplacian. Specifically such algorithms take as input a parameter $\sigma \geq 1$, and return a sparser graph $\wc{G}$, $E(\wc{G}) \subseteq E(G)$, with a Laplacian $\wc{L}_{n,\varepsilon}$ satisfying
\begin{equation*}
\frac{1}{\sigma} \cdot u^{\top} \wc{L}_{n,\varepsilon} u \leq u^{\top} L_{n,\varepsilon} u \leq \sigma \cdot u^{\top} \wc{L}_{n,\varepsilon}u \quad \textrm{for all $u \in \Reals^n$.}
\end{equation*}
Let $\wc{f}$ be the PCR-LE estimator computed using the eigenvectors of the sparsified graph Laplacian $\wc{L}_{n,\varepsilon}$ . Because $\wc{G}$ is sparser than $G$, it can be (much) faster to compute the eigenvectors of $\wc{L}_{n,\varepsilon}$ than the eigenvectors of $L_{n,\varepsilon}$, and consequently much faster to compute $\wc{f}$ than $\wh{f}$. Statistically speaking, letting $\wc{\lambda}_k$ be the $k$th eigenvalue of $\wc{L}_{n,\varepsilon}$, we have that conditional on $\{X_1,\ldots,X_n\}$,
\begin{equation*}
\|\wc{f} - f_0\|_n^2 \leq \frac{\dotp{\wc{L}_{n,\varepsilon}^s f_0}{f_0}_n}{\wc{\lambda}_{K + 1}^s} + \frac{5K}{n} \leq \sigma^{2s} \frac{\dotp{\wc{L}_{n,\varepsilon}^s f_0}{f_0}_n}{\wc{\lambda}_{K + 1}^s} + \frac{5K}{n},
\end{equation*}
with probability at least $1 - \exp(-K)$. Consequently $\|\wt{f} - f_0\|_n^2$is at most $\sigma^{2s} \cdot \|\wh{f} - f_0\|_n^2$, and for any choice of $\sigma$ that is constant in $n$ the estimator $\wc{f}$ will also be rate-optimal. 

In fact the aforementioned edge sparsification algorithms are overkill for our needs. For one thing, they are designed to work when $\sigma$ is very close to $1$, whereas in order for $\wc{f}$ to be rate-optimal, setting $\sigma$ to be any constant greater than $1$, say $\sigma = 2$, is sufficient. Additionally, edge sparsification algorithms are traditionally designed to work in the worst-case, where no assumptions are made on the structure of the graph $G$. But the geometric graphs we consider in this paper exhibit a special structure, in which very roughly speaking no single edge is a bottleneck. As pointed out by~\citet{sadhanala16b}, in this special case there are far simpler and faster methods for sparsification, which at least empirically seem to do the job.